\documentclass[smallextended]{svjour3}       
\smartqed  
\usepackage{graphicx}
\usepackage{amsfonts}
\usepackage{algorithm}
\usepackage{amsmath,booktabs,ctable,threeparttable}
\usepackage{amssymb,amsfonts,boxedminipage}
\usepackage{mathptmx}      
\usepackage{latexsym}

\journalname{myjournal}

\begin{document}

\title{The regularization theory of the Krylov iterative solvers LSQR and CGLS
for linear discrete ill-posed problems, part I: the simple singular
value case\thanks{This work was supported in part by
the National Science Foundation of China (No. 11371219)}}

\titlerunning{Regularization Theory of LSQR and CGLS} 

\author{Zhongxiao Jia}
\institute{Zhongxiao Jia \at Department of Mathematical Sciences, Tsinghua
University, 100084 Beijing, China. \\
              \email{jiazx@tsinghua.edu.cn}}
\date{Received: date / Accepted: date}
\maketitle

\begin{abstract}
For the large-scale linear discrete ill-posed problem $\min\|Ax-b\|$ or $Ax=b$
with $b$ contaminated by a white noise, the Lanczos bidiagonalization based
LSQR method and its mathematically equivalent Conjugate Gradient (CG) method for
$A^TAx=A^Tb$ are most commonly used. They have intrinsic regularizing effects,
where the number $k$ of iterations plays the role
of regularization parameter. However, there has been no
answer to the long-standing fundamental concern by Bj\"{o}rck and
Eld\'{e}n in 1979: {\em for which kinds of
problems LSQR and CGLS can find best possible regularized solutions}?
Here a best possible regularized solution means that
it is at least as accurate as the best regularized
solution obtained by the truncated singular value decomposition (TSVD) method
or standard-form Tikhonov regularization. In this paper,
assuming that the singular values of $A$ are simple, we analyze the
regularization of LSQR for severely, moderately and mildly ill-posed problems.
We establish
accurate estimates for the 2-norm distance between the underlying $k$-dimensional
Krylov subspace and the $k$-dimensional dominant right singular subspace
of $A$. For the first two kinds of
problems, we then prove that LSQR finds a best possible regularized solution at
semi-convergence occurring at iteration $k_0$ and that, for
$k=1,2,\ldots,k_0$, (i) the $k$-step Lanczos
bidiagonalization always generates a near best rank $k$ approximation to $A$;
(ii) the $k$ Ritz values always approximate the first $k$ large singular values
in natural order; (iii) the $k$-step LSQR
always captures the $k$ dominant SVD components of $A$.
For the third kind of problem, we prove that LSQR generally cannot find a best
possible regularized solution. We derive estimates for
the entries of the bidiagonal matrices generated by
Lanczos bidiagonalization, which can be practically exploited to identify
if LSQR finds a best possible regularized solution at semi-convergence.
Numerical experiments confirm our theory.
\end{abstract}


\section{Introduction and Preliminaries}\label{intro}

Consider the linear discrete ill-posed problem
\begin{equation}
  \min\limits_{x\in \mathbb{R}^{n}}\|Ax-b\| \mbox{\,\ or \ $Ax=b$,}
  \ \ \ A\in \mathbb{R}^{m\times n}, \label{eq1}
  \ b\in \mathbb{R}^{m},
\end{equation}
where the norm $\|\cdot\|$ is the 2-norm of a vector or matrix, and
$A$ is extremely ill conditioned with its singular values decaying
to zero without a noticeable gap. \eqref{eq1} mainly arises
from the discretization of the first kind Fredholm integral equation
\begin{equation}\label{eq2}
Kx=(Kx)(t)=\int_{\Omega} k(s,t)x(t)dt=g(s)=g,\ s\in \Omega
\subset\mathbb{R}^q,
\end{equation}
where the kernel $k(s,t)\in L^2({\Omega\times\Omega})$ and
$g(s)$ are known functions, while $x(t)$ is the
unknown function to be sought. If $k(s,t)$ is non-degenerate
and $g(s)$ satisfies the Picard condition, there exists the unique squares
integrable solution
$x(t)$; see \cite{engl00,hansen98,hansen10,kirsch,mueller}. Here for brevity
we assume that $s$ and $t$ belong to the same set $\Omega\subset
\mathbb{R}^q$ with $q\geq 1$.
Applications include image deblurring, signal processing, geophysics,
computerized tomography, heat propagation, biomedical and optical imaging,
groundwater modeling, and many others; see, e.g.,
\cite{aster,engl93,engl00,hansen10,ito15,kaipio,kern,kirsch,mueller,natterer,vogel02}.
The theory and numerical treatments of integral
equations can be found in \cite{kirsch,kythe}.
The right-hand side $b=\hat{b}+e$ is noisy and assumed to be
contaminated by a white noise $e$, caused by measurement, modeling
or discretization errors, where $\hat{b}$
is noise-free and $\|e\|<\|\hat{b}\|$.
Because of the presence of noise $e$ and the extreme
ill-conditioning of $A$, the naive
solution $x_{naive}=A^{\dagger}b$ of \eqref{eq1} bears no relation to
the true solution $x_{true}=A^{\dagger}\hat{b}$, where
$\dagger$ denotes the Moore-Penrose inverse of a matrix.
Therefore, one has to use regularization to extract a
best possible approximation to $x_{true}$.

In principle, regularizing an ill-posed problem is to replace it by a well-posed one,
such that the error is compensated by the gain in stability. In other words,
regularization is to compromise the error and stability as best as possible.
For a white noise $e$,
throughout the paper, we always assume that $\hat{b}$ satisfies the discrete
Picard condition $\|A^{\dagger}\hat{b}\|\leq C$ with some constant $C$ for $n$
arbitrarily large \cite{aster,gazzola15,hansen90,hansen90b,hansen98,hansen10,kern}.
It is an analog of the Picard condition in the  finite dimensional case;
see, e.g., \cite{hansen90}, \cite[p.9]{hansen98},
\cite[p.12]{hansen10} and \cite[p.63]{kern}.
The two dominating regularization approaches are
to solve the following two essentially equivalent problems
\begin{equation}\label{posed}
\min\limits_{x\in \mathbb{R}^{n}}\|Lx\| \ \ \mbox{subject to}\ \ \|Ax-b\|=\min
\end{equation}
and general-form Tikhonov regularization (cf. \cite{phillips,tikhonov63,tikhonov77})
\begin{equation}\label{tikhonov}
  \min\limits_{x\in \mathbb{R}^{n}}\{\|Ax-b\|^2+\lambda^2\|Lx\|^2\}
\end{equation}
with $\lambda>0$ the regularization parameter
for regularized solutions \cite{hansen98,hansen10}.
A suitable choice of the matrix $L$ is based on a-prior information on $x_{true}$,
and typically $L$ is either the identity matrix, a diagonal weighting matrix, or
a $p\times n$ discrete approximation of a first or second order derivative operator.
Particularly, if $L=I$, the identity matrix,
\eqref{tikhonov} is standard-form Tikhonov regularization.

The case
$L=I$ is of most common interests and our concern in this paper. From now on,
we always assume $L=I$, for which the solutions to \eqref{eq1},
\eqref{posed} and \eqref{tikhonov} can be
fully analyzed by the singular value decomposition (SVD) of $A$. Let
\begin{equation}\label{eqsvd}
  A=U\left(\begin{array}{c} \Sigma \\ \mathbf{0} \end{array}\right) V^{T}
\end{equation}
be the SVD of $A$,
where $U = (u_1,u_2,\ldots,u_m)\in\mathbb{R}^{m\times m}$ and
$V = (v_1,v_2,\ldots,v_n)\in\mathbb{R}^{n\times n}$ are orthogonal,
$\Sigma = {\rm diag} (\sigma_1,\sigma_2,\ldots,\sigma_n)\in\mathbb{R}^{n\times n}$
with the singular values
$\sigma_1>\sigma_2 >\cdots >\sigma_n>0$ assumed to be simple
throughout the paper, and the superscript $T$
denotes the transpose of a matrix or vector. Then
\begin{equation}\label{eq4}
  x_{naive}=\sum\limits_{i=1}^{n}\frac{u_i^{T}b}{\sigma_i}v_i =
  \sum\limits_{i=1}^{n}\frac{u_i^{T}\hat{b}}{\sigma_i}v_i +
  \sum\limits_{i=1}^{n}\frac{u_i^{T}e}{\sigma_i}v_i
  =x_{true}+\sum\limits_{i=1}^{n}\frac{u_i^{T}e}{\sigma_i}v_i
\end{equation}
with $\|x_{true}\|=\|A^{\dagger}\hat{b}\|=
\left(\sum_{i=1}^n\frac{|u_i^T\hat{b}|^2}{\sigma_i^2}\right)^{1/2}$.

The discrete Picard condition means that, on average, the Fourier coefficients
$|u_i^{T}\hat{b}|$ decay faster than $\sigma_i$ and enables
regularization to compute useful approximations to $x_{true}$, which
results in the following popular model that is used throughout Hansen's books
\cite{hansen98,hansen10} and the current paper:
\begin{equation}\label{picard}
  | u_i^T \hat b|=\sigma_i^{1+\beta},\ \ \beta>0,\ i=1,2,\ldots,n,
\end{equation}
where $\beta$ is a model parameter that controls the decay rates of
$| u_i^T \hat b|$.
Hansen \cite[p.68]{hansen10} points out, ``while this is a crude model,
it reflects the overall behavior often found in real problems."
One precise definition of the discrete Picard condition is
$| u_i^T \hat b|=\tau_i\sigma_i^{1+\zeta_i}$ with certain
constants $\tau_i\geq 0,\ \zeta_i>0,\ i=1,2,\ldots,n$. We remark that once
the $\tau_i>0$ and $\zeta_i$ do not differ greatly, such discrete Picard
condition does not affect our claims, rather it
complicates derivations and forms of the results.

The white noise $e$ has a number of attractive properties which
play a critical role in the regularization analysis: Its covariance matrix
is $\eta^2 I$, the expected values ${\cal E}(\|e\|^2)=m \eta^2$ and
${\cal E}(|u_i^Te|)=\eta,\,i=1,2,\ldots,n$, and
$\|e\|\approx \sqrt{m}\eta$ and $|u_i^Te|\approx \eta,\
i=1,2,\ldots,n$; see, e.g., \cite[p.70-1]{hansen98} and \cite[p.41-2]{hansen10}.
The noise $e$ thus affects $u_i^Tb,\ i=1,2,\ldots,n,$ {\em more or less equally}.
With \eqref{picard}, relation \eqref{eq4} shows that for large singular values
$|{u_i^{T}\hat{b}}|/{\sigma_i}$ is dominant relative to
$|u_i^{T}e|/{\sigma_i}$. Once
$| u_i^T \hat b| \leq | u_i^T e|$ from some $i$ onwards, the small singular
values magnify $|u_i^{T}e|/{\sigma_i}$, and the noise
$e$ dominates $| u_i^T b|/\sigma_i$ and must be suppressed. The
transition point $k_0$ is such that
\begin{equation}\label{picard1}
| u_{k_0}^T b|\approx | u_{k_0}^T \hat{b}|> | u_{k_0}^T e|\approx
\eta, \ | u_{k_0+1}^T b|
\approx | u_{k_0+1}^Te|
\approx \eta;
\end{equation}
see \cite[p.42, 98]{hansen10} and a similar description \cite[p.70-1]{hansen98}.
The $\sigma_k$ are then divided into the $k_0$ large ones and the $n-k_0$ small
ones.

The truncated SVD (TSVD) method \cite{hansen90,hansen98,hansen10} deals
with \eqref{posed} by solving
\begin{equation}\label{tsvd}
\min\|x\| \ \ \mbox{subject to}\ \
 \|A_kx-b\|=\min,\ k=1,2,\ldots,n,
\end{equation}
where $A_k=U_k\Sigma_k V_k^T$
is the best rank approximation $k$ to $A$ with respect to the 2-norm
with $U_k=(u_1,\ldots,u_k)$, $V_k=(v_1,\ldots,v_k)$ and $\Sigma_k=
{\rm diag}(\sigma_1,\ldots,\sigma_k)$; it holds that
$\|A-A_k\|=\sigma_{k+1}$ (cf. \cite[p.12]{bjorck96}).
 and
$
x_{k}^{tsvd}=A_k^{\dagger}b,
$
called the TSVD solution,
solves \eqref{tsvd}. An crucial observation is that
$x_k^{tsvd}$ is the minimum-norm least squares solution to
$
  \min\limits_{x\in \mathbb{R}^{n}}\|A_kx-b\|
$
that perturbs $A$ to $A_k$ in \eqref{eq1}, and we will frequently exploit this
interpretation later.

Based on the above properties of the white noise $e$, it is known
from \cite[p.70-1]{hansen98} and
\cite[p.71,86-8,95]{hansen10} that the TSVD solutions
\begin{equation}\label{solution}
  x^{tsvd}_k=A_k^{\dagger}b=\left\{\begin{array}{ll}
  \sum\limits_{i=1}^{k}\frac{u_i^{T}b}{\sigma_i}{v_i}\approx
  \sum\limits_{i=1}^{k}\frac{u_i^{T}\hat{b}}
{\sigma_i}{v_i},\ \ \ &k\leq k_0;\\ \sum\limits_{i=1}^{k}\frac{u_i^{T}b}
{\sigma_i}{v_i}\approx
\sum\limits_{i=1}^{k_0}\frac{u_i^{T}\hat{b}}{\sigma_i}{v_i}+
\sum\limits_{i=k_0+1}^{k}\frac{u_i^{T}e}{\sigma_i}{v_i},\ \ \ &k>k_0,
\end{array}\right.
\end{equation}
and $x_{k_0}^{tsvd}$ is
the best TSVD regularized solution to \eqref{eq1}, which balances the
regularization and perturbation errors optimally and stabilizes the
residual norms $\|Ax_k^{tsvd}-b\|$ for $k$ not close to $n$
after $k>k_0$. The index $k$ plays the role of the
regularization parameter that determines how many large SVD components of
$A$ are used to compute a regularized solution $x_k^{tsvd}$ to
\eqref{eq1}.

The solution $x_{\lambda}$ of the Tikhonov regularization
has a filtered SVD expansion
\begin{equation}\label{eqfilter}
  x_{\lambda} = \sum\limits_{i=1}^{n}f_i\frac{u_i^{T}b}{\sigma_i}v_i,
\end{equation}
where the $f_i=\frac{\sigma_i^2}{\sigma_i^2+\lambda^2}$ are called filters.
The TSVD method is a special parameter filtered method, where, in $x_k^{tsvd}$,
we take $f_i=1,\ i=1,2,\ldots,k$ and $f_i=0,\ i=k+1,\ldots,n$. The error
$x_{\lambda}-x_{true}$ can be written as the sum of the regularization and
perturbation errors, and an optimal $\lambda_{opt}$ aims to balance
these two errors and make the sum of their norms minimized
\cite{hansen98,hansen10,kirsch,vogel02}. The best possible regularized
solution $x_{\lambda_{opt}}$ retains the $k_0$ dominant SVD components
and dampens the other $n-k_0$ small SVD components as much as
possible \cite{hansen98,hansen10}. Apparently, the ability to acquire {\em only}
the largest SVD components of $A$ is fundamental in solving \eqref{eq1}.

A number of parameter-choice methods have been developed for finding
$\lambda_{opt}$ or $k_0$, such as the discrepancy principle \cite{morozov},
the L-curve criterion, whose use goes back to Miller \cite{miller} and
Lawson and Hanson \cite{lawson} and is termed much later and studied in detail
in \cite{hansen92,hansen93}, and the generalized cross validation
(GCV) \cite{golub79,wahba}; see, e.g.,
\cite{bauer11,hansen98,hansen10,kern,kilmer03,kindermann,neumaier98,reichel13,vogel02}
for numerous comparisons. All parameter-choice methods aim to
make $f_i/\sigma_i$ not small for $i=1,2,\ldots,k_0$
and $f_i/\sigma_i\approx 0$ for $i=k_0+1,\ldots,n$. Each of these methods
has its own merits and disadvantages, and
no one is absolutely reliable for all ill-posed problems.
For example, some of the mentioned parameter-choice methods
may fail to find accurate approximations to $\lambda_{opt}$;
see \cite{hanke96a,vogel96} for an analysis
on the L-curve method and \cite{hansen98} for some other parameter-choice
methods. A further investigation on paramater-choice methods is not our concern
in this paper.

The TSVD method and the standard-form Tikhonov regularization
produce very similar solutions with essentially the minimum
2-norm error, i.e., the worst-case error \cite[p.13]{kirsch};
see \cite{varah79}, \cite{hansen90b}, \cite[p.109-11]{hansen98} and
\cite[Sections 4.2 and 4.4]{hansen10}. Indeed, for an underlying linear
compact equation $Kx=g$, e.g.,\eqref{eq2}, with the noisy $g$ and true solution
$x_{true}(t)$, under the source condition that its solution $x_{true}(t)
\in {\cal R}(K^*)$ or $x_{true}(t)\in {\cal R}(K^*K)$, the range of
the adjoint $K^*$ of $K$ or
that of $K^*K$, which amounts to assuming that $x_{true}(t)$ or its
derivative is squares integrable, the errors of
the best regularized solutions by the TSVD method and the Tikhonov
regularization are {\em order optimal, i.e., the same order
as the worst-case error} \cite[p.13,18,20,32-40]{kirsch},
\cite[p.90]{natterer} and \cite[p.7-12]{vogel02}. These conclusions
carries over to \eqref{eq1} \cite[p.8]{vogel02}. Therefore,
both $x_{\lambda_{opt}}$ and $x_{k_0}^{tsvd}$
are best possible solutions to \eqref{eq1} under the above
assumptions, and any of them can be taken as the reference standard
when assessing the regularizing effects of an iterative
solver. For the sake of clarity and analysis, we will take
$x_{k_0}^{tsvd}$ as the standard reference.

For \eqref{eq1} large, the TSVD method and the Tikhonov regularization
method are generally too demanding, and only iterative regularization
methods are computationally viable. A major class of methods has been
Krylov iterative solvers that project \eqref{eq1} onto a sequence of
low dimensional Krylov subspaces
and computes iterates to approximate $x_{true}$
\cite{aster,engl00,gilyazov,hanke95,hansen98,hansen10,kirsch}.
Of Krylov iterative solvers, the CGLS (or CGNR) method,
which implicitly applies the CG
method \cite{golub89,hestenes} to $A^TAx=A^Tb$,
and its mathematically equivalent LSQR algorithm \cite{paige82}
have been most commonly used. The Krylov solvers CGME
(or CGNE) \cite{bjorck96,bjorck15,craig,hanke95,hanke01} and
LSMR \cite{bjorck15,fong} are also choices, which amount to the
CG method applied to $\min\|AA^Ty-b\|$ or $AA^Ty=b$
with $x=A^Ty$ and MINRES \cite{paige75}
applied to $A^TAx=A^Tb$, respectively. These Krylov solvers have been
intensively studied and known to have general regularizing
effects \cite{aster,eicke,gilyazov,hanke95,hanke01,hansen98,hansen10,hps16,hps09}
and exhibit semi-convergence \cite[p.89]{natterer};
see also \cite[p.314]{bjorck96}, \cite[p.733]{bjorck15},
\cite[p.135]{hansen98} and \cite[p.110]{hansen10}: The iterates
converge to $x_{true}$ and their norms increase steadily,
and the residual norms decrease in an initial stage; afterwards the
noise $e$ starts to deteriorate the iterates so that they start to diverge
from $x_{true}$ and instead converge to $x_{naive}$,
while their norms increase considerably and the residual norms stabilize.
If we stop at the right time, then, in principle,
we have a regularization method, where the iteration number plays the
role of the regularization parameter.
Semi-convergence is due to the fact that the projected problem starts to
inherit the ill-conditioning of \eqref{eq1} from some iteration
onwards, and a small singular
value of the projected problem amplifies the noise considerably.

The regularizing effects of CG type methods were noticed by
Lanczos \cite{lanczos} and were rediscovered in \cite{johnsson,squire,tal}.
Based on these works and motivated by a heuristic explanation on good
numerical results with very few iterations using CGLS
in \cite{johnsson}, and realizing that such an excellent performance
can only be expected if convergence to the regular part
of the solution, i.e., $x_{k_0}^{tsvd}$, takes place before the effects of
ill-posedness show up, on page 13 of \cite{bjorck79},
Bj\"{o}rck and Eld\'{e}n in 1979 foresightedly expressed a fundamental concern
on CGLS (and LSQR): {\em More
research is needed to tell for which problems this approach will work, and
what stopping criterion to choose.} See also \cite[p.145]{hansen98}.
As remarked by Hanke and Hansen \cite{hanke93}, the paper \cite{bjorck79}
was the only extensive survey on algorithmic details until that time,
and a strict proof of the regularizing properties of conjugate gradients is
extremely difficult. An enormous effort has long been made to the study of
regularizing effects of LSQR and CGLS (cf.
\cite{firro97,gazzola16,gilyza86,hanke95,hanke01,hansen98,hansen10,hps16,hps09,huangjia,nemi,nolet,paige06,scales,vorst90}),
but hitherto there has been no definitive answer to the above long-standing
fundamental question, and the same is for CGME and LSMR.

For $A$ symmetric, MINRES and MR-II applied to $Ax=b$ directly are alternatives
and have been shown to have regularizing effects
\cite{calvetti01,hanke95,hanke96,hansen10,jensen07,kilmer99}, but MR-II seems
preferable since
the noisy $b$ is excluded in the underlying subspace \cite{huang15,jensen07}.
For $A$ nonsymmetric or multiplication with $A^{T}$ difficult to compute,
GMRES and RRGMRES are candidate methods
\cite{baglama07,calvetti02,calvetti02c,neuman12}, and the latter
may be better \cite{jensen07}. The hybrid
approaches based on the Arnoldi process have been first
proposed in \cite{calvetti00b} and studied in
\cite{calvetti01,calvetti03,lewis09,novati13}.
Gazzola and her coauthors
\cite{gazzola14,gazzola15,gazzola14b,gazzola-online}
have described a general framework of the hybrid methods
and presented various Krylov-Tikhonov methods with different parameter-choice
strategies. The regularizing effects of these methods are highly problem
dependent, and it appears that they require that the mixing of the left and
right singular vectors of $A$ be weak, that is, $V^TU$ is
close to a diagonal matrix; for more details, see,
e.g., \cite{jensen07} and \cite[p.126]{hansen10}.

The behavior of ill-posed problems critically depends on the decay rate of
$\sigma_j$. The following characterization of the degree of ill-posedness
of \eqref{eq1} was introduced in \cite{hofmann86}
and has been widely used \cite{aster,engl00,hansen98,hansen10,mueller}:
If $\sigma_j=\mathcal{O}(j^{-\alpha})$, then \eqref{eq1}
is mildly or moderately ill-posed for $\frac{1}{2}<\alpha\le1$ or $\alpha>1$.
If $\sigma_j=\mathcal{O}(\rho^{-j})$ with $\rho>1$,
$j=1,2,\ldots,n$, then \eqref{eq1} is severely ill-posed.
Here for mildly ill-posed problems we add the requirement
$\alpha>\frac{1}{2}$, which does not appear in \cite{hofmann86}
but must be met for $k(s,t)\in L^2({\Omega\times\Omega})$
in \eqref{eq1} \cite{hanke93,hansen98}.
In the one-dimensional case, i.e., $q=1$, \eqref{eq1}
is severely ill-posed with $k(s,t)$ sufficiently smooth, and
it is moderately ill-posed with $\sigma_j=\mathcal{O}(j^{-p-1/2})$,
where $p$ is the highest order of continuous derivatives of
$k(s,t)$; see, e.g., \cite[p.8]{hansen98} and \cite[p.10-11]{hansen10}.
Clearly, the singular values $\sigma_j$ for a
severely ill-posed problem decay at the same rate $\rho^{-1}$,
while those of a moderately or mildly ill-posed problem decay
at the decreasing rate $\left(\frac{j}{j+1}\right)^{\alpha}$
that approaches one more quickly with $j$ for the mildly ill-posed problem
than for the moderately ill-posed problem.

If a regularized solution to \eqref{eq1} is at least as accurate as
$x_{k_0}^{tsvd}$, then it is called a best possible regularized solution.
Given \eqref{eq1}, if the regularized solution of an iterative regularization
solver at semi-convergence is such a best possible one,
then, by the words of Bj\"{o}rck and Eld\'{e}n,
the solver {\em works} for the problem and is said to have
the {\em full} regularization. Otherwise, the solver is said to have
only the {\em partial} regularization.

Because it has long been unknown whether or not LSQR, CGLS, LSMR and CGME
have the full regularization for a given \eqref{eq1},
one commonly combines them with some explicit
regularization, hoping that the resulting hybrid variants
find best possible regularized solutions
\cite{aster,hansen98,hansen10}. A hybrid CGLS is to run CGLS
for several trial regularization parameters $\lambda$ and
picks up the best one among the candidates \cite{aster}. Its
disadvantages are
that regularized solutions cannot be updated with different $\lambda$
and there is no guarantee that the selected regularized solution
is a best possible one.
The hybrid LSQR variants have been advocated by Bj\"{o}rck and Eld\'{e}n
\cite{bjorck79} and O'Leary and Simmons \cite{oleary81}, and improved and
developed by Bj\"orck \cite{bjorck88} and Bj\"{o}rck, Grimme and
van Dooren \cite{bjorck94}.
A hybrid LSQR first projects \eqref{eq1} onto Krylov
subspaces and then regularizes the projected problems explicitly.
It aims to remove the effects
of small Ritz values and expands Krylov subspaces until they
captures the $k_0$ dominant SVD components of $A$
\cite{bjorck88,bjorck94,hanke93,oleary81}. The explicit
regularization for projected problems should be introduced and play into
effects only after semi-convergence rather than from the very first iteration.
If it works, the error norms of regularized
solutions and the residual norms further decrease until they ultimately
stabilize. The hybrid LSQR and CGME have been intensively studied in, e.g.,
\cite{bazan10,bazan14,berisha,chung08,hanke01,hanke93,lewis09,neuman12,renaut}
and \cite{aster,hansen10,hansen13}.
Within the framework of such hybrid solvers, however,
it is hard to find a near-optimal
regularization parameter \cite{bjorck94,renaut}.

In contrast,
if an iterative solver is theoretically proved and practically identified to
have the full regularization, one simply stops it after
semi-convergence, and no complicated hybrid variant and further
iterations are needed. Obviously, we cannot emphasize too much the
importance of proving the full or partial
regularization of LSQR, CGLS, LSMR and CGME. By the definition of
the full or partial regularization,
we now modify the concern of Bj\"{o}rck and Eld\'{e}n as:
{\em  Do LSQR, CGLS, LSMR and CGME have the full or partial regularization for
severely, moderately and mildly ill-posed problems? How to identify
their full or partial regularization in practice?}

In this paper, we focus on LSQR and analyze its regularization for severely,
moderately and mildly ill-posed problems. Due to the mathematical equivalence
of CGLS and LSQR, the assertions on the full or partial regularization of
LSQR apply to CGLS as well.
We prove that LSQR has the full regularization for severely and
moderately ill-posed problems once $\rho>1$ and $\alpha>1$ suitably,
and it generally has only the partial regularization for mildly ill-posed
problems. In Section \ref{lsqr}, we describe the
Lanczos bidiagonalization process and LSQR, and make an
introductory analysis. In Section \ref{sine},
we establish $\sin\Theta$ theorems for the 2-norm
distance between the underlying $k$-dimensional Krylov subspace and the
$k$-dimensional dominant right singular subspace of $A$. We then derive
some follow-up results that play a central role in analyzing
the regularization of LSQR. In Section \ref{rankapp}, for the first
two kinds of problems
we prove that a $k$-step Lanczos bidiagonalization always
generates a near best rank $k$ approximation to $A$, and the $k$ Ritz values
always approximate the first $k$ large singular values in natural order,
and no small Ritz value appears for $k=1,2,\ldots,k_0$.
This will show that LSQR has the full regularization.
For mildly ill-posed problems,
we prove that, for some $k\leq k_0$, the $k$ Ritz values generally
do not approximate the first $k$ large singular values in natural order
and LSQR generally has only the partial regularization.
In Section \ref{alphabeta}, we derive bounds for the entries of bidiagonal
matrices generated by Lanczos bidiagonalization, showing how fast they
decay and how to use them to identify if
LSQR has the full regularization when
the degree of ill-posedness of \eqref{eq1} is unknown in advance.
In Section \ref{numer},
we report numerical experiments to confirm our theory on LSQR. Finally,
we summarize the paper with further remarks in Section \ref{concl}.

Throughout the paper, we denote by
$\mathcal{K}_{k}(C, w)= span\{w,Cw,\ldots,C^{k-1}w\}$
the $k$-dimensional Krylov subspace generated
by the matrix $\mathit{C}$ and the vector $\mathit{w}$, and by $I$ and the
bold letter $\mathbf{0}$ the identity matrix
and the zero matrix with orders clear from the context, respectively.
For the matrix $B=(b_{ij})$, we define $|B|=(|b_{ij}|)$,
and for $|C|=(|c_{ij}|)$, $|B|\leq |C|$ means
$|b_{ij}|\leq |c_{ij}|$ componentwise.

\section{The LSQR algorithm}\label{lsqr}

The LSQR algorithm is based on the Lanczos bidiagonalization process,
Algorithm~\ref{alg:lb}, that
computes two orthonormal bases $\{q_1,q_2,\dots,q_k\}$ and
$\{p_1,p_2,\dots,p_{k+1}\}$  of $\mathcal{K}_{k}(A^{T}A,A^{T}b)$ and
$\mathcal{K}_{k+1}(A A^{T},b)$  for $k=1,2,\ldots,n$,
respectively.
\begin{algorithm}
\begin{description}
   \item[1.] Take $ p_1=b/\|b\| \in \mathbb{R}^{m}$, and define $\beta_1{q_0}=0$.

   \item[2.] For $j=1,2,\ldots,k$
\begin{description}
  \item[(i)]
  $r = A^{T}p_j - \beta_j{q_{j-1}}$
  \item[(ii)] $\alpha_j = \|r\|;q_j = r/\alpha_j$
  \item[(iii)]
   $   z = Aq_j - \alpha_j{p_{j}}$
  \item[(iv)]
  $\beta_{j+1} = \|z\|;p_{j+1} = z/\beta_{j+1}.$
\end{description}
\end{description}
\caption{\ $k$-step Lanczos bidiagonalization process}
\label{alg:lb}
\end{algorithm}

Algorithm~\ref{alg:lb} can be written in the matrix form
\begin{align}
  AQ_k&=P_{k+1}B_k,\label{eqmform1}\\
  A^{T}P_{k+1}&=Q_{k}B_k^T+\alpha_{k+1}q_{k+1}e_{k+1}^{T}.\label{eqmform2}
\end{align}
where $e_{k+1}$ is the $(k+1)$-th canonical basis vector of
$\mathbb{R}^{k+1}$, $P_{k+1}=(p_1,p_2,\ldots,p_{k+1})$,
$Q_k=(q_1,q_2,\ldots,q_k)$ and
\begin{equation}\label{bk}
  B_k = \left(\begin{array}{cccc} \alpha_1 & & &\\ \beta_2 & \alpha_2 & &\\ &
  \beta_3 &\ddots & \\& & \ddots & \alpha_{k} \\ & & & \beta_{k+1}
  \end{array}\right)\in \mathbb{R}^{(k+1)\times k}.
\end{equation}
It is known from \eqref{eqmform1} that
\begin{equation}\label{Bk}
B_k=P_{k+1}^TAQ_k.
\end{equation}
We remind that the singular values of $B_k$, called the Ritz values of $A$ with
respect to the left and right subspaces $span\{P_{k+1}\}$ and $span\{Q_k\}$,
are all simple.

At iteration $k$, LSQR solves the problem
$\|Ax^{(k)}-b\|=\min_{x\in \mathcal{K}_k(A^TA,A^Tb)}
\|Ax-b\|$ and computes the iterates $x^{(k)}=Q_ky^{(k)}$ with
\begin{equation}\label{yk}
  y^{(k)}=\arg\min\limits_{y\in \mathbb{R}^{k}}\|B_ky-\|b\|e_1^{(k+1)}\|
  =\|b\| B_k^{\dagger} e_1^{(k+1)},
\end{equation}
where $e_1^{(k+1)}$ is the first canonical basis vector of $\mathbb{R}^{k+1}$,
and the residual norm $\|Ax^{(k)}-b\|$ decreases monotonically with respect to
$k$. We have $\|Ax^{(k)}-b\|=
\|B_k y^{(k)}-\|b\|e_1^{(k+1)}\|$ and $\|x^{(k)}\|=\|y^{(k)}\|$,
both of which can be cheaply computed.

Note that $\|b\|e_1^{(k+1)}=P_{k+1}^T b$. We have
\begin{equation}\label{xk}
x^{(k)}=Q_k B_k^{\dagger} P_{k+1}^Tb,
\end{equation}
that is, the iterate $x^{(k)}$ by LSQR is the minimum-norm least
squares solution to the perturbed
problem that replaces $A$ in \eqref{eq1} by its rank $k$ approximation
$P_{k+1}B_k Q_k^T$. Recall that the best rank $k$ approximation
$A_k$ to $A$ satisfies $\|A-A_k\|=\sigma_{k+1}$. Furthermore, analogous
to \eqref{tsvd}, LSQR now solves
\begin{equation}\label{lsqrreg}
\min\|x\| \ \ \mbox{ subject to }\ \ \|P_{k+1}B_kQ_k^Tx-b\|=\min,\ k=1,2,\ldots,n
\end{equation}
for the regularized solutions $x^{(k)}$ to \eqref{eq1}.
If $P_{k+1}B_k Q_k^T$ is a near best rank $k$ approximation
to $A$ with an approximate accuracy $\sigma_{k+1}$ and
the $k$ singular values of $B_k$ approximate the first $k$
large ones of $A$ in natural order for $k=1,2,\ldots,k_0$,
these two facts relate LSQR and the TSVD method naturally and
closely in two ways: (i) $x_k^{tsvd}$ and $x^{(k)}$ are the regularized
solutions to the two perturbed problems of \eqref{eq1} that replace $A$ by
its two rank $k$ approximations with the same quality, respectively;
(ii) $x_k^{tsvd}$ and $x^{(k)}$ solve almost the same two regularization
problems \eqref{tsvd} and \eqref{lsqrreg}, respectively.
As a consequence,
the LSQR iterate $x^{(k_0)}$ is as accurate as $x_{k_0}^{tsvd}$,
and LSQR has the full regularization. Otherwise,
as will be clear later, under the discrete Picard condition \eqref{picard},
$x^{(k_0)}$ cannot be as accurate as $x_{k_0}^{tsvd}$
if either $P_{k+1}B_k Q_k^T$ is not a near best rank $k$
approximation to $A$, $k=1,2,\ldots,k_0$, or $B_k$ has at least one
singular value smaller than $\sigma_{k_0+1}$ for some $k\leq k_0$.
Precisely, if either of them is violated for some $k\leq k_0$
and $\theta_k^{(k)}<\sigma_{k_0+1}$, $x^{(k)}$ has been deteriorated by the
noise $e$, and LSQR has only the partial regularization.
We will give a precise
definition of a near best rank $k$ approximation to $A$ soon.

\section{$\sin\Theta$ theorems for the distances between $\mathcal{K}_k(A^TA,A^Tb)$
and $span\{V_k\}$ as well as the others related} \label{sine}

van der Sluis and van der Vorst \cite{vorst86} prove the following result,
which has been used in Hansen \cite{hansen98} and the references therein
to illustrate the regularizing effects of LSQR and CGLS. We
will also investigate it further in our paper.

\begin{proposition}\label{help}
LSQR with the starting vector $p_1=b/\|b\|$ and CGLS
applied to $A^TAx=A^Tb$ with the starting vector
$x^{(0)}=0$ generate the same iterates
\begin{equation}\label{eqfilter2}
  x^{(k)}=\sum\limits_{i=1}^nf_i^{(k)}\frac{u_i^{T}b}{\sigma_i}v_i,\
  k=1,2,\ldots,n,
\end{equation}
where
\begin{equation}\label{filter}
f_i^{(k)}=1-\prod\limits_{j=1}^k\frac{(\theta_j^{(k)})^2-\sigma_i^2}
{(\theta_j^{(k)})^2},\ i=1,2,\ldots,n,
\end{equation}
and the $\theta_j^{(k)}$ are the singular values of $B_k$
labeled as $\theta_1^{(k)}>\theta_2^{(k)}>\cdots>\theta_k^{(k)}$.
\end{proposition}

\eqref{eqfilter2} shows that $x^{(k)}$ has a filtered SVD expansion of
form \eqref{eqfilter}. If all the Ritz values $\theta_j^{(k)}$
approximate the first $k$ singular values
$\sigma_j$ of $A$ in natural order, the filters $f_i^{(k)}\approx 1,\,
i=1,2,\ldots,k$ and the other $f_i^{(k)}$ monotonically approach zero
for $i=k+1,\ldots,n$. This indicates that if the $\theta_j^{(k)}$
approximate the first $k$ singular values
$\sigma_j$ of $A$ in natural order for $k=1,2,\ldots,k_0$ then the
$k_0$-step LSQR has the full regularization.
However, if a small Ritz value appears before some
$k\leq k_0$, i.e., $\theta_{k-1}^{(k)}>\sigma_{k_0+1}$ and
$\sigma_{j^*}<\theta_k^{(k)}\leq\sigma_{k_0+1}$ with the smallest integer
$j^*>k_0+1$, then $f_i^{(k)}\in (0,1)$ tends to zero
monotonically for $i=j^*,j^*+1,\ldots,n$; on the other hand,
we have
$$
\prod\limits_{j=1}^k\frac{(\theta_j^{(k)})^2-\sigma_i^2}
{(\theta_j^{(k)})^2}=\frac{(\theta_k^{(k)})^2-\sigma_i^2}
{(\theta_k^{(k)})^2}\prod\limits_{j=1}^{k-1}
\frac{(\theta_j^{(k)})^2-\sigma_i^2}{(\theta_j^{(k)})^2}\leq 0,
\ i=k_0+1,\ldots,j^*-1
$$
since the first factor is non-positive and the second factor is positive.
Then we get $f_i^{(k)}\geq 1,\ i=k_0+1,\ldots,j^*-1$, indicating that $x^{(k)}$
is already deteriorated and LSQR has only the partial regularization.

The standard $k$-step Lanczos bidiagonalization method computes the
$k$ Ritz values $\theta_j^{(k)}$, which are used to approximate some
singular values of $A$. It is mathematically equivalent to the
symmetric Lanczos method for
the eigenvalue problem of $A^TA$ starting with $q_1=A^Tb/\|A^Tb\|$;
see \cite{bai,bjorck96,bjorck15,parlett,vorst02} or \cite{reichel05,jia03,jia10} for
several variations that are based on standard, harmonic,
and refined projection \cite{bai,stewart01,vorst02}
or a combination of them \cite{jia05}.
It is known that, for general
singular value distribution and $b$, some Ritz
values become good approximations to the extreme
singular values of $A$ as $k$ increases.
If large singular values are well separated but
small singular values are clustered, large Ritz values
converge fast but small Ritz values converge slowly.

For \eqref{eq1}, $A^Tb$ contains more information on dominant right
singular vectors than on the ones corresponding to small singular values.
Therefore, $\mathcal{K}_k(A^TA,A^Tb)$ hopefully contains richer
information on the first $k$ right singular vectors $v_i$ than on the
other $n-k$ ones, at least for $k$ small.
Furthermore, note that $A$ has many small singular values clustered at zero.
Due to these two basic facts, all the Ritz values are expected to
approximate the large singular values of $A$ in natural order until some
iteration $k$, at which a small Ritz value shows up.
In this case, the iterates $x^{(k)}$ by LSQR capture only the largest $k$
dominant SVD components of $A$, and they
are deteriorated by the noise $e$ dramatically after that iteration.
This is why LSQR and CGLS have general regularizing effects; see,
e.g., \cite{aster,hansen98,hansen08,hansen10,hansen13} and the references
therein. Unfortunately, these arguments cannot help us
draw any definitive conclusion on the full or partial
regularization of LSQR because there has been no quantitative result on
the size of such $k$ for any kind of ill-posed problem and the noise $e$.
For a severely ill-posed example from
seismic tomography, it is reported in \cite{vorst90} that the desired
convergence of the Ritz values actually holds as long as the discrete
Picard condition is satisfied and there is a good separation
among the large singular values of $A$. Yet,
there has been no mathematical justification on these observations.

A complete understanding of the regularization
of LSQR includes accurate solutions of the following problems:
How accurately does $\mathcal{K}_{k}(A^{T}A, A^{T}b)$
approximate the $k$-dimensional dominant right singular subspace of
$A$? How accurate is the rank $k$ approximation $P_{k+1}B_kQ_k^T$ to $A$?
Can it be a near best rank $k$ approximation to $A$?
How does the noise level $\|e\|$ affects the approximation accuracy of
$\mathcal{K}_{k}(A^{T}A, A^{T}b)$ and $P_{k+1}B_kQ_k^T$ for $k\leq k_0$ and $k>k_0$,
respectively? What sufficient conditions on $\rho$ and $\alpha$ are needed
to guarantee that $P_{k+1}B_kQ_k^T$ is a near best rank $k$ approximation to $A$?
When do the $\theta_i^{(k)}$
approximate $\sigma_i,\ i=1,2,\ldots,k$ in natural order? When does at least
a small Ritz value appear, i.e., $\theta_k^{(k)}<\sigma_{k_0+1}$ before some
$k\leq k_0$? We will make a rigorous and detailed analysis on these problems,
present our results, and draw definitive
assertions on the regularization of LSQR for the three kinds of ill-posed
problems.

In terms of the canonical angles $\Theta(\mathcal{X},\mathcal{Y})$ between
two subspaces $\mathcal{X}$ and $\mathcal{Y}$ of equal
dimension \cite[p.43]{stewartsun}, we first present the following $\sin\Theta$
theorem, showing how $\mathcal{K}_{k}(A^{T}A, A^{T}b)$ approximates the
$k$-dimensional dominant right singular subspace $ span\{V_k\}$ of $A$ for
severely ill-posed problems.

\begin{theorem}\label{thm2}
Let the SVD of $A$ be as \eqref{eqsvd}. Assume that \eqref{eq1} is severely
ill-posed with $\sigma_j=\mathcal{O}(\rho^{-j})$ and $\rho>1$, $j=1,2,\ldots,n$,
and the discrete Picard condition \eqref{picard} is satisfied.
Let $\mathcal{V}_k=span\{V_k\}$ be the $k$-dimensional dominant right singular
subspace of $A$ spanned by the columns of $V_k=(v_1,v_2,\ldots,v_k)$
and $\mathcal{V}_k^R=\mathcal{K}_{k}(A^{T}A, A^{T}b)$. Then for
$k=1,2,\ldots,n-1$ we have
\begin{align}
\|\sin\Theta(\mathcal{V}_k,\mathcal{V}_k^R)\|&=
\frac{\|\Delta_k\|}{\sqrt{1+\|\Delta_k\|^2}}, \label{deltabound}\\
\|\tan\Theta(\mathcal{V}_k,\mathcal{V}_k^R)\|&=
\|\Delta_k\| \label{tantheta}
\end{align}

with $\Delta_k \in \mathbb{R}^{(n-k)\times k}$ to be defined by \eqref{defdelta}
and
\begin{equation}\label{k1}
\|\Delta_1\|\leq \frac{\sigma_{2}}{\sigma_1}\frac{|u_2^Tb|}{|u_1^Tb|}
\left(1+\mathcal{O}(\rho^{-2})\right),
\end{equation}
\begin{equation}\label{eqres1}
  \|\Delta_k\|\leq
  \frac{\sigma_{k+1}}{\sigma_k}\frac{|u_{k+1}^Tb|}{|u_{k}^Tb|}
  \left(1+\mathcal{O}(\rho^{-2})\right)
  |L_{k_1}^{(k)}(0)|,\ k=2,3,\ldots,n-1,
\end{equation}
where
\begin{equation}\label{lk}
|L_{k_1}^{(k)}(0)|=\max_{j=1,2,\ldots,k}|L_j^{(k)}(0)|,
\ |L_j^{(k)}(0)|=\prod\limits_{i=1,i\ne j}^k\frac{\sigma_i^2}{|\sigma_j^2-
\sigma_i^2|},\,j=1,2,\ldots,k.
\end{equation}
In particular, we have
\begin{align}
\|\Delta_1\|&\leq\frac{\sigma_2^{2+\beta}}{\sigma_1^{2+\beta}}
\left(1+\mathcal{O}(\rho^{-2})\right), \label{case5}\\
\|\Delta_k\|&\leq\frac{\sigma_{k+1}^{2+\beta}}{\sigma_k^{2+\beta}}
\left(1+\mathcal{O}(\rho^{-2})\right)|L_{k_1}^{(k)}(0)|,\
k=2,3,\ldots,k_0, \label{case1}\\
\|\Delta_k\|&\leq\frac{\sigma_{k+1}}{\sigma_k}\left(1+\mathcal{O}(\rho^{-2})\right)
|L_{k_1}^{(k)}(0)|, \ k=k_0+1,\ldots, n-1.\label{case2}
\end{align}
\end{theorem}

{\em Proof}.
Let $U_n=(u_1,u_2,\ldots,u_n)$ whose columns are the
first $n$ left singular vectors of $A$ defined by \eqref{eqsvd}.
Then the Krylov subspace $\mathcal{K}_{k}(\Sigma^2,
\Sigma U_n^Tb)=span\{DT_k\}$ with
\begin{equation*}\label{}
  D={\rm diag}(\sigma_i u_i^Tb)\in\mathbb{R}^{n\times n},\ \
  T_k=\left(\begin{array}{cccc} 1 &
  \sigma_1^2&\ldots & \sigma_1^{2k-2}\\
1 &\sigma_2^2 &\ldots &\sigma_2^{2k-2} \\
\vdots & \vdots&&\vdots\\
1 &\sigma_n^2 &\ldots &\sigma_n^{2k-2}
\end{array}\right).
\end{equation*}
Partition the diagonal matrix $D$ and the matrix $T_k$ as follows:
\begin{equation*}\label{}
  D=\left(\begin{array}{cc} D_1 & 0 \\ 0 & D_2 \end{array}\right),\ \ \
  T_k=\left(\begin{array}{c} T_{k1} \\ T_{k2} \end{array}\right),
\end{equation*}
where $D_1, T_{k1}\in\mathbb{R}^{k\times k}$. Since $T_{k1}$ is
a Vandermonde matrix with $\sigma_j$ distinct for $j=1,2,\ldots,k$,
it is nonsingular. Therefore, from $\mathcal{K}_{k}(A^{T}A, A^{T}b)=span\{VDT_k\}$
we have
\begin{equation}\label{kry}
\mathcal{V}_k^R=\mathcal{K}_{k}(A^{T}A, A^{T}b)=span
  \left\{V\left(\begin{array}{c} D_1T_{k1} \\ D_2T_{k2} \end{array}\right)\right\}
  =span\left\{V\left(\begin{array}{c} I \\ \Delta_k \end{array}\right)\right\},
\end{equation}
where
\begin{equation}\label{defdelta}
\Delta_k=D_2T_{k2}T_{k1}^{-1}D_1^{-1}\in \mathbb{R}^{(n-k)\times k}.
\end{equation}
Write $V=(V_k, V_k^{\perp})$, and define
\begin{equation}\label{zk}
Z_k=V\left(\begin{array}{c} I \\ \Delta_k \end{array}\right)
=V_k+V_k^{\perp}\Delta_k.
\end{equation}
Then $Z_k^TZ_k=I+\Delta_k^T\Delta_k$, and the columns of
$\hat{Z}_k=Z_k(Z_k^TZ_k)^{-\frac{1}{2}}$
form an orthonormal basis of $\mathcal{V}_k^R$. So we get an
orthogonal direct sum decomposition of $\hat{Z}_k$:
\begin{equation}\label{decomp}
\hat{Z}_k=(V_k+V_k^{\perp}\Delta_k)(I+\Delta_k^T\Delta_k)^{-\frac{1}{2}}.
\end{equation}
By definition and \eqref{decomp},  we obtain
$$
   \|\sin\Theta(\mathcal{V}_k,\mathcal{V}_k^R)\|
   =\|(V_k^{\perp})^T\hat{Z}_k\|
   =\|\Delta_k(I+\Delta_k^T\Delta_k)^{-\frac{1}{2}}\|
   =\frac{\|\Delta_k\|}{\sqrt{1+\|\Delta_k\|^2}},
$$
which is \eqref{deltabound}. From it, we get \eqref{tantheta} directly.

Next we estimate $\|\Delta_k\|$. For $k=2,3,\ldots,n-1$,
it is easily justified that the $j$-th column of $T_{k1}^{-1}$ consists of
the coefficients of the $j$-th Lagrange polynomial
\begin{equation*}\label{}
  L_j^{(k)}(\lambda)=\prod\limits_{i=1,i\neq j}^k
  \frac{\lambda-\sigma_i^2}{\sigma_j^2-\sigma_i^2}
\end{equation*}
that interpolates the elements of the $j$-th canonical basis vector
$e_j^{(k)}\in \mathbb{R}^{k}$ at the abscissas $\sigma_1^2,\sigma_2^2
\ldots, \sigma_k^2$. Consequently, the $j$-th column of $T_{k2}T_{k1}^{-1}$ is
\begin{equation}\label{tk12i}
  T_{k2}T_{k1}^{-1}e_j^{(k)}=(L_j^{(k)}(\sigma_{k+1}^2),\ldots,L_j^{(k)}
  (\sigma_{n}^2))^T, \ j=1,2,\ldots,k,
\end{equation}
from which we obtain
\begin{equation}\label{tk12}
  T_{k2}T_{k1}^{-1}=\left(\begin{array}{cccc} L_1^{(k)}(\sigma_{k+1}^2)&
  L_2^{(k)}(\sigma_{k+1}^2)&\ldots & L_k^{(k)}(\sigma_{k+1}^2)\\
L_1^{(k)}(\sigma_{k+2}^2)&L_2^{(k)}(\sigma_{k+2}^2) &\ldots &
L_k^{(k)}(\sigma_{k+2}^2) \\
\vdots & \vdots&&\vdots\\
L_1^{(k)}(\sigma_{n}^2)&L_2^{(k)}(\sigma_{n}^2) &\ldots &L_k^{(k)}(\sigma_{n}^2)
\end{array}\right)\in \mathbb{R}^{(n-k)\times k}.
\end{equation}
Since $|L_j^{(k)}(\lambda)|$ is monotonically
decreasing for $0\leq \lambda<\sigma_k^2$, it is bounded by $|L_j^{(k)}(0)|$.
With this property and the definition of $L_{k_1}^{(k)}(0)$, we get
\begin{align}
|\Delta_k|&=|D_2T_{k2}T_{k1}^{-1}D_1^{-1}| \notag \\
&\leq
\left(\begin{array}{cccc}
\frac{\sigma_{k+1}}{\sigma_1}\left|\frac{u_{k+1}^Tb}
{u_1^Tb}\right||L_{k_1}^{(k)}(0)| &\frac{\sigma_{k+1}}{\sigma_2}
\left|\frac{u_{k+1}^Tb}
{u_2^Tb}\right||L_{k_1}^{(k)}(0)| &\ldots&\frac{\sigma_{k+1}}{\sigma_k}
\left|\frac{u_{k+1}^Tb}{u_k^Tb}\right||L_{k_1}^{(k)}(0)| \\
\frac{\sigma_{k+2}}{\sigma_1}\left|\frac{ u_{k+2}^Tb}
{u_1^Tb}\right| |L_{k_1}^{(k)}(0)| &\frac{\sigma_{k+2}}{\sigma_2}
\left|\frac{u_{k+2}^Tb}
{u_2^Tb}\right| |L_{k_1}^{(k)}(0)|&
\ldots &\frac{\sigma_{k+2}}{\sigma_k}\left|\frac{u_{k+2}^Tb}
{u_k^Tb}\right| |L_{k_1}^{(k)}(0)| \\
\vdots &\vdots & &\vdots\\
\frac{\sigma_n}{\sigma_1}\left|\frac{u_n^Tb}
{u_1^Tb}\right| |L_{k_1}^{(k)}(0)| &\frac{\sigma_n}{\sigma_2}\left|\frac{u_n^Tb}
{u_2^Tb}\right| |L_{k_1}^{(k)}(0)|& \ldots &
\frac{\sigma_n}{\sigma_k}\left|\frac{u_n^Tb}{u_k^Tb}\right| |L_{k_1}^{(k)}(0)|
\end{array}
\right) \notag\\
&= |L_{k_1}^{(k)}(0)||\tilde\Delta_k|, \label{amplify}
\end{align}
where
\begin{equation}
|\tilde\Delta_k|=\left|(\sigma_{k+1} u_{k+1}^T b,\sigma_{k+2}u_{k+2}^Tb,
\ldots,\sigma_n u_n^T b)^T
\left(\frac{1}{\sigma_1 u_1^Tb},\frac{1}{\sigma_2 u_2^Tb},\ldots,
\frac{1}{\sigma_k u_k^Tb}\right)\right|  \label{delta1}
\end{equation}
is a rank one matrix. Therefore, by $\|C\|\leq \||C|\|$
(cf. \cite[p.53]{stewart98}), we get
\begin{align}
\|\Delta_k\| &\leq \||\Delta_k|\|\leq |L_{k_1}^{(k)}(0)|
\left\||\tilde\Delta_k|\right\| \notag\\
&=|L_{k_1}^{(k)}(0)|\left(\sum_{j=k+1}^n\sigma_j^2| u_j^Tb|^2\right)^{1/2}
\left(\sum_{j=1}^k \frac{1}{\sigma_j^2| u_j^Tb|^2}\right)^{1/2}.
\label{delta2}
\end{align}

By the discrete Picard condition \eqref{picard}, \eqref{picard1} and the
description between them,
for the white noise $e$, it is known from \cite[p.70-1]{hansen98} and
\cite[p.41-2]{hansen10} that
$
| u_j^T b|\approx | u_j^T \hat{b}|=\sigma_j^{1+\beta}
$
decrease as $j$ increases up to $k_0$ and then become stabilized as
$
| u_j^T b|\approx | u_j^T e |\approx \eta \approx \frac{\|e\|}{\sqrt{m}},
$
a small constant for $j>k_0$. In order to simplify the
derivation and present our results compactly, in terms of these
assumptions and properties, in later proofs we will use the following strict
equalities and inequalities:
\begin{align}
|u_j^T b|&= |u_j^T\hat{b}|=\sigma_j^{1+\beta},\ j=1,2,\ldots,k_0,\label{ideal3}\\
|u_j^T b|&=|u_j^T e|=\eta, \ j=k_0+1,\ldots,n,\label{ideal2}\\
|u_{j+1}^Tb|&\leq |u_j^Tb |,\ j=1,2,\ldots,n-1.\label{ideal}
\end{align}
From \eqref{ideal} and $\sigma_j=\mathcal{O}(\rho^{-j}),\ j=1,2,\ldots,n$,
for $k=1,2,\ldots,n-1$ we obtain
\begin{align}
\left(\sum_{j=k+1}^n\sigma_j^2| u_j^Tb|^2\right)^{1/2}
&= \sigma_{k+1}| u_{k+1}^Tb| \left(\sum_{j=k+1}^n
\frac{\sigma_j^2| u_j^Tb|^2}{\sigma_{k+1}^2| u_{k+1}^Tb|^2}\right)^{1/2}
\notag\\
&\leq \sigma_{k+1}| u_{k+1}^Tb| \left(\sum_{j=k+1}^n
\frac{\sigma_j^2}{\sigma_{k+1}^2}\right)^{1/2} \notag\\
&=\sigma_{k+1}| u_{k+1}^Tb|\left(1+\sum_{j=k+2}^n\mathcal{O}
(\rho^{2(k-j)+2})\right)^{1/2}
\notag \\
&=\sigma_{k+1}| u_{k+1}^Tb|\left(1+\mathcal{O}\left(\sum_{j=k+2}^n
\rho^{2(k-j)+2}\right)\right)^{1/2}
\notag \\
&=\sigma_{k+1}| u_{k+1}^Tb|\left(1+ \mathcal{O}\left(\frac{\rho^{-2}}
    {1-\rho^{-2}}\left(1-\rho^{-2(n-k-1)}\right)\right)\right)^{1/2}
\notag \\
&=\sigma_{k+1}| u_{k+1}^Tb| \left(1+\mathcal{O}(\rho^{-2})\right)^{1/2}\notag\\
&=\sigma_{k+1}| u_{k+1}^Tb| \left(1+\mathcal{O}(\rho^{-2})\right)
\label{severe1}
\end{align}
with $1+\mathcal{O}(\rho^{-2})=1$ for $k=n-1$.
For $k=2,3,\ldots,n-1$, from \eqref{ideal} we get
\begin{align*}
\left(\sum_{j=1}^k \frac{1}{\sigma_j^2| u_j^Tb|^2}\right)^{1/2}
&=\frac{1}{\sigma_k | u_k^T b|}\left(\sum_{j=1}^k\frac{\sigma_k^2| u_k^Tb|^2}
{\sigma_j^2| u_j^Tb|^2}\right)^{1/2}
\leq \frac{1}{\sigma_k | u_k^T b|}\left(\sum_{j=1}^k\frac{\sigma_k^2}
{\sigma_j^2}\right)^{1/2} \\
&=\frac{1}{\sigma_k | u_k^T b|}\left(1+\mathcal{O}\left(\sum_{j=1}^{k-1}
\rho^{2(j-k)}\right)\right)^{1/2} \\
&=\frac{1}{\sigma_k | u_k^T b|}\left(1+\mathcal{O}(\rho^{-2})\right).
\end{align*}
From the above and \eqref{delta2}, we finally obtain \eqref{eqres1} by noting
$$
\|\Delta_k\|\leq \frac{\sigma_{k+1}}{\sigma_k}\frac{| u_{k+1}^T b|}
{| u_k^T b|}\left(1+\mathcal{O}(\rho^{-2})\right)|L_{k_1}^{(k)}(0)|,
\ k=2,3,\ldots,n-1.
$$

Note that the Lagrange polynomials $L_j^{(k)}(\lambda)$ require $k\geq 2$. So,
we need to treat the case $k=1$ independently.
Observe from \eqref{defdelta} and \eqref{ideal} that
$$
T_{k2}=(1,1,\ldots,1)^T,\ D_2T_{k2}=(\sigma_2u_2^Tb,\sigma_3 u_3^Tb,
\ldots,\sigma_n u_n^Tb)^T,\
T_{k1}^{-1}=1,\ D_1^{-1}=\frac{1}{\sigma_1 u_1^Tb}.
$$
Therefore, we have
\begin{equation}\label{deltaexp}
\Delta_1=(\sigma_2u_2^Tb,\sigma_3 u_3^Tb,
\ldots,\sigma_n u_n^Tb)^T\frac{1}{\sigma_1 u_1^Tb},
\end{equation}
from which and \eqref{severe1} for $k=1$ it is direct to get \eqref{k1}.

In terms of the discrete Picard condition \eqref{picard},
\eqref{picard1}, \eqref{ideal3} and \eqref{ideal2}, we have
\begin{equation}\label{ratio1}
\frac{|u_{k+1}^Tb|}{| u_k^T b|}= \frac{|u_{k+1}^T\hat{b}|}{| u_k^T\hat{b}|}
=\frac{\sigma_{k+1}^{1+\beta}}{\sigma_k^{1+\beta}}, \ k\leq k_0,
\end{equation}
\begin{equation}\label{ratio2}
\frac{|u_{k+1}^Tb|}{| u_k^T b|}= \frac{|u_{k+1}^T e|}{| u_k^T e|}
=1,\  k>k_0.
\end{equation}
Applying them to \eqref{k1} and \eqref{eqres1} establishes \eqref{case5},
\eqref{case1} and \eqref{case2}, respectively.
\qed

We next estimate the factor $|L_{k_1}^{(k)}(0)|$ and all $|L_j^{(k)}(0)|$,
$j=1,2,\ldots,k$ accurately.

\begin{theorem}\label{estlk}
For the severely ill-posed problem, we have
\begin{align}
|L_k^{(k)}(0)|&=1+\mathcal{O}(\rho^{-2}), \label{lkkest}\\
|L_j^{(k)}(0)|&=\frac{1+\mathcal{O}(\rho^{-2})}
{\prod\limits_{i=j+1}^k\left(\frac{\sigma_{j}}{\sigma_i}\right)^2}
=\frac{1+\mathcal{O}(\rho^{-2})}{\mathcal{O}(\rho^{(k-j)(k-j+1)})},
\ j=1,2,\ldots,k-1, \label{lj0}\\
|L_{k_1}^{(k)}(0)|&=\max_{j=1,2,\ldots,k}|L_j^{(k)}(0)|
=1+\mathcal{O}(\rho^{-2}). \label{lkk}
\end{align}
\end{theorem}

{\em Proof}.
Exploiting the Taylor series expansion and
$\sigma_i=\mathcal{O}(\rho^{-i})$ for $i=1,2,\ldots,n$,
by definition, for $j=1,2,\ldots,k-1$ we have
\begin{align}
|L_j^{(k)}(0)|&=\prod\limits_{i=1,i\neq j}^k
\left|\frac{\sigma_i^2}{\sigma_i^2-\sigma_j^2}\right|
  =\prod\limits_{i=1}^{j-1}\frac{\sigma_i^2}{\sigma_i^2-\sigma_j^2}
   \cdot\prod\limits_{i=j+1}^{k}\frac{\sigma_i^2}{\sigma_j^2-\sigma_{i}^2}
   \notag\\
& =\prod\limits_{i=1}^{j-1}\frac{1}
{1-\mathcal{O}(\rho^{-2(j-i)})}
\prod\limits_{i=j+1}^{k}\frac{1}
{1-\mathcal{O}(\rho^{-2(i-j)})}\frac{1}
{\prod\limits_{i=j+1}^{k}\mathcal{O}(\rho^{2(i-j)})} \notag\\
&=\frac{\left(1+\sum\limits_{i=1}^j \mathcal{O}(\rho^{-2i})\right)
\left(1+\sum\limits_{i=1}^{k-j+1} \mathcal{O}(\rho^{-2i})\right)}
{\prod\limits_{i=j+1}^{k}\mathcal{O}(\rho^{2(i-j)})} \label{lik}
\end{align}
by absorbing those higher order terms into $\mathcal{O}(\cdot)$
in the numerator. For $j=k$, we get
\begin{align*}
|L_k^{(k)}(0)|&=\prod\limits_{i=1}^{k-1}
\left|\frac{\sigma_i^2}{\sigma_i^2-\sigma_{k}^2}\right|
=\prod\limits_{i=1}^{k-1}\frac{1}
{1-\mathcal{O}(\rho^{-2(k-i)})}=
\prod\limits_{i=1}^{k-1}\frac{1}
{1-\mathcal{O}(\rho^{-2i})}\\
&=1+\sum\limits_{i=1}^k \mathcal{O}(\rho^{-2i})
=1+\mathcal{O}\left(\sum\limits_{i=1}^k\rho^{-2i}\right)\\
&=1+ \mathcal{O}\left(\frac{\rho^{-2}}
    {1-\rho^{-2}}(1-\rho^{-2k})\right)
=1+\mathcal{O}(\rho^{-2}),
\end{align*}
which is \eqref{lkkest}.

Note that for the numerator of \eqref{lik} we have
  $$
  1+\sum\limits_{i=1}^j \mathcal{O}(\rho^{-2i})
    =1+ \mathcal{O}\left(\sum\limits_{i=1}^j\rho^{-2i}\right)
    =1+ \mathcal{O}\left(\frac{\rho^{-2}}
    {1-\rho^{-2}}(1-\rho^{-2j})\right),
  $$
  and
  $$
    1+\sum\limits_{i=1}^{k-j+1} \mathcal{O}(\rho^{-2i})
    =1+ \mathcal{O}\left(\sum\limits_{i=1}^{k-j+1}\rho^{-2i}\right)
    =1+ \mathcal{O}\left(\frac{\rho^{-2}}{1-\rho^{-2}}
    (1-\rho^{-2(k-j+1)})\right),
  $$
whose product for any $k$ is
  $$
  1+ \mathcal{O}\left(\frac{2\rho^{-2}}{1-\rho^{-2}}\right)
  +\mathcal{O}\left(\left(\frac{\rho^{-2}}{1-\rho^{-2}}\right)^2\right)=
  1+ \mathcal{O}\left(\frac{2\rho^{-2}}{1-\rho^{-2}}\right)
  =  1+\mathcal{O}(\rho^{-2}).
  $$
On the other hand, note that the denominator of \eqref{lik} is defined by
$$
\prod\limits_{i=j+1}^k\left(\frac{\sigma_{j}}{\sigma_i}\right)^2
=\prod\limits_{i=j+1}^{k}\mathcal{O}(\rho^{2(i-j)})
=\mathcal{O}((\rho\cdot\rho^2\cdots\rho^{k-j})^2)
=\mathcal{O}(\rho^{(k-j)(k-j+1)}),
$$
which, together with the above estimate
for the numerator of \eqref{lik}, proves \eqref{lj0}.
Notice that the above quantity
is always {\em bigger than one} for $j=1,2,\ldots,k-1$.
Therefore, for any $k$, combining \eqref{lkkest} and \eqref{lj0}
gives \eqref{lkk}.
\qed

\begin{remark}\label{severerem}
From \eqref{lkk}, the results in Theorem~\ref{thm2} are simplified as
\begin{align}
\|\Delta_k\|&\leq\frac{\sigma_{k+1}^{2+\beta}}{\sigma_k^{2+\beta}}
\left(1+\mathcal{O}(\rho^{-2})\right),\ k=1,2,\ldots, k_0,
\label{case3}\\
\|\Delta_k\|&\leq\frac{\sigma_{k+1}}{\sigma_k}\left(1+\mathcal{O}(\rho^{-2})\right),
\ k=k_0+1,\ldots,n-1. \label{case4}
\end{align}
\end{remark}

\begin{remark}
It is seen from the proof that $k_1$ must be close to $k$ or equals $k$.
\eqref{lj0} illustrates that $|L_j^{(k)}(0)|$ increases fast, up to
$1+\mathcal{O}(\rho^{-2})$,
with $j$ increasing, and the smaller $j$, the smaller $|L_j^{(k)}(0)|$.
\eqref{case3} and \eqref{case4} indicate that $\mathcal{V}_k^R$ captures
$\mathcal{V}_k$ better for $k\leq k_0$ than for $k>k_0$. That is,
after iteration $k_0$, the noise $e$ starts to
impair the ability of $\mathcal{V}_k^R$ to capture
$\mathcal{V}_k$.
\end{remark}

Next we estimate $\|\sin\Theta(\mathcal{V}_k,\mathcal{V}_k^R)\|$ for
moderately and mildly ill-posed problems.

\begin{theorem}\label{moderate}
Assume that \eqref{eq1} is moderately or mildly ill-posed with $\sigma_j=
\zeta j^{-\alpha},\ j=1,2,\ldots,n$,
where $\alpha>\frac{1}{2}$ and $\zeta>0$ is some constant.
Then \eqref{deltabound} and \eqref{tantheta}  hold with
\begin{align}
\|\Delta_1\|&\leq \frac{\sigma_2^{1+\beta}}{\sigma_1^{1+\beta}}
\sqrt{\frac{1}{2\alpha-1}},\label{mod1}\\
\|\Delta_k\|&\leq \frac{\sigma_{k+1}^{1+\beta}}{\sigma_k^{1+\beta}}
\sqrt{\frac{k^2}{4\alpha^2-1}+\frac{k}{2\alpha-1}}|L_{k_1}^{(k)}(0)|, \
k=2,3,\ldots, k_0, \label{modera2} \\
\|\Delta_k\|&\leq
\sqrt{\frac{k^2}{4\alpha^2-1}+\frac{k}{2\alpha-1}}|L_{k_1}^{(k)}(0)|,\
k=k_0+1,\ldots, n-1.
\label{modera3}
\end{align}
\end{theorem}

{\em Proof}.
Following the proof of Theorem~\ref{thm2}, we know that $|\Delta_k|\leq
|L_{k_1}^{(k)}(0)||\tilde\Delta_k|$ still holds with
$\tilde{\Delta}_k$ defined by \eqref{delta1}. So we only need to
bound the right-hand side of \eqref{delta2}. For $k=1,2,\ldots, n-1$, from
\eqref{ideal} we get
\begin{align}
\left(\sum_{j=k+1}^n\sigma_j^2| u_j^Tb|^2\right)^{1/2}
&= \sigma_{k+1}| u_{k+1}^Tb| \left(\sum_{j=k+1}^n
\frac{\sigma_j^2| u_j^Tb|^2}{\sigma_{k+1}^2| u_{k+1}^Tb|^2}\right)^{1/2}
\notag\\
&\leq \sigma_{k+1}| u_{k+1}^Tb| \left(\sum_{j=k+1}^n
\frac{\sigma_j^2}{\sigma_{k+1}^2}\right)^{1/2} \notag\\
&= \sigma_{k+1}| u_{k+1}^Tb| \left(\sum_{j=k+1}^n \left(\frac{j}{k+1}
\right)^{-2\alpha}\right)^{1/2} \notag \\
&=\sigma_{k+1}| u_{k+1}^Tb|
\left((k+1)^{2\alpha}\sum_{j=k+1}^n \frac{1}{j^{2\alpha}}\right)^{1/2}
\notag\\
&< \sigma_{k+1}| u_{k+1}^Tb| (k+1)^{\alpha}\left(\int_k^{\infty}
\frac{1}{x^{2\alpha}} dx\right)^{1/2}
\notag \\
&= \sigma_{k+1}| u_{k+1}^Tb|\left(\frac{k+1}{k}\right)^{\alpha}
\sqrt{\frac{k}{2\alpha-1}} \notag\\
&=\sigma_{k+1}| u_{k+1}^Tb|\frac{\sigma_k}{\sigma_{k+1}} \sqrt{\frac{k}
{2\alpha-1}} \notag\\
&=\sigma_k | u_{k+1}^Tb|\sqrt{\frac{k}
{2\alpha-1}}.
\label{modeest}
\end{align}
Since the function $x^{2\alpha}$ with any $\alpha> \frac{1}{2}$
is convex over the interval $[0,1]$, for $k=2,3,\ldots, n-1$, from \eqref{ideal}
we obtain
\begin{align}
\left(\sum_{j=1}^k \frac{1}{\sigma_j^2| u_j^Tb|^2}\right)^{1/2}
&=\frac{1}{\sigma_k | u_k^T b|}\left(\sum_{j=1}^k
\frac{\sigma_k^2| u_k^Tb|^2}
{\sigma_j^2| u_j^Tb|^2}\right)^{1/2}
\leq\frac{1}{\sigma_k | u_k^T b|}\left(\sum_{j=1}^k\frac{\sigma_k^2}
{\sigma_j^2}\right)^2 \notag \\
&=\frac{1}{\sigma_k | u_k^T b|}
\left(\sum_{j=1}^k \left(\frac{j}{k}
\right)^{2\alpha}\right)^{1/2} \notag\\
&=\frac{1}{\sigma_k | u_k^T b|}
\left(k\sum_{j=1}^{k} \frac{1}{k}\left(\frac{j-1}{k}
\right)^{2\alpha}+1\right)^{1/2} \label{sum1} \\
&< \frac{1}{\sigma_k | u_k^T b|} \left(k\int_0^1
x^{2\alpha}dx+1\right)^{1/2} \notag\\
&=\frac{1}{\sigma_k | u_k^T b|} \sqrt{\frac{k}{2\alpha+1}+1}. \label{estimate2}
\end{align}
Substituting the above and \eqref{modeest} into \eqref{delta2}
and exploiting \eqref{ratio1} and \eqref{ratio2}, we obtain
\eqref{modera2} and \eqref{modera3}. For $k=1$,
it follows from \eqref{deltaexp} and \eqref{ratio1} that \eqref{mod1} holds.
\qed

\begin{remark}
For the sake of precise presentation,
we have used the simplified singular value
model $\sigma_j=\zeta j^{-\alpha}$ to replace the general form
$\sigma_j=\mathcal{O}(j^{-\alpha})$, where the constant in each
$\mathcal{O}(\cdot)$ is implicit. This model,
though simple, reflects the essence of moderately and mildly ill-posed
problems and avoids some non-transparent formulations.
\end{remark}


Unlike the severely ill-posed problem case, for moderately and mildly ill-posed
problems it appears impossible to estimate $|L_{k_1}^{(k)}(0)|$ both elegantly
and accurately. We present the following results on
$|L_j^{(k)}(0)|,\ j=1,2,\ldots,k$ and $|L_{k_1}^{(k)}(0)|$.

\begin{proposition}\label{estlk2}
For the moderately and mildly ill-posed problems with
$\sigma_i=\zeta i^{-\alpha},\ i=1,2,\ldots,n$ and $\alpha>\frac{1}{2}$, we have
\begin{align}
|L_k^{(k)}(0)|&\approx 1+\frac{k}{2\alpha+1}, \label{lkkmoderate}\\
|L_j^{(k)}(0)|&\approx  \left(1+\frac{j}{2\alpha+1}\right)\left(1+\frac{j-j^{2\alpha}
k^{-2\alpha+1}}{2\alpha-1}\right)
\prod_{i=j+1}^k\left(\frac{j}{i}\right)^{2\alpha},\,j=1,2,\ldots,k-1. \label{lkjmod}
\end{align}
For $\alpha>1$, we have
\begin{align}
|L_j^{(k)}(0) |&\approx\left(1+\frac{j}{2\alpha+1}\right)
\prod_{i=j+1}^k\left(\frac{j}{i}\right)^{2\alpha},\,j=1,2,\ldots,k-1,
\label{lkjmod1}\\
\frac{k}{2\alpha+1}&<|L_{k_1}^{(k)}(0)|\approx 1+\frac{k}{2\alpha+1}
\label{lk1size}
\end{align}
with the lower bound requiring $k$ satisfying $\frac{2\alpha+1}{k}\leq 1$;
for $\frac{1}{2}<\alpha\leq 1$ and $k$ satisfying $\frac{2\alpha+1}{k}\leq 1$,
we have
\begin{equation}\label{lk1sizemild}
\frac{k}{2\alpha+1}<|L_{k_1}^{(k)}(0)|.
\end{equation}
\end{proposition}

{\em Proof}.
Exploiting the first order Taylor expansion, we obtain
estimate
\begin{align*}
|L_k^{(k)}(0)|&=
\prod\limits_{i=1}^{k-1}\frac{\sigma_i^2}{\sigma_i^2-\sigma_k^2}
=\prod\limits_{i=1}^{k-1}\frac{1}
{1-(\frac{i}{k})^{2\alpha}} \notag\\
&\approx 1+\sum\limits_{i=1}^{k-1}\left(\frac{i}{k}\right)^{2\alpha}
=1+k\sum\limits_{i=1}^k\frac{1}{k}\left(\frac{i-1}{k}\right)^{2\alpha}
\notag\\
&\approx1+k\int_0^1 x^{2\alpha}dx=1+\frac{k}{2\alpha+1},
\end{align*}
which proves \eqref{lkkmoderate}.

For $j=1,2,\ldots,k-1$, by the definition of $\sigma_i$,
since $\alpha\geq \frac{1}{2}$, we have
\begin{align*}
|L_j^{(k)}(0)|&=\prod\limits_{i=1,i\neq j}^k
\left|\frac{\sigma_i^2}{\sigma_i^2-\sigma_j^2}\right|
  =\prod\limits_{i=1}^{j-1}\frac{\sigma_i^2}{\sigma_i^2-\sigma_j^2}
   \cdot\prod\limits_{i=j+1}^{k}\frac{\sigma_i^2}{\sigma_j^2-\sigma_{i}^2}
   \\
& =\prod\limits_{i=1}^{j-1}\frac{1}
{1-\left(\frac{i}{j}\right)^{2\alpha}}
\prod\limits_{i=j+1}^{k}\frac{1}
{1-\left(\frac{j}{i}\right)^{2\alpha}}\frac{1}
{\prod\limits_{i=j+1}^{k}\left(\frac{i}{j}\right)^{2\alpha}} \\
&\approx \left(1+\sum\limits_{i=1}^{j-1}\left(\frac{i}{j}\right)^{2\alpha}\right)
\left(1+\sum\limits_{i=j+1}^{k} \left(\frac{j}{i}\right)^{2\alpha}\right)
{\prod\limits_{i=j+1}^{k}\left(\frac{j}{i}\right)^{2\alpha}}\\
&\leq \left(1+\int_0^1 x^{2\alpha} dx\right)\left(1+j^{2\alpha}\int_j^k \frac{1}
{x^{2\alpha}}dx\right){\prod\limits_{i=j+1}^{k}\left(\frac{j}{i}\right)^{2\alpha}}\\
&= \left(1+\frac{j}{2\alpha+1}\right)\left(1+\frac{j-j^{2\alpha}
k^{-2\alpha+1}}{2\alpha-1}\right)
\prod_{i=j+1}^k\left(\frac{j}{i}\right)^{2\alpha}.
\end{align*}
Note that $\prod_{i=j+1}^k\left(\frac{j}{i}\right)^{2\alpha}$ are always smaller
than one for $j=1,2,\ldots,k-1$, and the smaller $j$ is, the smaller this
factor is. Furthermore, exploiting
$$
\left(\frac{j}{k}\right)^{k-j}<\prod_{i=j+1}^k\frac{j}{i}
<\left(\frac{j}{j+1}\right)^{k-j}
$$
and by some elementary manipulation, for $\alpha>1$ we can justify the
estimates
$$
\frac{j-j^{2\alpha}
k^{-2\alpha+1}}{2\alpha-1}
\prod_{i=j+1}^k\left(\frac{j}{i}\right)^{2\alpha}\approx 0, \ j=1,2,\ldots,k-1.
$$
As a result, for $\alpha>1$ we have
$$
|L_j^{(k)}(0) |\approx\left(1+\frac{j}{2\alpha+1}\right)
\prod_{i=j+1}^k\left(\frac{j}{i}\right)^{2\alpha},\,j=1,2,\ldots,k-1,
$$
which establishes \eqref{lkjmod1}. A combination of it and \eqref{lkkmoderate}
gives the right-hand part of \eqref{lk1size}.

On the other hand, once $k$ is such that $\frac{2\alpha+1}{k}\leq 1$, we
always have
\begin{align}
 |L_{k_1}^{(k)}(0)|&\geq |L_{k}^{(k)}(0)|=
\prod\limits_{i=1}^{k-1}\frac{\sigma_i^2}{\sigma_i^2-\sigma_{k}^2}
=\prod\limits_{i=1}^{k-1}\frac{1}
{1-(\frac{i}{k})^{2\alpha}} \notag\\
&> 1+\sum\limits_{i=1}^{k-1}\left(\frac{i}{k}\right)^{2\alpha}
>1+k\int_0^{\frac{k-1}{k}} x^{2\alpha}dx \notag\\
&=1+\frac{k\left(\frac{k-1}{k}\right)^{2\alpha+1}}{2\alpha+1}
\approx 1+\frac{k}{2\alpha+1}\left(1-\frac{2\alpha+1}{k}\right)=
\frac{k}{2\alpha+1}, \label{lowerbound}
\end{align}
which yields the lower bound of \eqref{lk1size} and \eqref{lk1sizemild}.
\qed

\begin{remark}
The inaccuracy source of \eqref{lkkmoderate} and \eqref{lkjmod}
consists in using $\sum$ to replace $\prod$ approximately in the proof. They
are considerable underestimates for $\frac{1}{2}<\alpha\leq 1$ but are accurate,
provided that $\alpha>1$ suitably;
the bigger $\alpha$ is, the more accurate the estimates \eqref{lkkmoderate} and
\eqref{lkjmod} are. The derivation of \eqref{lowerbound}
indicates that $|L_{k_1}^{(k)}(0)|$ can be bigger than
$\frac{k}{2\alpha+1}$ substantially for $\frac{1}{2}<\alpha\leq 1$,
particularly when $\alpha$ is close to $\frac{1}{2}$; in this case,
we cannot bound $|L_{k_1}^{(k)}(0)|$ from above
since \eqref{lkjmod} is a considerable underestimate and
the denominator $2\alpha-1$ in \eqref{lkjmod} can be very small.
\end{remark}

\begin{remark}
It is easily seen from \eqref{deltabound} that
$\|\sin\Theta(\mathcal{V}_k,\mathcal{V}_k^R)\|$ increases monotonically
with respect to $\|\Delta_k\|$. For $\|\Delta_k\|$ reasonably small and
$\|\Delta_k\|$ large we have
$$
\|\sin\Theta(\mathcal{V}_k,\mathcal{V}_k^R)\|\approx \|\Delta_k\|
\ \mbox{ and } \
\|\sin\Theta(\mathcal{V}_k,\mathcal{V}_k^R)\|\approx 1,
$$
respectively. From \eqref{picard} and \eqref{picard1}, we obtain
$k_0=\lfloor\eta^{-\frac{1}{\alpha(1+\beta)}}\rfloor-1$,
where $\lfloor\cdot\rfloor$ is the Gaussian function.
As a result, for $\alpha>1$,
$k_0$ is typically small and at most modest for a practical noise $e$ with
$\|e\|\approx\sqrt{m}\eta$ since $\frac{\|e\|}{\|\hat{b}\|}$ typically ranges
from $10^{-4}$ to $10^{-2}$. This means that for a moderately ill-posed problem
$\|\Delta_k\|$ is at most modest and cannot be large, so that
$\|\sin\Theta(\mathcal{V}_k,\mathcal{V}_k^R)\|<1$ fairly.
\end{remark}

\begin{remark}
For severely ill-posed problems, since all the $\frac{\sigma_{k+1}}{\sigma_k}
\sim \rho^{-1}$, \eqref{case3} and \eqref{case4}
indicate that $\|\sin\Theta(\mathcal{V}_k,\mathcal{V}_k^R)\|$
is essentially unchanged for $k=1,2,\ldots,k_0$ and $k=k_0+1,\ldots,n-1$,
respectively, meaning that $\mathcal{V}_k^R$ captures $\mathcal{V}_k$ with almost
the same accuracy for $k\leq k_0$ and $k>k_0$, respectively. However,
the situation is different for moderately ill-posed
problems. For them,  $\frac{\sigma_{k+1}}{\sigma_k}=
\left(\frac{k}{k+1}\right)^{\alpha}$ increases slowly as $k$ increases, and
$\sqrt{\frac{k^2}{4\alpha^2-1}+\frac{k}{2\alpha-1}}
|L_{k_1}^{(k)}(0)|$ increases as $k$ grows. Therefore, \eqref{modera2} and
\eqref{modera3} illustrate that $\|\sin\Theta(\mathcal{V}_k,\mathcal{V}_k^R)\|$
increases slowly with $k\leq k_0$ and $k> k_0$, respectively.
This means that $\mathcal{V}_k^R$ may not capture
$\mathcal{V}_k$ so well as it does for severely ill-posed problems as $k$
increases. In particular, starting with some $k>k_0$,
$\|\sin\Theta(\mathcal{V}_k,\mathcal{V}_k^R)\|$ starts to approach one,
which indicates that, for $k$ big, $\mathcal{V}_k^R$ will contain substantial
information on the right singular vectors corresponding to the $n-k$ small
singular values of $A$.
\end{remark}

\begin{remark}\label{mildrem}
For mildly ill-posed problems with $\frac{1}{2}<
\alpha\leq 1$, there are some distinctive features. Note
from \eqref{picard} and \eqref{picard1} that
$k_0$ is now considerably bigger than
that for a severely or moderately ill-posed problem
with the same noise level $\|e\|$ and $\beta$. As a result,
firstly, for $\alpha\leq 1$ and the same $k$, the factor $\frac{\sigma_{k+1}}
{\sigma_{k}}=\left(\frac{k}{k+1}\right)^{\alpha}$ is
bigger than that for the moderately ill-posed problem;
secondly, $\sqrt{\frac{k^2}{4\alpha^2-1}+\frac{k}{2\alpha-1}} \sim k$
if $\alpha\approx 1$ and is much bigger than $k$ and can be
arbitrarily large if $\alpha\approx \frac{1}{2}$; thirdly,
\eqref{lowerbound} and the comment on it
indicate that $|L_{k_1}^{(k)}(0)|$ is bigger than one
considerably for $\frac{1}{2}<\alpha\leq 1$ as
$k$ increases up to $k_0$.
The bound \eqref{modera2} thus becomes increasingly large as $k$
increases up to $k_0$ for mildly ill-posed
problems, causing that $\|\Delta_k\|$ is large and
$\|\sin\Theta(\mathcal{V}_k,\mathcal{V}_k^R)\|\approx 1$ starting with
some $k\leq k_0$. Consequently, $\mathcal{V}_{k_0}^R$ cannot
effectively capture $\mathcal{V}_{k_0}$ and
contains substantial information on the right singular vectors
corresponding to the $n-k_0$ small singular values.
\end{remark}

Before proceeding, we tentatively investigate how
$\|\sin\Theta(\mathcal{V}_k,\mathcal{V}_k^R)\|$ affects the smallest Ritz
value $\theta_k^{(k)}$. This problem is of central importance for
understanding the regularizing effects of LSQR. We aim to lead the reader
to a first manifestation that (i) we may have $\theta_k^{(k)}>\sigma_{k+1}$,
that is, no small Ritz value may appear when
$\|\sin\Theta(\mathcal{V}_k,\mathcal{V}_k^R)\|<1$ suitably, and (ii) we must have
$\theta_k^{(k)}\leq\sigma_{k+1}$, that is, $\theta_k^{(k)}$ cannot approximate
$\sigma_k$ in natural order, meaning that $\theta_k^{(k)}\leq\sigma_{k_0+1}$
no later than iteration $k_0$, once
$\|\sin\Theta(\mathcal{V}_k,\mathcal{V}_k^R)\|$ is sufficiently close to one.

\begin{theorem}\label{initial}
Let $\|\sin\Theta(\mathcal{V}_k,\mathcal{V}_k^R)\|^2=1-\varepsilon_k^2$ with
$0< \varepsilon_k< 1$, $k=1,2,\ldots,n-1$, and let the unit-length
$\tilde{q}_k\in\mathcal{V}_k^R$
be a vector that has the smallest acute angle with $span\{V_k^{\perp}\}$, i.e.,
the closest to $span\{V_k^{\perp}\}$, where $V_k^{\perp}$ is the matrix consisting
of the last $n-k$ columns of $V$ defined by \eqref{eqsvd}. Then it holds that
\begin{equation}\label{rqi}
\varepsilon_k^2\sigma_k^2+
(1-\varepsilon_k^2)\sigma_n^2< \tilde{q}_k^TA^TA\tilde{q}_k<
\varepsilon_k^2\sigma_{k+1}^2+
(1-\varepsilon_k^2)\sigma_1^2.
\end{equation}
If $\varepsilon_k\geq \frac{\sigma_{k+1}}{\sigma_k}$,
then
\begin{equation}
\sqrt{\tilde{q}_k^TA^TA\tilde{q}_k}>\sigma_{k+1};
\label{est1}
\end{equation}
if $\varepsilon_k^2\leq\frac{\delta}
{(\frac{\sigma_1}{\sigma_{k+1}})^2-1}$ for a given arbitrarily small
$\delta>0$, then
\begin{equation}\label{thetasigma}
\theta_k^{(k)}<(1+\delta)^{1/2}\sigma_{k+1}.
\end{equation}
\end{theorem}

{\em Proof}.
Since the columns of $Q_k$ generated by Lanczos bidiagonalization form an
orthonormal basis of $\mathcal{V}_k^R$, by definition and the assumption on
$\tilde{q}_k$ we have
\begin{align}
\|\sin\Theta(\mathcal{V}_k,\mathcal{V}_k^R)\|&=\|(V_k^{\perp})^TQ_k\|
=\|V_k^{\perp}(V_k^{\perp})^TQ_k\| \notag\\
&=\max_{\|c\|=1}\|V_k^{\perp}(V_k^{\perp})^TQ_kc\|
=\|V_k^{\perp}(V_k^{\perp})^T Q_kc_k\| \notag\\
&=\|V_k^{\perp}(V_k^{\perp})^T\tilde{q}_k\|=\|(V_k^{\perp})^T\tilde{q}_k\|
=\sqrt{1-\varepsilon_k^2}
\label{qktilde}
\end{align}
with $\tilde{q}_k=Q_kc_k\in\mathcal{V}_k^R$ and $\|c_k\|=1$.
Since $\mathcal{V}_k$ is the orthogonal complement of $span\{V_k^{\perp}\}$,
by definition we know that $\tilde{q}_k\in \mathcal{V}_k^R$ has the largest acute
angle with $\mathcal{V}_k$, that is, it is the vector in $\mathcal{V}_k^R$
that contains the least information on $\mathcal{V}_k$.

Expand $\tilde{q}_k$ as the following orthogonal direct sum decomposition:
\begin{equation}\label{decompqk}
\tilde{q}_k=V_k^{\perp}(V_k^{\perp})^T\tilde{q}_k+V_kV_k^T\tilde{q}_k.
\end{equation}
Then from $\|\tilde{q}_k\|=1$ and \eqref{qktilde} we obtain
\begin{align}\label{angle2}
\|V_k^T\tilde{q}_k\|&=\|V_kV_k^T\tilde{q}_k\|=
\sqrt{1-\|V_k^{\perp}(V_k^{\perp})^T\tilde{q}_k\|^2}=\sqrt{1-(1-\varepsilon_k^2)}
=\varepsilon_k.
\end{align}
From \eqref{decompqk}, we next bound the Rayleigh quotient of $\tilde{q}_k$
with respect to $A^TA$ from below. By the SVD \eqref{eqsvd} of $A$ and
$V=(V_k,V_k^{\perp})$, we partition
$$
\Sigma=\left(\begin{array}{cc}
\Sigma_k &\\
&\Sigma_k^{\perp}
\end{array}
\right),
$$
where $\Sigma_k={\rm diag}(\sigma_1,\sigma_2,\ldots,\sigma_k)$ and
$\Sigma_k^{\perp}={\rm diag}(\sigma_{k+1},\sigma_{k+2},\ldots,\sigma_n)$.
Making use of $A^TAV_k=V_k\Sigma_k^2$ and $A^TAV_k^{\perp}=
V_k^{\perp}(\Sigma_k^{\perp})^2$ as well as $V_k^TV_k^{\perp}=\mathbf{0}$,
we obtain
\begin{align}
\tilde{q}_k^TA^TA\tilde{q}_k&=\left(V_k^{\perp}(V_k^{\perp})^T\tilde{q}_k+V_kV_k^T
\tilde{q}_k\right)^TA^TA \left(V_k^{\perp}(V_k^{\perp})^T\tilde{q}_k+
V_kV_k^T\tilde{q}_k\right) \notag\\
&=\left(\tilde{q}_k^TV_k^{\perp}(V_k^{\perp})^T+\tilde{q}_k^TV_kV_k^T\right)
\left(V_k^{\perp}(\Sigma_k^{\perp})^2(V_k^{\perp})^T\tilde{q}_k+V_k\Sigma_k^2V_k^T
\tilde{q}_k\right) \notag\\
&=\tilde{q}_k^TV_k^{\perp}(\Sigma_k^{\perp})^2(V_k^{\perp})^T\tilde{q}_k
+\tilde{q}_k^TV_k\Sigma_k^2V_k^T\tilde{q}_k. \label{expansion}
\end{align}
Observe that it is impossible for $(V_k^{\perp})^T\tilde{q}_k$ and
$V_k^T\tilde{q}_k$ to be the eigenvectors of $(\Sigma_k^{\perp})^2$
and $\Sigma_k^2$ associated with their respective smallest eigenvalues
$\sigma_n^2$ and $\sigma_k^2$ simultaneously, which are
the $(n-k)$-th canonical vector $e_{n-k}$ of $\mathbb{R}^{n-k}$ and
the $k$-th canonical vector $e_k$ of $\mathbb{R}^{k}$, respectively;
otherwise, we have $\tilde{q}_k=v_n$ and
$\tilde{q}_k=v_k$ simultaneously, which are impossible as $k<n$. Therefore,
from \eqref{expansion}, \eqref{qktilde} and \eqref{angle2},
we obtain the strict inequality
\begin{align*}
\tilde{q}_k^TA^TA\tilde{q}_k&> \|(V_k^{\perp})^T\tilde{q}_k\|^2
\sigma_n^2+\|V_k^T\tilde{q}_k\|^2\sigma_k^2
=(1-\varepsilon_k^2)\sigma_n^2+\varepsilon_k^2 \sigma_k^2,
\end{align*}
from which it follows that the lower bound of \eqref{rqi} holds. Similarly,
from \eqref{expansion}  and \eqref{qktilde}, \eqref{angle2}
we obtain the upper bound of \eqref{rqi}:
$$
\tilde{q}_k^TA^TA\tilde{q}_k <\|(V_k^{\perp})^T\tilde{q}_k\|^2
\|(\Sigma_k^{\perp})^2\|+\|V_k^T\tilde{q}_k\|^2\|\Sigma_k^2\|
=(1-\varepsilon_k^2)\sigma_{k+1}^2+\varepsilon_k^2 \sigma_1^2.
$$

From the lower bound of \eqref{rqi}, we see that if
$\varepsilon_k$ satisfies $\varepsilon_k^2 \sigma_k^2\geq \sigma_{k+1}^2$,
i.e., $\varepsilon_k\geq \frac{\sigma_{k+1}}{\sigma_k}$,
then $\sqrt{\tilde{q}_k^TA^TA\tilde{q}_k}>\sigma_{k+1}$, i.e.,
\eqref{est1} holds.

From \eqref{Bk}, we obtain $B_k^TB_k=Q_k^TA^TAQ_k$.
Note that $(\theta_k^{(k)})^2$ is the smallest eigenvalue
of the symmetric positive definite matrix $B_k^TB_k$.
Therefore, we have
\begin{equation}\label{rqi2}
(\theta_k^{(k)})^2=\min_{\|c\|=1} c^TQ_k^TA^TAQ_kc=
\min_{q\in \mathcal{V}_k^R,\ \|q\|=1} q^TA^TAq
=\hat{q}_k^TA^TA\hat{q}_k,
\end{equation}
where $\hat{q}_k$ is, in fact, the Ritz vector of $A^TA$ from
$\mathcal{V}_k^R$ corresponding to the smallest Ritz value $(\theta_k^{(k)})^2$.
Therefore, for $\tilde{q}_k$ defined in Theorem~\ref{initial} we have
$$
\theta_k^{(k)}\leq \sqrt{\tilde{q}_k^TA^TA\tilde{q}_k},
$$
from which it follows from \eqref{rqi} that
$(\theta_k^{(k)})^2<(1-\varepsilon_k^2)\sigma_{k+1}^2+\varepsilon_k^2 \sigma_1^2$.
As a result, for any $\delta>0$, we can choose $\varepsilon_k\geq 0$ such that
$$
(\theta_k^{(k)})^2<(1-\varepsilon_k^2)\sigma_{k+1}^2+\varepsilon_k^2 \sigma_1^2
\leq (1+\delta)\sigma_{k+1}^2,
$$
i.e., \eqref{thetasigma} holds,
solving which for $\varepsilon_k^2$ gives $\varepsilon_k^2\leq\frac{\delta}
{(\frac{\sigma_1}{\sigma_{k+1}})^2-1}$.
\qed

\begin{remark}
We analyze $\theta_k^{(k)}$ when $\varepsilon_k\geq \frac{\sigma_{k+1}}{\sigma_k}$.
In the sense of $\min$
in \eqref{rqi2}, $\hat{q}_k\in \mathcal{V}_k^R$ is the optimal vector that
extracts the least information from $\mathcal{V}_k$ and the richest information
from $span\{V_k^{\perp}\}$.
From Theorem~\ref{initial}, since $\mathcal{V}_k$ is the orthogonal complement
of $span\{V_k^{\perp}\}$, we know that $\tilde{q}_k\in \mathcal{V}_k^R$
has the largest acute angle with $\mathcal{V}_k$, that is,
it contains the least information from $\mathcal{V}_k$ and the richest information
from $span\{V_k^{\perp}\}$. Therefore, $\hat{q}_k$ and $\tilde{q}_k$ have a similar
optimality, so that we have
\begin{equation}\label{approxeq}
\theta_k^{(k)}\approx \sqrt{\tilde{q}_k^TA^TA\tilde{q}_k}.
\end{equation}
Combining this estimate with \eqref{est1}, we may have
$\theta_k^{(k)}>\sigma_{k+1}$ when
$\varepsilon_k\geq \frac{\sigma_{k+1}}{\sigma_k}$.
\end{remark}

\begin{remark}
We inspect the condition $\varepsilon_k\geq\frac{\sigma_{k+1}}
{\sigma_k}$ for \eqref{est1} and get insight into whether or not the true
$\varepsilon_k$ resulting from the three kinds of ill-posed problems satisfies
it. For severely ill-posed problems, the lower bound $\frac{\sigma_{k+1}}
{\sigma_k}$ is basically $\rho^{-1}$; for moderately ill-posed
problems with $\alpha>1$, the bound increases with increasing
$k\leq k_0$, and it cannot be close to one provided that $\alpha>1$ suitably
or $k_0$ not big; for mildly ill-posed problems with $\alpha<1$, the bound
increases faster than it does for moderately ill-posed problems, and
it may well approach one for $k\leq k_0$. Therefore,
the condition for \eqref{est1} requires that
$\|\sin\Theta(\mathcal{V}_k,\mathcal{V}_k^R)\|$
be not close to one for severely and moderately ill-posed problems,
but $\|\sin\Theta(\mathcal{V}_k,\mathcal{V}_k^R)\|$
must be close to zero for mildly ill-posed problems.
In view of \eqref{deltabound} and
$\|\sin\Theta(\mathcal{V}_k,\mathcal{V}_k^R)\|^2=1-\varepsilon_k^2$,
we have $\|\Delta_k\|^2=\frac{1-\varepsilon_k^2}{\varepsilon_k^2}$.
Thus, the condition $\varepsilon_k\geq\frac{\sigma_{k+1}}
{\sigma_k}$ for \eqref{est1} amounts to requiring
that $\|\Delta_k\|$ be at most modest and cannot
be large for severely and moderately ill-posed problems but
it must be fairly small for mildly ill-posed problems. Unfortunately,
Theorems~\ref{thm2}--\ref{moderate} and the remarks followed
indicate that $\|\Delta_k\|$ increases with $k$ increasing and is generally
large for a mildly ill-posed problem, while it
increases slowly with $k\leq k_0$ for a moderately
ill-posed problem with $\alpha>1$ suitably, and  by \eqref{case3}
it is approximately $\rho^{-(2+\beta)}$, considerably smaller
than one for a
severely ill-posed problem with $\rho>1$ not close to one.
Consequently, for mildly ill-posed problems,
because the actual $\|\Delta_k\|$ can hardly be small and is generally
large, the true $\varepsilon_k$ is small and may well be close to zero,
so that the condition $\varepsilon_k\geq\frac{\sigma_{k+1}}{\sigma_k}$
generally fails to meet as $k$ increases, while it is satisfied for severely
or moderately ill-posed problems with $\rho>1$ or $\alpha>1$ suitably.
\end{remark}

\begin{remark}\label{appear}
\eqref{thetasigma} shows that there is
at least one $\theta_k^{(k)}\leq\sigma_{k+1}$ if
$\|\sin\Theta(\mathcal{V}_k,\mathcal{V}_k^R)\|$ is sufficiently close
to one since we can choose $\delta$ small enough such that
$(1+\delta)^{1/2}\sigma_{k+1}$ is close to $\sigma_{k+1}$ arbitrarily.
As we have shown, $\|\sin\Theta(\mathcal{V}_k,\mathcal{V}_k^R)\|$ cannot be
close to one for severely or moderately ill-posed problems with $\rho>1$ or
$\alpha>1$ suitably, but it is generally so for mildly ill-posed problems.
This means that for some $k\leq k_0$ it is very likely that
$\theta_k^{(k)}\leq\sigma_{k+1}$ for mildly ill-posed problems.
\end{remark}

We must be aware that our above analysis on
$\theta_k^{(k)}>\sigma_{k+1}$ is not rigorous because
we cannot quantify {\em how small}
$\sqrt{\tilde{q}_k^TA^TA\tilde{q}_k}-\theta_k^{(k)}$ is.
From $\theta_k^{(k)}\leq\sqrt{\tilde{q}_k^TA^TA\tilde{q}_k}$, it is apparent
that the condition $\varepsilon_k\geq\frac{\sigma_{k+1}}{\sigma_k}$
may not be sufficient for $\theta_k^{(k)}>\sigma_{k+1}$.
We delay our detailed and rigorous analysis
to Section \ref{rankapp}.

Theorems~\ref{thm2}--\ref{moderate} establish necessary
background for answering the fundamental
concern by Bj\"{o}rck and Eld\'{e}n, and their proof approaches
also provide key ingredients for some of the later results.
We now present the following results,
which will play a central role in our later analysis.

\begin{theorem}\label{thm3}
Assume that the discrete Picard condition \eqref{picard} is satisfied,
let $\Delta_k\in \mathbb{R}^{(n-k)\times k}$ be defined as
\eqref{defdelta} and $L_j^{(k)}(0)$ and $L_{k_1}^{(k)}(0)$
defined as \eqref{lk}, and write
$\Delta_k=(\delta_1,\delta_2,\ldots,\delta_k)$.  Then for severely
ill-posed problems and $k=1,2,\ldots,n-1$ we have
\begin{align}
\|\delta_j\|&\leq \frac{\sigma_{k+1}}{\sigma_j}\frac{| u_{k+1}^Tb|}
{| u_j^T b|}\left(1+\mathcal{O}(\rho^{-2})\right)| L_j^{(k)}(0)|,
\ k>1,\ j=1,2,\ldots,k, \label{columndelta} \\
\|\delta_1\|& \leq \frac{\sigma_{2}}{\sigma_1}\frac{| u_2^Tb|}
{| u_1^T b|}\left(1+\mathcal{O}(\rho^{-2})\right)|,\ k=1
\label{columndelta1}
\end{align}
and
\begin{equation} \label{prodnorm}
\|\Sigma_k\Delta_k^T\|\leq \left\{\begin{array}{ll}
\sigma_{k+1}\frac{| u_{k+1}^Tb|}{| u_k^T b|}
\left(1+\mathcal{O}(\rho^{-2})\right)
& \mbox{ for } 1\leq k\leq k_0,\\
\sigma_{k+1}\sqrt{k-k_0+1}\left(1+\mathcal{O}(\rho^{-2})\right)
& \mbox{ for } k_0<k\leq n-1;
\end{array}
\right.
\end{equation}
for moderately or mild ill-posed problems with the singular values
$\sigma_j=\zeta j^{-\alpha}$ and $\zeta$ a positive constant we have
\begin{align}
\|\delta_j\|&\leq \frac{\sigma_k}{\sigma_j}\frac{| u_{k+1}^Tb|}
{| u_j^T b|} \sqrt{\frac{k}{2\alpha-1}}|L_j^{(k)}(0)|,
\ k>1,\ j=1,2,\ldots,k, \label{columndelta2} \\
\|\delta_1\|&\leq \frac{| u_2^Tb|}
{| u_1^T b|} \sqrt{\frac{1}{2\alpha-1}},\ k=1 \label{columnnorm}
\end{align}
and
\begin{equation}\label{prodnorm2}
\|\Sigma_k\Delta_k^T\|\leq \left\{\begin{array}{ll}
\sigma_1\frac{| u_2^Tb|}{| u_1^T b|}\sqrt{\frac{1}{2\alpha-1}} &
 \mbox{ for } k=1,\\
\sigma_k\frac{| u_{k+1}^Tb|}{| u_k^T b|}\sqrt{\frac{k^2}{4\alpha^2-1}+
\frac{k}{2\alpha-1}}
|L_{k_1}^{(k)}(0)|& \mbox{ for } 1<k\leq k_0, \\
\sigma_k \sqrt{\frac{k k_0}{4\alpha^2-1}+
\frac{k(k-k_0+1)}{2\alpha-1}}|L_{k_1}^{(k)}(0)|
& \mbox{ for } k_0<k\leq n-1.
\end{array}
\right.
\end{equation}
\end{theorem}

{\em Proof}.
From \eqref{defdelta} and \eqref{amplify}, for $j=1,2,\ldots,k$ and $k>1$
we have
\begin{equation}\label{deltaj}
\|\delta_j\|^2\leq |L_j^{(k)}(0)|^2 \sum_{i=k+1}^n\frac{\sigma_{i}^2}
{\sigma_j^2}\frac{| u_i^T b|^2}{| u_j^T b|^2}
\end{equation}
and from \eqref{deltaexp}, for $k=1$ we have
\begin{equation}\label{deltas}
\|\delta_1\|^2=\sum_{i=2}^n\frac{\sigma_{i}^2}
{\sigma_1^2}\frac{| u_i^T b|^2}{| u_1^T b|^2}.
\end{equation}
For severely ill-posed problems, $k=1,2,\ldots,n-1$ and
$j=1,2,\ldots,k$, from \eqref{severe1} we obtain
\begin{align*}
\sum_{i=k+1}^n\frac{\sigma_{i}^2}
{\sigma_j^2}\frac{| u_i^T b|^2}{| u_j^T b|^2}&=
\frac{1}{\sigma_j^2| u_j^Tb|^2}
\sum_{i=k+1}^n\sigma_{i}^2 |u_i^T b|^2 \\
&\leq \frac{\sigma_{k+1}^2}
{\sigma_j^2}\frac{| u_{k+1}^T b|^2}{| u_j^T b|^2}
\left(1+\mathcal{O}(\rho^{-2})\right).
\end{align*}
For moderately or mildly ill-posed problems, $k=1,2,\ldots,n-1$ and
$j=1,2,\ldots,k$, from \eqref{modeest} we obtain
\begin{align*}
\sum_{i=k+1}^n\frac{\sigma_{i}^2}
{\sigma_j^2}\frac{| u_i^T b|^2}{| u_j^T b|^2}
&=\frac{1}{\sigma_j^2| u_j^Tb|^2}
\sum_{i=k+1}^n\sigma_{i}^2| u_i^T b|^2\\
&\leq \frac{\sigma_k^2}
{\sigma_j^2}\frac{| u_{k+1}^T b|^2}{| u_j^T b|^2}
\frac{k}{2\alpha-1}.
\end{align*}
Combining the above with \eqref{deltaj}, \eqref{lkk} and
$
\left(1+\mathcal{O}(\rho^{-2})\right)|L_{k_1}^{(k)}(0)|= 1+
\mathcal{O}(\rho^{-2}), \ k=2,3,\ldots,n-1,
$
we obtain \eqref{columndelta}, while \eqref{columndelta2} follows
from the above and \eqref{deltaj} directly. For $k=1$, from \eqref{deltas}
and the above we get \eqref{columndelta1} and \eqref{columnnorm}, respectively.

By \eqref{delta1}, for $k>1$ we have
$$
|\Delta_k\Sigma_k|\leq |L_{k_1}^{(k)}(0)|\left|(\sigma_{k+1} u_{k+1}^T b,
\sigma_{k+2}u_{k+2}^Tb,\ldots,\sigma_n u_n^T b)^T
\left(\frac{1}{u_1^Tb},\frac{1}{u_2^Tb},\ldots,
\frac{1}{u_k^Tb}\right)\right|.
$$
Therefore, we get
\begin{align}
\|\Sigma_k\Delta_k^T\|&=\|\Delta_k\Sigma_k\|\leq \left\||\Delta_k\Sigma_k|\right\|
\notag\\
&\leq |L_{k_1}^{(k)}(0)|\left(\sum_{j=k+1}^n\sigma_j^2| u_j^Tb|^2\right)^{1/2}
\left(\sum_{j=1}^k \frac{1}{| u_j^Tb|^2}\right)^{1/2}. \label{sigdel}
\end{align}
By \eqref{deltaexp}, for $k=1$ we have
$$
\|\Delta_1\Sigma_1\|=\left(\sum_{j=2}^n\sigma_j^2| u_j^Tb|^2\right)^{1/2}
\frac{1}{| u_1^Tb|}.
$$
We have derived the bounds \eqref{severe1} and
\eqref{modeest} for $\left(\sum_{j=k+1}^n\sigma_j^2| u_j^Tb|^2\right)^{1/2}$
for severely and moderately or mildly ill-posed problems, respectively, from
which we obtain \eqref{prodnorm} and \eqref{prodnorm2} for $k=1$. In order
to bound $\|\Sigma_k\Delta_k^T\|$ for $k>1$, we need to estimate
$\left(\sum_{j=1}^k\frac{1}{| u_j^Tb|^2}\right)^{1/2}$.
We next carry out this task for severely and moderately or mildly ill-posed
problems, respectively, for each kind of which we consider the cases of
$k\leq k_0$ and $k>k_0$ separately.

Case of $k\leq k_0$ for severely ill-posed problems: From the discrete
Picard condition \eqref{picard} and \eqref{ideal3}, we obtain
\begin{align*}
\sum_{j=1}^k
\frac{1}{| u_j^Tb|^2}&= \frac{1}{| u_k^Tb|^2} \sum_{j=1}^k
\frac{| u_k^Tb|^2}{| u_j^Tb|^2}
=\frac{1}{| u_k^Tb|^2} \left(1+\mathcal{O}\left(\sum_{j=1}^{k-1}\rho^{2(j-k)
(1+\beta)}\right)\right)\\
&=\frac{1}{| u_k^Tb|^2} \left(1+\mathcal{O}(\rho^{-2(1+\beta)}) \right).
\end{align*}

Case of $k> k_0$ for severely ill-posed problems: From \eqref{ideal3} and
\eqref{ideal2}, we obtain
\begin{align*}
\sum_{j=1}^k
\frac{1}{| u_j^Tb|^2}&= \frac{1}{| u_k^Tb|^2}\left( \sum_{j=1}^{k_0}
\frac{| u_k^Tb|^2}{| u_j^Tb|^2}+\sum_{j=k_0+1}^{k}
\frac{| u_k^Tb|^2}{| u_j^Tb|^2}\right)\\
&=\frac{1}{| u_k^Tb|^2}\left(1+\mathcal{O}\left(\sum_{j=1}^{k_0-1}
\rho^{2(j-k_0)(1+\beta)}\right)+k-k_0\right)\\
&=\frac{1}{| u_k^Tb|^2}\left(1+\mathcal{O}(\rho^{-2(1+\beta)})+k-k_0\right).
\end{align*}
Substituting the above two relations for the two cases into \eqref{sigdel} and
combining them with \eqref{severe1} and \eqref{lkk}, we
get \eqref{prodnorm}.

Case of $k\leq k_0$ for moderately or mildly ill-posed problems: From \eqref{ideal3}
we have
\begin{align*}
\sum_{j=1}^k
\frac{1}{| u_j^Tb|^2}&= \frac{1}{| u_k^Tb|^2} \sum_{j=1}^k
\frac{| u_k^Tb|^2}{| u_j^Tb|^2}= \frac{1}{| u_k^Tb|^2}\sum_{j=1}^k
\left(\frac{j}{k}\right)^{2\alpha (1+\beta)}\\
&<\frac{1}{| u_k^Tb|^2} \sum_{j=1}^k
\left(\frac{j}{k}\right)^{2\alpha}
=\frac{1}{| u_k^T b|^2} k\sum_{j=1}^k \frac{1}{k}\left(\frac{j}{k}
\right)^{2\alpha} \notag \\
&< \frac{1}{| u_k^T b|^2} \left(k \int_0^1
x^{2\alpha}dx+1 \right)=\frac{1}{| u_k^Tb|^2}\left(\frac{k}{2\alpha+1}+1\right).
\end{align*}

Case of $k> k_0$ for moderately or mildly ill-posed problems: From \eqref{ideal3} and
\eqref{ideal2} we have
\begin{align*}
\sum_{j=1}^k\frac{1}{| u_j^Tb|^2}
&= \frac{1}{| u_k^Tb|^2} \left(\sum_{j=1}^{k_0}
\frac{| u_k^Tb|^2}{| u_j^Tb|^2}
+\sum_{j=k_0+1}^{k}
\frac{| u_k^Tb|^2}{| u_j^Tb|^2}\right)\\
&= \frac{1}{| u_k^Tb|^2} \left(\sum_{j=1}^{k_0}
\left(\frac{j}{k_0}\right)^{2\alpha (1+\beta)}+k-k_0\right)\\
&<\frac{1}{| u_k^Tb|^2} \left(\sum_{j=1}^{k_0}
\left(\frac{j}{k_0}\right)^{2\alpha}+k-k_0\right)\\
&\leq \frac{1}{| u_k^Tb|^2} \left(\frac{k_0}{2\alpha+1}+1+k-k_0\right).
\end{align*}
Substituting the above two bounds for the two cases into \eqref{sigdel} and
combining them with \eqref{modeest}, we get \eqref{prodnorm2}.
\qed

\eqref{prodnorm} and \eqref{prodnorm2} indicate that $\|\Sigma_k\Delta_k^T\|$
decays swiftly as $k$ increases. As has been seen, we must take some cares to
accurately bound $\|\Sigma_k\Delta_k^T\|$. Indeed,
for $1<k\leq k_0$, if we had simply estimated it by
\begin{equation}\label{rough}
\|\Sigma_k\Delta_k^T\|\leq \|\Sigma_k\|\|\Delta_k^T\|=\sigma_1\|\Delta_k\|,
\end{equation}
we would have obtained a bound, which not only does not decay but also
increases for moderately and mildly ill-posed problems as $k$ increases.
Such bound is useless to precisely analyze the regularization of
LSQR for ill-posed problems and makes us
impossible to get those predictively accurate results to be presented
in Sections \ref{rankapp}--\ref{alphabeta}.

\section{The rank $k$ approximation $P_{k+1}B_kQ_k^T$ to $A$, the Ritz values
$\theta_i^{(k)}$ and the regularization of LSQR}\label{rankapp}

Making use of Theorems~\ref{thm2}--\ref{thm3}, we are able to solve those
key problems stated before Theorem~\ref{thm2} and give definitive
answers to the fundamental concern by Bj\"{o}rck and Eld\'{e}n, proving
that LSQR has the full regularization
for severely or moderately ill-posed problems with $\rho>1$ or
$\alpha>1$ suitably and it, in general, has only the partial regularization
for mildly ill-posed problems.

Define
\begin{equation}\label{gammak}
\gamma_k = \|A-P_{k+1}B_kQ_k^T\|,
\end{equation}
which measures the accuracy of the rank $k$ approximation $P_{k+1}B_kQ_k^T$ to $A$
generated by Lanczos bidiagonalization. Recall \eqref{xk} and
the comments followed. It is known that the full
or partial regularization of LSQR uniquely depends on whether or not
$\gamma_k\approx \sigma_{k+1}$ holds, where we will make the precise meaning
`$\approx$' clear by introducing the definition of near best
rank $k$ approximation to $A$, and on whether or not the $k$ Ritz values
$\theta_i^{(k)}$ approximate
the $k$ large singular values $\sigma_i$ of $A$ in natural order
for $k=1,2,\ldots, k_0$. If both of them hold, LSQR has the full regularization;
if either of them is not satisfied, LSQR has only the partial
regularization.


\subsection{Accuracy of the rank $k$ approximation
$P_{k+1}B_kQ_k^T$ to $A$}\label{rankaccur}

We first present one of the main results in this paper.

\begin{theorem}\label{main1}
Assume that the discrete Picard condition \eqref{picard} is
satisfied. Then for $k=1,2,\ldots,n-1$ we have
\begin{equation}\label{final}
  \sigma_{k+1}\leq \gamma_k\leq \sqrt{1+\eta_k^2}\sigma_{k+1}
\end{equation}
with
\begin{equation} \label{const1}
\eta_k\leq \left\{\begin{array}{ll}
\xi_k\frac{| u_{k+1}^Tb|}{| u_k^T b|}
\left(1+\mathcal{O}(\rho^{-2})\right)
& \mbox{ for } 1\leq k\leq k_0,\\
\xi_k\sqrt{k-k_0+1}\left(1+\mathcal{O}(\rho^{-2})\right)
& \mbox{ for } k_0<k \leq n-1
\end{array}
\right.
\end{equation}
for severely ill-posed problems and
\begin{equation}\label{const2}
\eta_k\leq \left\{\begin{array}{ll}
\xi_1\frac{\sigma_1}{\sigma_2}\frac{| u_2^Tb|}{| u_1^Tb|}
\sqrt{\frac{1}{2\alpha-1}} & \mbox{ for } k=1, \\
\xi_k\frac{\sigma_k}{\sigma_{k+1}}\frac{|u_{k+1}^T b|}{|u_k^T b|}
\sqrt{\frac{k^2}{4\alpha^2-1}+\frac{k}{2\alpha-1}}
|L_{k_1}^{(k)}(0)|& \mbox{ for } 1< k\leq k_0, \\
\xi_k\frac{\sigma_k}{\sigma_{k+1}}\sqrt{\frac{k k_0}{4\alpha^2-1}+
\frac{k(k-k_0+1)}{2\alpha-1}}|L_{k_1}^{(k)}(0)|
& \mbox{ for } k_0<k\leq n-1
\end{array}
\right.
\end{equation}
for moderately or mildly ill-posed problems with $\sigma_j=\zeta j^{-\alpha},\
j=1,2,\ldots,n$, where
$\xi_k=\sqrt{\left(\frac{\|\Delta_k\|}{1+\|\Delta_k\|^2}\right)^2+1}$ for
$\|\Delta_k\|<1$ and $\xi_k\leq\frac{\sqrt{5}}{2}$ for $\|\Delta_k\|\geq 1$.
\end{theorem}

{\em Proof}.
Since $A_k$ is the best rank $k$ approximation to
$A$ with respect to the 2-norm and $\|A-A_k\|=\sigma_{k+1}$,
the lower bound in \eqref{final} holds. Next we prove the upper bound.

From \eqref{eqmform1}, we obtain
\begin{align}
\gamma_k
&= \|A-P_{k+1}B_kQ_k^T\|= \|A-AQ_kQ_k^T\|= \|A(I-Q_kQ_k^T)\|. \label{gamma2}
\end{align}
From Algorithm~\ref{alg:lb}, \eqref{kry}, \eqref{zk} and \eqref{decomp}, we obtain
$$
\mathcal{V}_k^R
=\mathcal{K}_{k}(A^{T}A,A^{T}b)=span\{Q_k\}=span\{\hat{Z}_k\}
$$
with $Q_k$ and $\hat{Z}_k$ being orthonormal, and the
orthogonal projector onto $\mathcal{V}_k^R$ is thus
\begin{equation}\label{twobasis}
Q_kQ_k^T=\hat{Z}_k\hat{Z}_k^T.
\end{equation}
Keep in mind that $A_k=U_k\Sigma_k V_k^T$. It is direct to justify that
$(U_k\Sigma_k V_k^T)^T(A-U_k\Sigma_k V_k^T)=\mathbf{0}$ for $k=1,2,\ldots,n-1$.
Therefore, exploiting this and noting that $\|I-\hat{Z}_k\hat{Z}_k^T\|=1$ and
$V_k^TV_k^{\perp}=\mathbf{0}$  for $k=1,2,\ldots,n-1$,
we get from \eqref{gamma2}, \eqref{twobasis} and \eqref{decomp} that
\begin{align}
\gamma_k^2 &= \|(A-U_k\Sigma_kV_k^T+U_k\Sigma_kV_k^T)(I-\hat{Z}_k\hat{Z}_k^T)\|^2
  \notag\\
  &=\max_{\|y\|=1}\|(A-U_k\Sigma_kV_k^T+U_k\Sigma_kV_k^T)
  (I-\hat{Z}_k\hat{Z}_k^T)y\|^2 \notag\\
  &=\max_{\|y\|=1}\|(A-U_k\Sigma_kV_k^T)(I-\hat{Z}_k\hat{Z}_k^T)y+
 U_k\Sigma_kV_k^T(I-\hat{Z}_k\hat{Z}_k^T)y\|^2\notag\\
 &=\max_{\|y\|=1}\left(\|(A-U_k\Sigma_kV_k^T)(I-\hat{Z}_k\hat{Z}_k^T)y\|^2+
 \| U_k\Sigma_kV_k^T(I-\hat{Z}_k\hat{Z}_k^T)y\|^2\right)\notag\\
 &\leq \|(A-U_k\Sigma_kV_k^T)(I-\hat{Z}_k\hat{Z}_k^T)\|^2+
 \| U_k\Sigma_kV_k^T(I-\hat{Z}_k\hat{Z}_k^T)\|^2 \notag \\
 &\leq \sigma_{k+1}^2+\| \Sigma_kV_k^T(I-\hat{Z}_k\hat{Z}_k^T)\|^2 \notag\\
 &\leq \sigma_{k+1}^2+\|\Sigma_kV_k^T\left(I-(V_k+V_k^{\perp}\Delta_k)(I+
 \Delta_k^T\Delta_k)^{-1}(V_k+V_k^{\perp}\Delta_k)^T\right)\|^2\notag\\
 &= \sigma_{k+1}^2 + \left\|\Sigma_k\left(V_k^T-(I+
 \Delta_k^T\Delta_k)^{-1}(V_k+V_k^{\perp}\Delta_k)^T\right)\right\|^2 \notag\\
  &= \sigma_{k+1}^2 + \left\|\Sigma_k(I+
 \Delta_k^T\Delta_k)^{-1}\left((I+\Delta_k^T\Delta_k)V_k^T-
 \left(V_k+V_k^{\perp}\Delta_k\right)^T\right)\right\|^2 \notag\\
  &= \sigma_{k+1}^2+ \|\Sigma_k(I+
 \Delta_k^T\Delta_k)^{-1}\left(\Delta_k^T\Delta_kV_k^T-\Delta_k^T
 (V_k^{\perp})^T\right)\|^2\notag\\
  &= \sigma_{k+1}^2 + \|\Sigma_k(I+
 \Delta_k^T\Delta_k)^{-1}\Delta_k^T\Delta_kV_k^T-\Sigma_k(I+
 \Delta_k^T\Delta_k)^{-1}\Delta_k^T(V_k^{\perp})^T\|^2 \label{twomatrix}\\
    &\leq \sigma_{k+1}^2 +\|\Sigma_k(I+
 \Delta_k^T\Delta_k)^{-1}\Delta_k^T\Delta_k\|^2+\|\Sigma_k(I+
 \Delta_k^T\Delta_k)^{-1}\Delta_k^T\|^2\notag\\
 &=\sigma_{k+1}^2+\epsilon_k^2,
 \label{estimate1}
\end{align}
where the last inequality follows by using $V_k^T V_k^{\perp}=\mathbf{0}$ and
the definition of the induced matrix 2-norm to amplify the second term
in \eqref{twomatrix}.

We estimate $\epsilon_k$ accurately below. To this end, we need to use two
key identities and some results related. By the SVD of $\Delta_k$, it is direct
to justify that
\begin{equation}\label{inden1}
(I+
 \Delta_k^T\Delta_k)^{-1}\Delta_k^T\Delta_k=\Delta_k^T\Delta_k(I+
 \Delta_k^T\Delta_k)^{-1}
\end{equation}
and
\begin{equation}\label{inden2}
(I+
 \Delta_k^T\Delta_k)^{-1}\Delta_k^T=\Delta_k^T(I+
 \Delta_k\Delta_k^T)^{-1}.
\end{equation}
Define the function $f(\lambda)=\frac{\lambda}{1+\lambda^2}$ with
$\lambda\in [0,\infty)$. Since the derivative
$f^{\prime}(\lambda)=\frac{1-\lambda^2}{(1+\lambda^2)^2}$,
$f(\lambda)$ is monotonically increasing for $\lambda\in [0,1]$
and decreasing for $\lambda\in [1,\infty)$, and the maximum of
$f(\lambda)$ over $\lambda\in [0,\infty)$ is $\frac{1}{2}$, which attains at
$\lambda=1$. Based on
these properties and exploiting the SVD of $\Delta_k$, for the matrix 2-norm
we get
\begin{equation}\label{compact}
\|\Delta_k(I+\Delta_k^T\Delta_k)^{-1}\|=\frac{\|\Delta_k\|}{1+\|\Delta_k\|^2}
\end{equation}
for $\|\Delta_k\|<1$ and
\begin{equation}\label{noncomp}
\|\Delta_k(I+\Delta_k^T\Delta_k)^{-1}\|\leq\frac{1}{2}
\end{equation}
for $\|\Delta_k\|\geq 1$ (Note: in this case, since $\Delta_k$ may have at least
one singular value smaller than one, we do not
have an expression like \eqref{compact}). It then follows
from \eqref{estimate1}, \eqref{compact}, \eqref{noncomp}
and $\|(1+\Delta_k\Delta_k^T)^{-1}\|\leq 1$ that
\begin{align}
\epsilon_k^2&=\|\Sigma_k \Delta_k^T\Delta_k(I+
 \Delta_k^T\Delta_k)^{-1}\|^2+\|\Sigma_k \Delta_k^T (I+
 \Delta_k\Delta_k^T)^{-1}\|^2 \label{separa}\\
 &\leq \|\Sigma_k\Delta_k^T\|^2\|\Delta_k(I+\Delta_k^T\Delta_k)^{-1}\|^2+
 \|\Sigma_k\Delta_k^T\|^2 \|(1+\Delta_k\Delta_k^T)^{-1}\|^2 \notag\\
 &\leq \|\Sigma_k\Delta_k^T\|^2\left(\|\Delta_k
 (I+\Delta_k^T\Delta_k)^{-1}\|^2+1\right)\notag\\
  &=\|\Sigma_k\Delta_k^T\|^2\left(\left(\frac{\|\Delta_k\|}
  {1+\|\Delta_k\|^2}\right)^2+1\right)
  =\xi_k^2\|\Sigma_k\Delta_k^T\|^2 \notag
\end{align}
for $\|\Delta_k\|<1$ and
$$
\epsilon_k\leq \|\Sigma_k\Delta_k^T\|\sqrt{\|\Delta_k
(I+\Delta_k^T\Delta_k)^{-1}\|^2+1}=\xi_k\|\Sigma_k\Delta_k^T\|
\leq \frac{\sqrt{5}}{2}\|\Sigma_k\Delta_k^T\|
$$
for $\|\Delta_k\|\geq 1$. Replace $\|\Sigma_k\Delta_k^T\|$ by
its bounds \eqref{prodnorm} and \eqref{prodnorm2} in the
above, insert the resulting bounds for $\epsilon_k$
into \eqref{estimate1}, and let $\epsilon_k=\eta_k\sigma_{k+1}$.
Then we obtain the upper bound in \eqref{final} with $\eta_k$
satisfying \eqref{const1} and \eqref{const2} for severely and moderately
or mildly ill-posed problems, respectively.
\qed

Note from \eqref{ideal3} that
$$
\frac{|u_{k+1}^T b|}{| u_k^T b|}=
\frac{\sigma_{k+1}^{1+\beta}}{\sigma_k^{1+\beta}},\ k\leq k_0.
$$
Therefore, for the right-hand side of \eqref{const2} and $k\leq k_0$
we have
$$
\frac{\sigma_k}{\sigma_{k+1}}\frac{| u_{k+1}^T b|}{| u_k^T b|}=
\left(\frac{\sigma_{k+1}}{\sigma_k}\right)^{\beta}<1.
$$

\begin{remark}\label{decayrate}
For severely ill-posed problems, from \eqref{case3}, \eqref{case4} and
the definition of $\xi_k$ we know that
$$
\xi_k(1+\mathcal{O}(\rho^{-2}))
=1+\mathcal{O}(\rho^{-2})
$$
for both $k\leq k_0$ and $k>k_0$. Therefore, from
\eqref{const1} and \eqref{ideal3}, for $k\leq k_0$ we have
\begin{equation}\label{etak0}
\eta_k\leq \xi_k
\frac{|u_{k+1}^Tb |}{|u_k^T b |}\left(1+\mathcal{O}(\rho^{-2})\right)
=\frac{|u_{k+1}^Tb |}{|u_k^T b |}=
\frac{\sigma_{k+1}^{1+\beta}}{\sigma_k^{1+\beta}}=\mathcal{O}(\rho^{-1-\beta})<1
\end{equation}
by ignoring the smaller term $\mathcal{O}(\rho^{-1-\beta})
\mathcal{O}(\rho^{-2})
=\mathcal{O}(\rho^{-3-\beta})$, and for $k>k_0$ we have
\begin{equation}\label{incres}
\eta_k\leq \xi_k\sqrt{k-k_0+1}\left(1+\mathcal{O}(\rho^{-2})\right)
=\sqrt{k-k_0+1}
\end{equation}
by ignoring the smaller term $\sqrt{k-k_0+1}\mathcal{O}(\rho^{-2})$,
which increases slowly with $k$.
\end{remark}

\begin{remark}
For the moderately or mildly ill-posed problems with
$\sigma_j=\zeta j^{-\alpha}$,
from the derivation on $\eta_k$ and its estimate \eqref{const2},
for $k\leq k_0$ we approximately have
\begin{equation}\label{etadelta}
\frac{\sigma_k}{\sigma_{k+1}}\|\Delta_k\|\leq
\eta_k\leq \frac{\sqrt{5}}{2}\frac{\sigma_k}{\sigma_{k+1}}\|\Delta_k\|,
\end{equation}
and for $k>k_0$, from \eqref{lk1size} and \eqref{lk1sizemild}
we approximately have
\begin{align}
\eta_k&< \frac{\sigma_k}{\sigma_{k+1}}\sqrt{\frac{k k_0}{4\alpha^2-1}+
\frac{k(k-k_0+1)}{2\alpha-1}}|L_{k_1}^{(k)}(0)| \notag\\
&\sim \frac{k^{3/2}\sqrt{k_0}}{(2\alpha+1)\sqrt{{4\alpha^2-1}}}+
\frac{k^{3/2}\sqrt{k-k_0+1}}{(2\alpha+1)\sqrt{{2\alpha-1}}}, \label{asym}
\end{align}
which increases faster than the right-hand side of \eqref{incres} with
respect to $k$.
\end{remark}

\begin{remark}\label{decayrate2}
From \eqref{final}, \eqref{const1} and \eqref{etak0},
for severely ill-posed problems we have
$$
1<\sqrt{1+\eta_k^2}<1+\frac{1}{2}{\eta_k^2}\leq
1+\frac{1}{2}\frac{\sigma_{k+1}^{2(1+\beta)}}{\sigma_k^{2(1+\beta)}}
\sim 1+\frac{1}{2}\rho^{-2(1+\beta)},
$$
and $\gamma_k$ is an accurate approximation to
$\sigma_{k+1}$ for $k\leq k_0$ and marginally less accurate for $k>k_0$.
Thus, the rank $k$ approximation $P_{k+1}B_kQ_k^T$ is as accurate as
the best rank $k$ approximation $A_k$ within the
factor $\sqrt{1+\eta_k^2}\approx 1$ for $k\leq k_0$ and $\rho>1$ suitably.
For moderately ill-posed problems, $\gamma_k$ is still an excellent
approximation to $\sigma_{k+1}$, and the rank $k$ approximation
$P_{k+1}B_kQ_k^T$ is almost as accurate as the best rank $k$ approximation
$A_k$ for $k\leq k_0$. Therefore, $P_{k+1}B_kQ_k^T$ plays the same role as
$A_k$ for these two kinds of ill-posed problems and $k\leq k_0$, it is known
from the clarification in Section \ref{lsqr} that LSQR may have the full
regularization. We will, afterwards, deepen this theorem and derive
more results, proving that LSQR must have the full regularization for
these two kinds of problems provided that $\rho>1$ and $\alpha>1$ suitably.

For both severely and moderately ill-posed problems, we note that the
situation is not so satisfying for increasing $k>k_0$. But at that time,
a possibly big $\eta_k$ does not do harm to our regularization purpose
since we will prove that, provided that $\rho>1$ and $\alpha>1$ suitably,
LSQR has the full regularization and has already found
a best possible regularized solution at semi-convergence occurring at
iteration $k_0$. If it is the case, we will simply stop performing it
after semi-convergence.
\end{remark}

\begin{remark}\label{mildre}
For mildly ill-posed problems, the situation is fundamentally different.
As clarified in Remark~\ref{mildrem}, we have
$\sqrt{\frac{k^2}{4\alpha^2-1}+\frac{k}{2\alpha-1}}>1$ and
$|L_{k_1}^{(k)}(0)|>1$ considerably as $k$ increases up to $k_0$
because of $\frac{1}{2}<\alpha\leq 1$,
leading to $\eta_k>1$ substantially. This means that
$\gamma_{k_0}$ is substantially bigger than $\sigma_{k_0+1}$ and can
well lie between $\sigma_{k_0}$ and $\sigma_1$, so that
the rank $k_0$ approximation $P_{k_0+1}B_{k_0}Q_{k_0}^T$ is much less accurate
than the best rank $k_0$ approximation $A_{k_0}$ and LSQR has only the partial
regularization.
\end{remark}

\begin{remark}
There are several subtle treatments in the proof of Theorem~\ref{main1}, each of
which turns out to be absolutely necessary. Ignoring or missing any one of them
would be fatal and make us fail to obtain accurate estimates for $\epsilon_k$
defined by \eqref{estimate1}. The first is the treatment of
$\|U_k\Sigma_kV_k^T(I-\hat{Z}_k\hat{Z}_k^T)\|$.
By the definition of $\|\sin\Theta(\mathcal{V}_k,\mathcal{V}_k^R)\|$, if we had
amplified it by
$$
\|U_k\Sigma_kV_k^T(I-\hat{Z}_k\hat{Z}_k^T)\|
\leq \|\Sigma_k\|\|V_k^T(I-\hat{Z}_k\hat{Z}_k^T)\|=
\sigma_1\|\sin\Theta(\mathcal{V}_k,\mathcal{V}_k^R)\|,
$$
we would have obtained a too large overestimate, which is almost a fixed
constant for severely
ill-posed problems and $k=1,2,\ldots,k_0$ and increases with
$k=1,2,\ldots,k_0$ for moderately and mildly ill-posed problems. Such rough
estimates are useless to get a meaningful bound for $\gamma_k$.
The second is the use of
\eqref{inden1} and \eqref{inden2}. The third is the extraction of
$\|\Sigma_k\Delta_k^T\|$ from \eqref{separa} as a whole other than
amplify it to $\|\Sigma_k\|\|\Delta_k\|=\sigma_1\|\Delta_k\|$.
The fourth is accurate estimates for it;
see \eqref{prodnorm} and \eqref{prodnorm2} in Theorem~\ref{thm3}.
For example, without using \eqref{inden1} and \eqref{inden2}, by
\eqref{deltabound} we would have
no way but to obtain
\begin{align*}
\epsilon_k^2 &\leq \|\Sigma_k\|^2\|(I+
 \Delta_k^T\Delta_k)^{-1}\Delta_k^T\Delta_k\|^2+\|\Sigma_k\|^2\|(I+
 \Delta_k^T\Delta_k)^{-1}\Delta_k^T\|^2\\
 &=\sigma_1^2\left(\frac{\|\Delta_k\|^2}{1+\|\Delta_k\|^2}\right)^2+\sigma_1^2
 \|(I+\Delta_k^T\Delta_k)^{-1}\Delta_k^T\|^2\\
 &=\sigma_1^2\|\sin\Theta(\mathcal{V}_k,\mathcal{V}_k^R)\|^4+\sigma_1^2
 \|\Delta_k(I+\Delta_k^T\Delta_k)^{-1}\|^2.
\end{align*}
From \eqref{noncomp} and the previous
estimates for $\|\Delta_k\|$, such bound is too pessimistic
and completely useless in our context, and
it even does not decrease and could not be small as $k$ increases, while
our estimates for $\epsilon_k=\eta_k\sigma_{k+1}$
in Theorem~\ref{main1} are much more accurate and
decay swiftly as $k$ increases, as indicated by \eqref{const1}
and \eqref{const2}.
\end{remark}

In order to prove the full or partial regularization of LSQR for
\eqref{eq1} completely and rigorously, besides Theorem~\ref{main1},
we need to introduce a precise definition of the near best
rank $k$ approximation $P_{k+1}B_kQ_k^T$ to $A$.
By definition \eqref{gammak}, the rank
$k$ matrix $P_{k+1}B_kQ_k^T$ is called a near best rank $k$ approximation
to $A$ if it satisfies
\begin{equation}\label{near}
\sigma_{k+1}\leq \gamma_k<\sigma_k \mbox{ and } \gamma_k-\sigma_{k+1}
<\sigma_k-\gamma_k,\mbox{ i.e., } \gamma_k<\frac{\sigma_k+\sigma_{k+1}}{2},
\end{equation}
that is, $\gamma_k$ lies between $\sigma_k$ and $\sigma_{k+1}$ and is closer to
$\sigma_{k+1}$. This definition is natural.
We mention in passing that a near best rank $k$ approximation to $A$ from
an ill-posed problem is much more stringent than it is for a matrix from
a numerically rank-deficient problem where the large singular values are
well separated from the small ones and there is a substantial gap
between two groups of singular values.

Based on Theorem~\ref{main1}, for the severely and moderately or mildly ill-posed
problems with the singular value models $\sigma_k=\zeta\rho^{-k}$ and
$\sigma_k=\zeta k^{-\alpha}$, we next derive the sufficient conditions on
$\rho$ and $\alpha$ that guarantee that $P_{k+1}B_kQ_k^T$ is a near best rank $k$
approximation to $A$ for $k=1,2,\ldots,k_0$. We analyze if and how the
sufficient conditions are satisfied for three kinds of ill-posed problems.

\begin{theorem}\label{nearapprox}
For a given \eqref{eq1}, assume that the discrete Picard condition
\eqref{picard} is satisfied. Then, in the sense of \eqref{near},
$P_{k+1}B_kQ_k^T$ is a near best rank $k$ approximation to $A$
for $k=1,2,\ldots,k_0$ if
\begin{equation}\label{condition}
\sqrt{1+\eta_k^2}<\frac{1}{2}\frac{\sigma_k}{\sigma_{k+1}}+\frac{1}{2}.
\end{equation}
For the severely ill-posed problems with $\sigma_k=\zeta\rho^{-k}$ and
the moderately or mildly ill-posed problems with $\sigma_k=\zeta k^{-\alpha}$,
$P_{k+1}B_kQ_k^T$ is a near best rank $k$ approximation to $A$
for $k=1,2,\ldots,k_0$ if $\rho>2$ and $\alpha$ satisfies
\begin{equation}\label{condition1}
2\sqrt{1+\eta_k^2}-1<\left(\frac{k_0+1}{k_0}\right)^{\alpha},
\end{equation}
respectively.
\end{theorem}

{\em Proof}.
By \eqref{final}, we see that $\gamma_k\leq \sqrt{1+\eta_k^2}\sigma_{k+1}$.
Therefore, $P_{k+1}B_kQ_k^T$ is a near best rank $k$ approximation to $A$ in
the sense of \eqref{near} provided that
$$
\sqrt{1+\eta_k^2}\sigma_{k+1}<\sigma_k
$$
and
$$
\sqrt{1+\eta_k^2}\sigma_{k+1}<\frac{\sigma_k+\sigma_{k+1}}{2},
$$
from which \eqref{condition} follows.

From \eqref{etak0}, for the severely ill-posed problems with
$\sigma_k=\zeta\rho^{-k}$ and $\rho>1$ we have
\begin{equation}\label{simp}
\sqrt{1+\eta_k^2}<1+\frac{1}{2}\eta_k^2\leq 1+\frac{1}{2}\rho^{-2(1+\beta)}
<1+\rho^{-1}, \ k=1,2,\ldots,k_0,
\end{equation}
from which it follows that
\begin{align}\label{ampli}
\sqrt{1+\eta_k^2}\sigma_{k+1}
&<(1+\rho^{-1})\sigma_{k+1}.
\end{align}
Since $\sigma_k/\sigma_{k+1}=\rho$, \eqref{condition} holds provided that
$$
1+\rho^{-1}<\frac{1}{2}\rho+\frac{1}{2},
$$
i.e., $\rho^2-\rho-2>0$, solving which for $\rho$ we get $\rho>2$. For the
moderately or mildly ill-posed problems with $\sigma_k=\zeta k^{-\alpha}$,
it is direct from \eqref{condition} to
get
$$
2\sqrt{1+\eta_k^2}-1<\left(\frac{k+1}{k}\right)^{\alpha}.
$$
Since $\left(\frac{k+1}{k}\right)^{\alpha}$ decreases
monotonically as $k$ increases, its minimum over $k=1,2,\ldots,k_0$
is $\left(\frac{k_0+1}{k_0}\right)^{\alpha}$.
Therefore, we obtain \eqref{condition1}.
\qed

\begin{remark}
Given the noise level $\|e\|$, the discrete Picard condition \eqref{picard}
and \eqref{picard1}, from the bound \eqref{const2} for
$\eta_k,\,k=1,2,\ldots,k_0$, we see that the bigger $\alpha>1$ is, the smaller
$k_0$ and $\eta_k$ are. Therefore,
there must be $\alpha>1$ such that \eqref{condition1} holds.
Here we should remind that it is more suitable to
regard the conditions on $\rho$ and $\alpha$ as an indication that
$\rho$ and $\alpha$ must not be close to one other than precise requirements
since we have used the bigger \eqref{simp} and simplified
models $\sigma_k=\zeta \rho^{-k}$ and $\sigma_k=\zeta k^{-\alpha}$.
\end{remark}

\begin{remark}
For the mildly ill-posed problems with $\sigma_k=\zeta k^{-\alpha}$,
Theorem~\ref{moderate} has shown that
$\|\Delta_k\|$ is generally not small and can be arbitrarily large
for $k=1,2,\ldots,k_0$. From \eqref{etadelta}, we see
that the size of $\eta_k$ is comparable to $\|\Delta_k\|$. Note that
the right-hand
side $\left(\frac{k_0+1}{k_0}\right)^{\alpha}\leq 2$ for
$\frac{1}{2}<\alpha\leq 1$ and any $k_0\geq 1$. Consequently, \eqref{condition1}
cannot be met generally for mildly ill-posed problems. The rare possible
exceptions are that $k_0$ is only very few and $\alpha$ is close to one since,
in such case, $\eta_k$ is not large for $k=1,2,\ldots,k_0$.
So, $P_{k+1}B_kQ_k^T$ is generally not a near best rank $k$ approximation
to $A$ for $k=1,2,\ldots, k_0$ for this kind of problem.
\end{remark}

\subsection{The approximation behavior of the Ritz values $\theta_i^{(k)}$}
\label{ritzapprox}

Starting with Theorem~\ref{main1}, we prove that, under
certain sufficient conditions on $\rho$ and $\alpha$ for the severely
and moderately ill-posed problems with the models $\sigma_i=\zeta\rho^{-i}$ and
$\sigma_i=\zeta i^{-\alpha}$, respectively,
the $k$ Ritz values $\theta_i^{(k)}$ approximate the first
$k$ large singular values $\sigma_i$ in natural order  for $k=1,2,\ldots,k_0$,
which means that no Ritz value smaller than $\sigma_{k_0+1}$ appears.
Combining this result with Theorem~\ref{nearapprox},
we can draw the definitive conclusion that LSQR must have the full
regularization for these two kinds of problems provided that $\rho>1$ and
$\alpha>1$ suitably. On the other hand, we will show why
LSQR generally has only the partial regularization for mildly ill-posed problems.

\begin{theorem}\label{ritzvalue}
Assume that \eqref{eq1} is severely ill-posed with
$\sigma_i=\zeta\rho^{-i}$ and $\rho>1$ or moderately ill-posed with
$\sigma_i=\zeta i^{-\alpha}$ and $\alpha>1$,
and the discrete Picard condition \eqref{picard} is
satisfied. Let the Ritz values $\theta_i^{(k)}$ be labeled
as $\theta_1^{(k)}>\theta_2^{(k)}>\cdots>\theta_{k}^{(k)}$.
Then
\begin{align}
0<\sigma_i-\theta_i^{(k)} &\leq \sqrt{1+\eta_k^2}\sigma_{k+1},\
i=1,2,\ldots,k.\label{error}
\end{align}
If $\rho\geq 1+\sqrt{2}$ or $\alpha>1$ satisfies
\begin{equation}\label{condm}
1+\sqrt{1+\eta_{k}^2}<\left(\frac{k_0+1}{k_0}\right)^{\alpha},\
k=1,2,\ldots,k_0,
\end{equation}
then the $k$ Ritz values $\theta_i^{(k)}$ strictly interlace
the first large $k+1$ singular values of $A$ and approximate
the first $k$ large ones in natural order for $k=1,2,\ldots,k_0$:
\begin{align}
\sigma_{i+1}&<\theta_i^{(k)}<\sigma_i,\,i=1,2,\ldots,k,
\label{error2}
\end{align}
meaning that there is no Ritz value $\theta_i^{(k)}$ smaller than $\sigma_{k_0+1}$
for $k=1,2,\ldots, k_0$.
\end{theorem}

{\em Proof}.
Note that for $k=1,2,\ldots,k_0$ the $\theta_i^{(k)},\ i=1,2,\ldots,k$ are
just the nonzero singular values of $P_{k+1}B_kQ_k^T$, whose other $n-k$
singular values are zeros. We write
$$
A=P_{k+1}B_k Q_k^T+(A-P_{k+1}B_k Q_k^T)
$$
with $\|A-P_{k+1}B_k Q_k^T\|=\gamma_k$
by definition \eqref{gammak}. Then by the Mirsky's theorem of
singular values \cite[p.204, Thm 4.11]{stewartsun}, we have
\begin{equation}\label{errbound}
| \sigma_i-\theta_i^{(k)}|\leq \gamma_k\leq
\sqrt{1+\eta_k^2}\sigma_{k+1},\ i=1,2,\ldots,k.
\end{equation}
Since the singular values of $A$ are simple and $b$ has components in all the
left singular vectors $u_1,u_2,\ldots, u_n$ of $A$, Lanczos bidiagonalization,
i.e., Algorithm~\ref{alg:lb}, can be run to completion, producing $P_{n+1},\ Q_n$ and
the lower bidiagonal $B_n\in \mathbb{R}^{(n+1)\times n}$ such that
\begin{equation}\label{fulllb}
P^TAQ_n=\left(\begin{array}{c}
B_n\\
\mathbf{0}
\end{array}
\right)
\end{equation}
with the $m\times m$ matrix $P=(P_{n+1},\hat{P})$ and $n\times n$ matrix $Q_n$
orthogonal and all
the $\alpha_i$ and $\beta_{i+1}$, $i=1,2,\ldots,n$, of $B_n$ being positive.
Note that the singular values of $B_k,\ k=1,2,\ldots,n,$
are all simple and that $B_k$ consists of the first $k$ columns of $B_n$
with the last $n-k$ {\em zero} rows deleted. Applying the Cauchy's {\em strict}
interlacing theorem \cite[p.198, Corollary 4.4]{stewartsun} to the singular
values of $B_k$ and $B_n$, we have
\begin{align}
\sigma_{n-k+i}< \theta_i^{(k)}&< \sigma_i,\ i=1,2,\ldots,k.
\label{interlace}
\end{align}
Therefore, \eqref{errbound} becomes
\begin{equation}\label{ritzapp}
0< \sigma_i-\theta_i^{(k)}\leq\gamma_{k}\leq
\sqrt{1+\eta_k^2}\sigma_{k+1},\ i=1,2,\ldots,k,
\end{equation}
which proves \eqref{error}.
That is, the $\theta_i^{(k)}$ approximate $\sigma_i$ from below
for $i=1,2,\ldots,k$ with the errors no more than
$\gamma_k\leq \sqrt{1+\eta_k^2}\sigma_{k+1}$.
For $i=1,2,\ldots,k$, notice that $\rho^{-k+i}\leq 1$. Then
from \eqref{ritzapp}, \eqref{simp} and $\sigma_i=\zeta\rho^{-i}$ we obtain
\begin{align*}
\theta_i^{(k)}&\geq \sigma_i-\gamma_k>\sigma_i-
(1+\rho^{-1})\sigma_{k+1}\\
&=\zeta\rho^{-i}-\zeta (1+\rho^{-1})\rho^{-(k+1)}\\
&=\zeta\rho^{-(i+1)}(\rho-(1+\rho^{-1})\rho^{-k+i})\\
&\geq \zeta\rho^{-(i+1)}(\rho-\rho^{-1}-1)\\
&\geq\zeta\rho^{-(i+1)}=\sigma_{i+1},
\end{align*}
provided that $\rho-\rho^{-1}\geq 2$, solving which we get $\rho\geq 1+\sqrt{2}$.
Together with the upper bound of \eqref{interlace}, we have proved \eqref{error2}.

For the moderately ill-posed problems with $\sigma_i=\zeta i^{-\alpha},\
i=1,2,\ldots,k$ and $k=1,2,\ldots,k_0$,
we get
\begin{align*}
\theta_i^{(k)}&\geq \sigma_i-\gamma_k\geq\sigma_i-\sqrt{1+\eta_k^2}
\sigma_{k+1}\\
&=\zeta i^{-\alpha}-\zeta \sqrt{1+\eta_k^2}(k+1)^{-\alpha}\\
&=\zeta (i+1)^{-\alpha}\left(\left(\frac{i+1}{i}\right)^{\alpha}
-\sqrt{1+\eta_k^2}\left(\frac{i+1}{k+1}\right)^{\alpha}\right)\\
&>\zeta (i+1)^{-\alpha}=\sigma_{i+1},
\end{align*}
i.e., \eqref{error2} holds, provided that $\eta_k>0$ and $\alpha>1$ are such that
$$
\left(\frac{i+1}{i}\right)^{\alpha}
-\sqrt{1+\eta_k^2}\left(\frac{i+1}{k+1}\right)^{\alpha}>1,
$$
which means that
$$
\sqrt{1+\eta_k^2}<\left(\left(\frac{i+1}{i}\right)^{\alpha}-1\right)
\left(\frac{k+1}{i+1}\right)^{\alpha}=
\left(\frac{k+1}{i}\right)^{\alpha}-\left(\frac{k+1}{i+1}\right)^{\alpha},\
i=1,2,\ldots,k.
$$
It is easily justified that the above right-hand side monotonically
decreases with respect to $i=1,2,\ldots,k$, whose minimum
attains at $i=k$ and equals $\left(\frac{k+1}{k}\right)^{\alpha}-1$.
Furthermore, since $\left(\frac{k+1}{k}\right)^{\alpha}-1$ decreases
monotonically as $k$ increases, its minimum over $k=1,2,\ldots,k_0$
is $\left(\frac{k_0+1}{k_0}\right)^{\alpha}-1$,
which is just the condition \eqref{condm}.
\qed

\begin{remark}
Similar to \eqref{condition1},
there must be $\alpha>1$ such that \eqref{condm} holds.
Comparing Theorem~\ref{nearapprox} with Theorem~\ref{ritzvalue}, we
find out that, as far as the severely or moderately ill-posed problems are
concerned, for $k=1,2,\ldots,k_0$ the near best rank approximation
$P_{k+1}B_kQ_k^T$ essentially means that the singular values
$\theta_i^{(k)}$ of $B_k$ approximate the first $k$
large singular values $\sigma_i$ of $A$ in natural order, provided that
$\rho>1$ or $\alpha>1$ suitably.
\end{remark}


\begin{remark}\label{extract}
In terms of the above remarks,
Theorems~\ref{main1}--\ref{ritzvalue} show that LSQR has the full
regularization for these two kinds of ill-posed
problems with $\rho>1$ and $\alpha>1$ suitably and can obtain
best possible regularized solutions $x^{(k_0)}$ at semi-convergence.
\end{remark}

For mildly ill-posed
problems. We observe that the sufficient condition \eqref{condm} for \eqref{error2}
is never met for this kind of problem because
$
\left(\frac{k_0+1}{k_0}\right)^{\alpha}\leq 2
$
for any $k_0$ and $\frac{1}{2}< \alpha\leq 1$. This indicates that,
for $k=1,2,\ldots,k_0$, the $k$ Ritz values $\theta_i^{(k)}$ may not approximate
the first $k$ large singular values $\sigma_i$ in natural order
and particularly there is at least one Ritz value
$\theta_{k_0}^{(k_0)}<\sigma_{k_0+1}$, causing that $x^{(k_0)}$ is already
deteriorated and cannot be as accurate as the best TSVD solution $x_{k_0}^{tsvd}$,
so that LSQR has only the partial regularization.
We can also make use of Theorem~\ref{initial} to explain the partial
regularization of LSQR: Theorem~\ref{moderate} has shown that
$\|\Delta_k\|$ is generally not small and
may become arbitrarily large as $k$ increases up to $k_0$ for mildly ill-posed
problems, meaning that $\|\sin\Theta(\mathcal{V}_k,\mathcal{V}_k^R)\|\approx 1$,
as the sharp bound \eqref{modera2} indicates, from which it follows that
a small Ritz value $\theta_{k_0}^{(k_0)}<\sigma_{k_0+1}$ generally appears.

\section{Decay rates of $\alpha_k$ and $\beta_{k+1}$ and their practical
importance}
\label{alphabeta}

In this section, we will present a number of results on the decay rates of
$\alpha_k$ and $\beta_{k+1}$.
The decay rates of $\alpha_k$ and $\beta_{k+1}$ are
particularly useful for practically detecting the degree of
ill-posedness of \eqref{eq1} and identifying the full or partial regularization
of LSQR. We
prove how $\alpha_k$ and $\beta_{k+1}$ decay by relating them to
$\gamma_k$ and the estimates established for it. Then we show how to
exploit the decay rate of $\alpha_k+\beta_{k+1}$ to identify
the degree of ill-posedness of \eqref{eq1} and the regularization of LSQR.

\begin{theorem}\label{main2}
With the notation defined previously, the following results hold:
\begin{eqnarray}
  \alpha_{k+1}&<&\gamma_k\leq \sqrt{1+\eta_k^2}\sigma_{k+1},
  \ k=1,2,\ldots,n-1,\label{alpha}\\
 \beta_{k+2}&<& \gamma_k\leq\sqrt{1+\eta_k^2}\sigma_{k+1}, \ k=1,2,\ldots,n-1,
 \label{beta}\\
\alpha_{k+1}\beta_{k+2}&\leq &
\frac{\gamma_k^2}{2}\leq
   \frac{(1+\eta_k^2)\sigma_{k+1}^2}{2}, \ k=1,2,\ldots,n-1,
  \label{prod2}\\
\gamma_{k+1}&<&\gamma_k,\  \ k=1,2,\ldots,n-2. \label{gammamono}
\end{eqnarray}
\end{theorem}

{\em Proof}.
From \eqref{fulllb}, since $P$ and $Q_n$ are orthogonal matrices,
we have
\begin{align}
\gamma_k &=\|A-P_{k+1}B_kQ_k^T\|=\|P^T(A-P_{k+1}B_kQ_k^T)Q_n\| \label{invar}\\
&=
\left\| \left(\begin{array}{c}
B_n \\
\mathbf{0}
\end{array}
\right)-(I,\mathbf{0} )^TB_k (I,\mathbf{0} )\right\|=\|G_k\| \label{gk}
\end{align}
with
\begin{align}\label{gk1}
G_k&=\left(\begin{array}{cccc}
\alpha_{k+1} & & & \\
\beta_{k+2}& \alpha_{k+2} & &\\
&
  \beta_{k+3} &\ddots & \\& & \ddots & \alpha_{n} \\
  & & & \beta_{n+1}
  \end{array}\right)\in \mathbb{R}^{(n-k+1)\times (n-k)}
\end{align}
resulting from deleting the $(k+1)\times k$ leading principal matrix of $B_n$
and the first $k$ zero rows and columns of the resulting matrix.
From the above, for $k=1,2,\ldots,n-1$ we have
\begin{align}
\alpha_{k+1}^2+\beta_{k+2}^2&=\|G_ke_1\|^2\leq \|G_k\|^2=\gamma_k^2,
\label{alphabetasum1}
\end{align}
which shows that $\alpha_{k+1}< \gamma_k$ and $\beta_{k+2}<\gamma_k$
since $\alpha_{k+1}>0$ and $\beta_{k+2}>0$.
So from \eqref{final}, we get \eqref{alpha} and \eqref{beta}. On the other
hand, noting that
\begin{align*}
2\alpha_{k+1}\beta_{k+2}&\leq \alpha_{k+1}^2+\beta_{k+2}^2\leq \gamma_k^2,
\end{align*}
we get \eqref{prod2}.

Note that $\alpha_k>0$ and $\beta_{k+1}>0,\ k=1,2,\ldots,n$.
By $\gamma_k=\|G_k\|$ and \eqref{gk1}, note that $\gamma_{k+1}=\|G_{k+1}\|$
equals the 2-norm of the submatrix deleting the first column of $G_k$.
Applying the Cauchy's strict interlacing theorem to the singular values
of this submatrix and $G_k$, we obtain \eqref{gammamono}.
\qed

\begin{remark}
For severely and moderately ill-posed problems, based on
the results in the last section,
\eqref{alpha} and \eqref{beta} show that $\alpha_{k+1}$ and $\beta_{k+2}$
decay as fast as $\sigma_{k+1}$ for $k\leq k_0$ and their decays may become slow
for $k>k_0$. For mildly ill-posed problems,
since $\eta_k$ are generally bigger than one considerably for $k\leq k_0$,
$\alpha_{k+1}$ and $\beta_{k+2}$ cannot generally decay as fast as
$\sigma_{k+1}$, and their decays become slower for $k>k_0$.
\end{remark}

We now shed light on \eqref{alpha} and \eqref{beta}.
For a given \eqref{eq1}, its degree of ill-posedness
is either known or unknown. If it is unknown, \eqref{alpha} is of
practical importance and can be exploited to identify whether or not LSQR has
the full regularization without extra cost in an automatic and
reliable way, so is \eqref{beta}.
From the proofs of \eqref{alpha} and \eqref{beta}, we find that
$\alpha_{k+1}$ and $\beta_{k+2}$ are as small as $\gamma_k$. Since our theory
and analysis in Section \ref{rankapp} have proved that $\gamma_k$
decays as fast as $\sigma_{k+1}$ for severely or moderately ill-posed problems
with $\rho>1$ or $\alpha>1$ suitably and it decays more slowly than
$\sigma_{k+1}$ for mildly il-posed problems, the decay rate
of $\sigma_k$ can be judged by that of $\alpha_k$
or $\beta_{k+1}$ or better judged by that of $\alpha_k+\beta_{k+1}$
reliably, as shown below.

Given \eqref{eq1}, run LSQR until semi-convergence
occurs at iteration $k^*$. Check how $\alpha_k+\beta_{k+1}$ decays as $k$
increases during the process. If, on average, it decays in
an obviously exponential way, then \eqref{eq1} is a severely ill-posed problem.
In this case, LSQR has the full regularization, and semi-convergence means
that we have found a best possible regularized solution. If, on average,
$\alpha_k$ decays as fast as $k^{-\alpha}$ with $\alpha>1$ considerably, then
\eqref{eq1} is surely a moderately ill-posed problem, and LSQR also has found a
best possible regularized solution at semi-convergence. If, on average,
it decays at most as fast as or more slowly than $k^{-\alpha}$
with $\alpha$ no more than one, \eqref{eq1} is a mildly ill-posed problem.
Notice that the noise $e$ does not deteriorate regularized solutions until
semi-convergence. Therefore, if a hybrid LSQR is used, then it is more reasonable
and also cheaper to apply regularization to projected problems only
from iteration $k^*+1$ onwards
other than from the first iteration, as done in the hybrid
Lanczos bidiagonalization/Tikhonov regularization scheme \cite{berisha}, until
a best possible regularized solution is found.

\section{Numerical experiments}\label{numer}

Huang and Jia \cite{huangjia} have numerically justified the full
regularization of LSQR for severely and moderately ill-posed problems and its
partial regularization for mildly ill-posed problems \cite{hansen07}, where
each $A$ is $1,024\times 1,024$.
In this section, we report numerical experiments to
confirm our theory and illustrate the full or partial regularization
of LSQR in much more detail. For the first two kinds of problems,
we demonstrate that $\gamma_k,\ \alpha_{k+1}$
and $\beta_{k+2}$ decay as fast as $\sigma_{k+1}$. We compare LSQR
and the hybrid LSQR with the TSVD method applied to projected
problems after semi-convergence.
For each of severely and moderately ill-posed problems, we show that the
regularized solution obtained by LSQR at semi-convergence is at least as
accurate as the best TSVD regularized solution, indicating
that LSQR has the full regularization. In the meantime, for mildly ill-posed
problems, we show that the regularized solution obtained by LSQR at
semi-convergence is considerably less accurate than $x_{k_0}^{tsvd}$,
demonstrating that LSQR has only the partial regularization.

We choose several ill-posed problems from Hansen's regularization
toolbox \cite{hansen07}, which include
the severely ill-posed problems $\mathsf{shaw,\ wing}$,
the moderately ill-posed problems $\mathsf{heat,\ phillips}$, and
the mildly ill-posed problem $\mathsf{deriv2}$ with the parameter
"example=3".
All the codes are from \cite{hansen07}, and the problems
arise from discretizations of \eqref{eq2}.
We remind that, as far as solving \eqref{eq1} is concerned,
our primary goal consists in justifying the regularizing effects
of iterative solvers for \eqref{eq1}, which are {\em unaffected by the size}
of \eqref{eq1} and only depends on the degree of
ill-posedness, the noise level $\|e\|$ and the actual discrete Picard
condition, provided that the condition number of \eqref{eq1},
measured by the ratio between the largest and smallest singular values
of each $A$, is large enough.
Therefore, for this purpose, as extensively done in the
literature (see, e.g., \cite{hansen98,hansen10} and the references therein
as well as many other papers), it is enough to report the results on small
and/or medium sized discrete ill-posed problems since the condition numbers of
these $A$ are already huge or large, which, in finite precision arithmetic,
are roughly $10^{16}, 10^{8}$ and $10^{6}$ for severely, moderately and
mildly ill-posed problems with $n=256$, respectively.
Indeed, for $n$ large, say, 10,000 or more, we have observed that
LSQR has the same behavior as for small $n$,
e.g., $n=256$, which is used in this paper. The only exception is
$\mathsf{deriv2}$, and we will test a larger one of $n=3,000$ whose
condition number is one order larger than that of $n=256$, so as to better
confirm the partial regularization of LSQR. Also,
an important reason is that such choice enables us to fully
justify the regularization effects of LSQR by comparing it with
the TSVD method, which suits only for small and/or medium sized problems
because of its computational complexity. For each example,
we generate $A$, $x_{true}$ and $\hat{b}$.
In order to simulate the noisy data, we generate white noise
vectors $e$ such that the relative noise levels
$\varepsilon=\frac{\|e\|}{\|\hat{b}\|}=10^{-2}, 10^{-3}, 10^{-4}$, respectively.
To simulate exact arithmetic,
LSQR uses full reorthogonalization in Lanczos bidiagonalization.
All the computations are carried out in Matlab 7.8 with the machine precision
$\epsilon_{\rm mach}= 2.22\times10^{-16}$ under the Miscrosoft
Windows 7 64-bit system.

\subsection{The accuracy of rank $k$ approximations}

\begin{figure}
\begin{minipage}{0.48\linewidth}
  \centerline{\includegraphics[width=6.0cm,height=4.5cm]{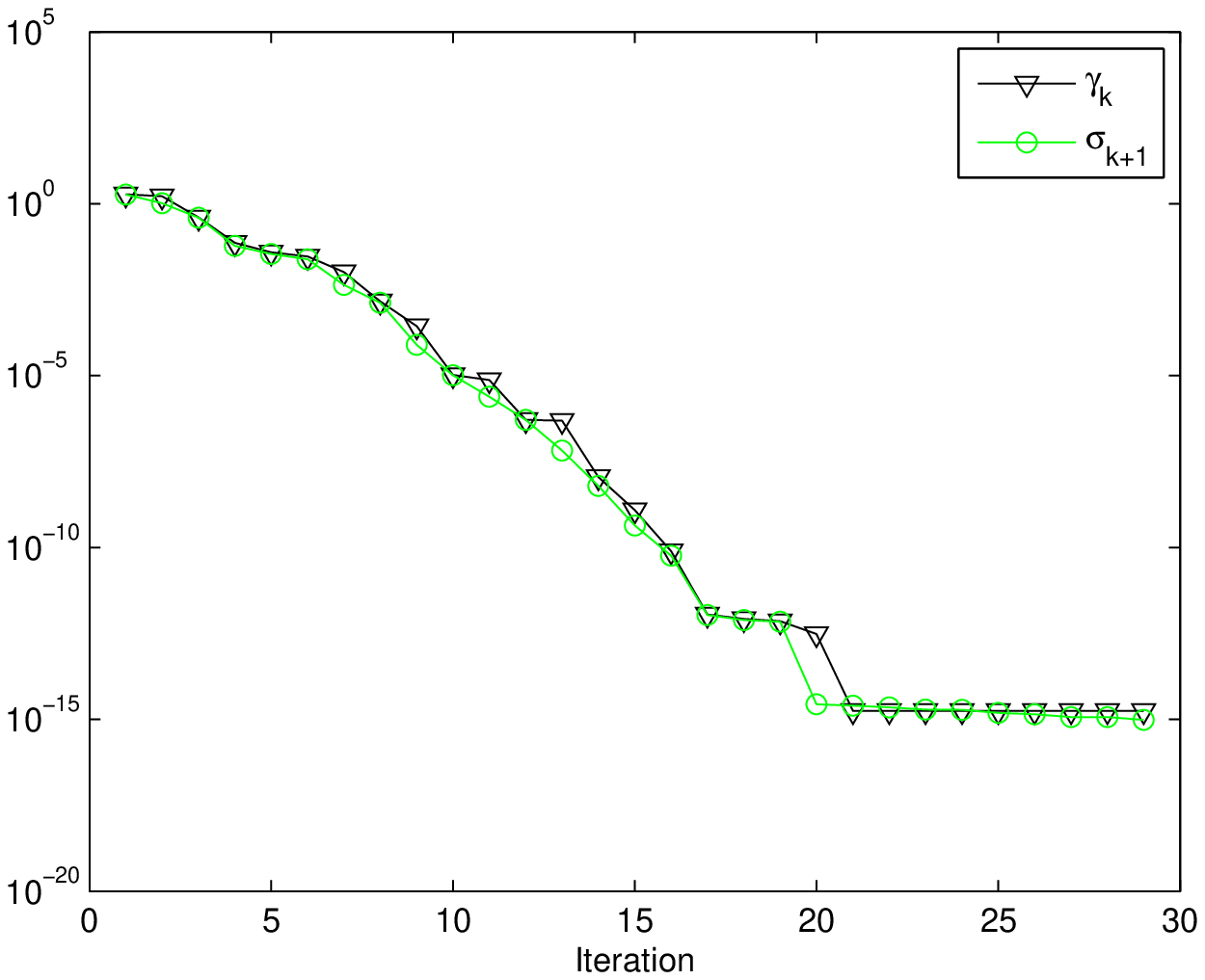}}
  \centerline{(a)}
\end{minipage}
\hfill
\begin{minipage}{0.48\linewidth}
  \centerline{\includegraphics[width=6.0cm,height=4.5cm]{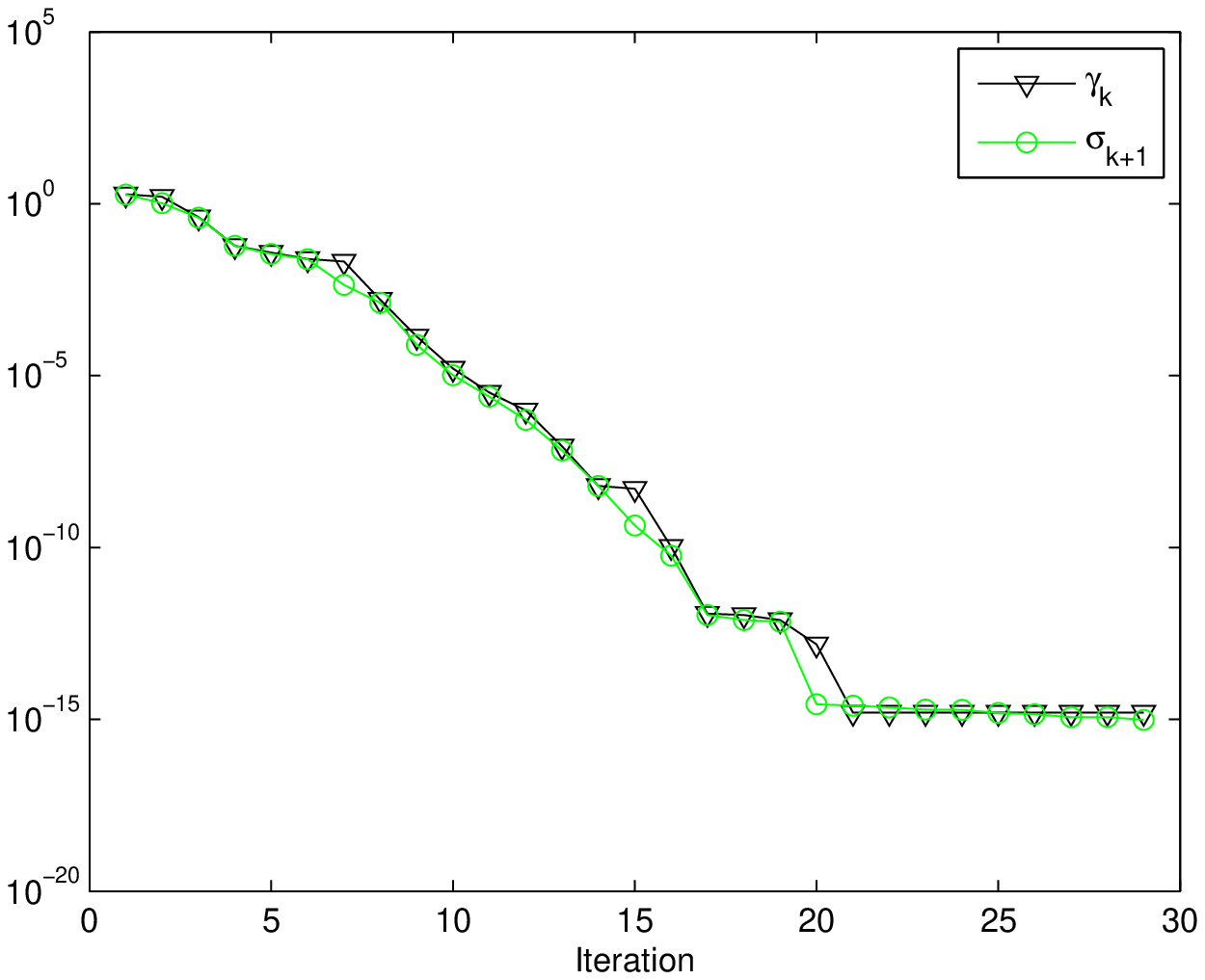}}
  \centerline{(b)}
\end{minipage}
\vfill
\begin{minipage}{0.48\linewidth}
  \centerline{\includegraphics[width=6.0cm,height=4.5cm]{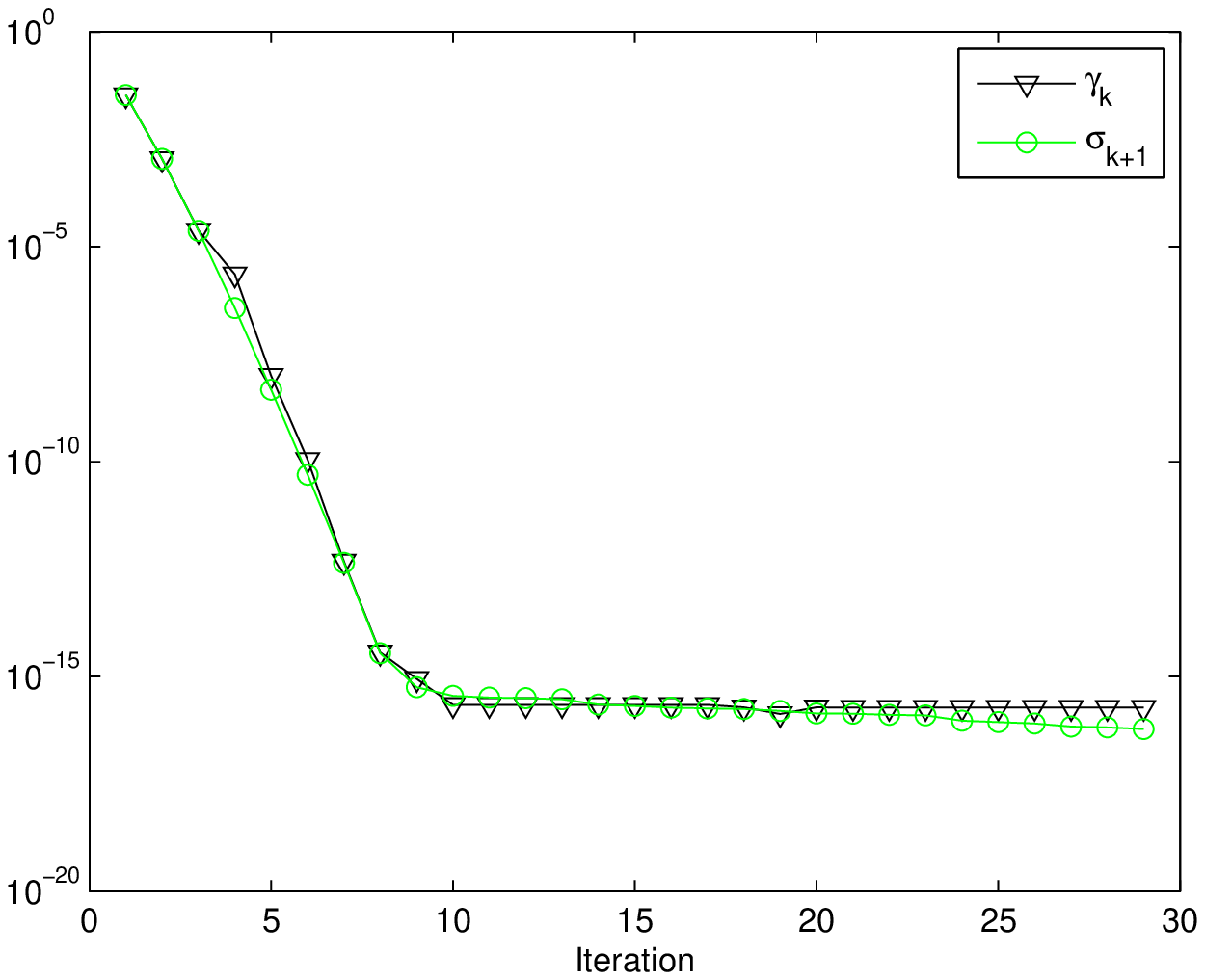}}
  \centerline{(c)}
\end{minipage}
\hfill
\begin{minipage}{0.48\linewidth}
  \centerline{\includegraphics[width=6.0cm,height=4.5cm]{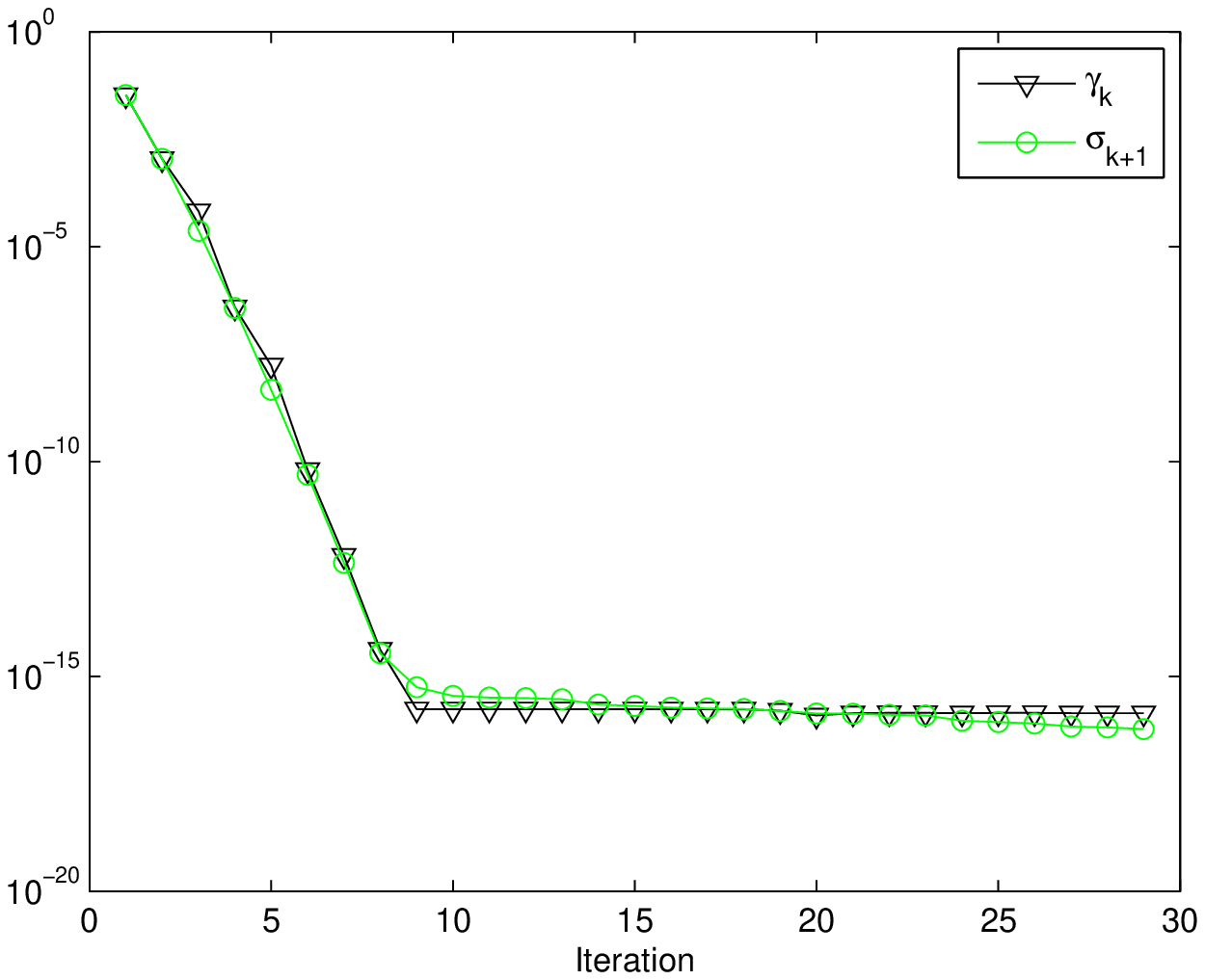}}
  \centerline{(d)}
\end{minipage}
\caption{(a)-(b): Decay curves of the sequences $\gamma_k$ and
$\sigma_{k+1}$ for $\mathsf{shaw}$ with $\varepsilon=10^{-2}$
(left) and $\varepsilon=10^{-3}$ (right); (c)-(d): Decay curves of the
sequences $\gamma_k$ and $\sigma_{k+1}$ for $\mathsf{wing}$ with
$\varepsilon=10^{-3}$ (left) and $\varepsilon=10^{-4}$ (right).} \label{fig1}
\end{figure}

In Figure~\ref{fig1}, we display the decay curves of the $\gamma_k$ for
$\mathsf{shaw}$ with $\varepsilon=10^{-2}, 10^{-3}$ and for $\mathsf{wing}$
with $\varepsilon=10^{-3}, 10^{-4}$, respectively.
We observe that the three curves with different $\varepsilon$ are almost
unchanged. This is in accordance with our Remark~\ref{decayrate}, where
it is stated that the decay rate of $\gamma_k$ is little affected
by noise levels for severely ill-posed problems, since
$\gamma_k$ primarily depends on the decay rate of $\sigma_{k+1}$
and different noise levels only affect the value of
$k_0$ other than the decay rate of $\gamma_k$. In addition,
we have observed that $\gamma_k$ and $\sigma_{k+1}$ decay
until they level off at $\epsilon_{\rm mach}$ due to round-off errors.
Most importantly, the results have clearly confirmed the theory
that $\gamma_k$ decreases as fast as $\sigma_{k+1}$, and we have
$\gamma_k\approx\sigma_{k+1}$, whose decay curves
are almost indistinguishable.

\begin{figure}
\begin{minipage}{0.48\linewidth}
  \centerline{\includegraphics[width=6.0cm,height=4.5cm]{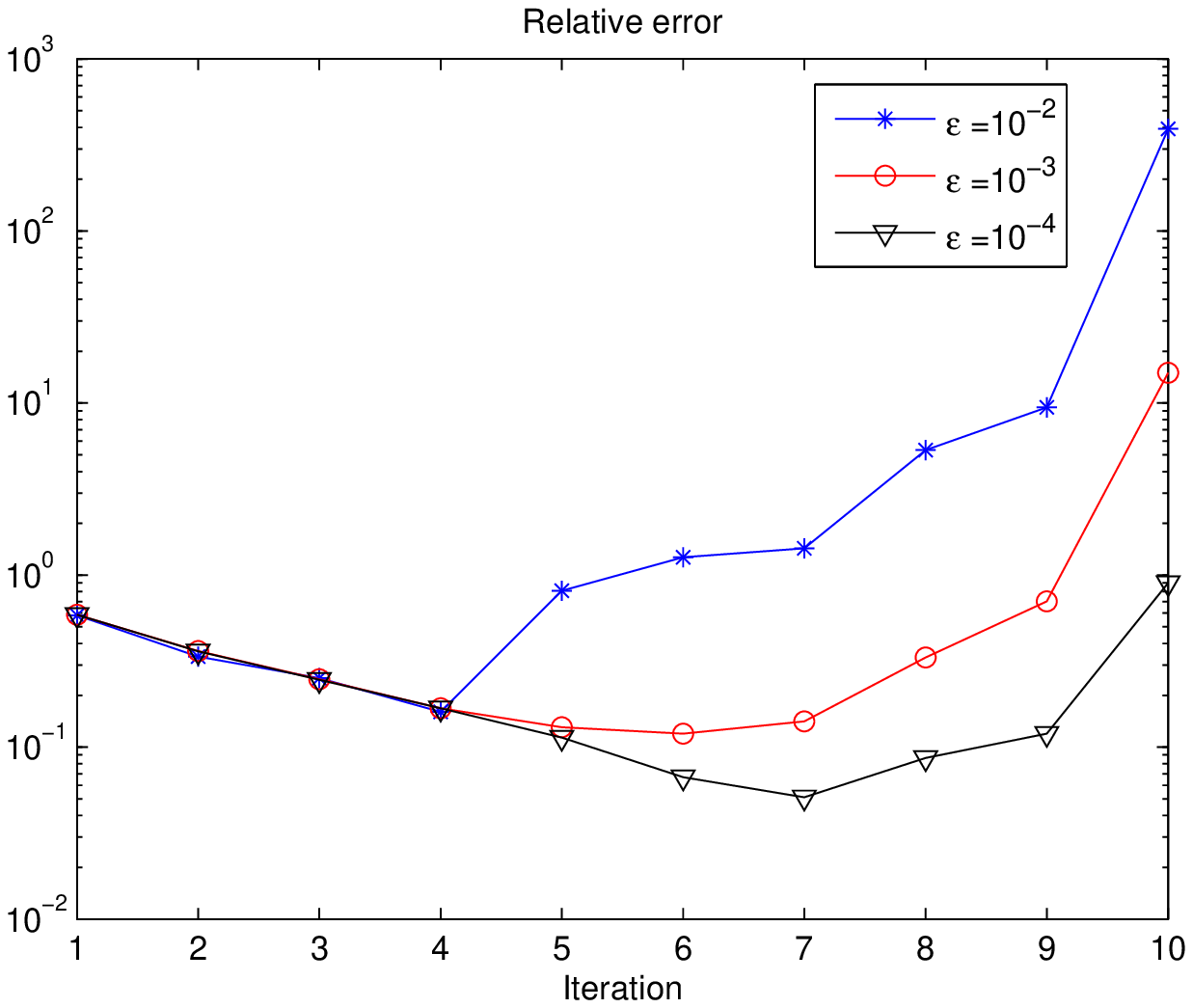}}
  \centerline{(a)}
\end{minipage}
\hfill
\begin{minipage}{0.48\linewidth}
  \centerline{\includegraphics[width=6.0cm,height=4.5cm]{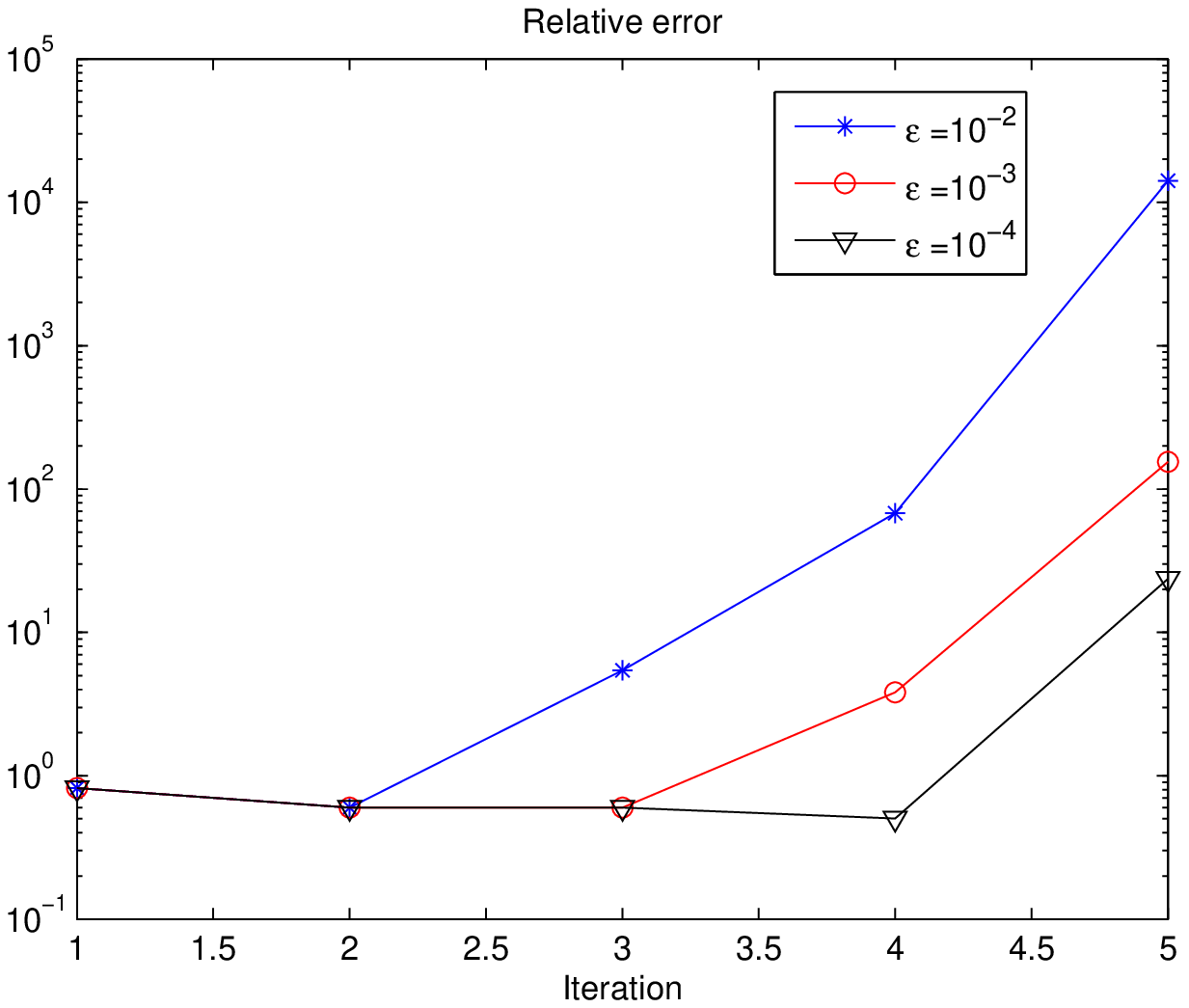}}
  \centerline{(b)}
\end{minipage}
\caption{ The relative errors $\|x^{(k)}-x_{true}\|/
\|x_{true}\|$ with $\varepsilon=10^{-2}, 10^{-3}, 10^{-4}$
for $\mathsf{shaw}$ (left) and $\mathsf{wing}$ (right).}
\label{fig2}
\end{figure}

In Figure~\ref{fig2}, we plot the relative errors $\|x^{(k)}
-x_{true}\|/\|x_{true}\|$ with different $\varepsilon$ for
these two problems. As we have seen, LSQR exhibits clear semi-convergence.
Moreover, for a smaller $\varepsilon$, we get a more accurate regularized
solution at cost of more iterations, as $k_0$ is bigger from \eqref{picard}
and \eqref{picard1}.

\begin{figure}
\begin{minipage}{0.48\linewidth}
  \centerline{\includegraphics[width=6.0cm,height=4.5cm]{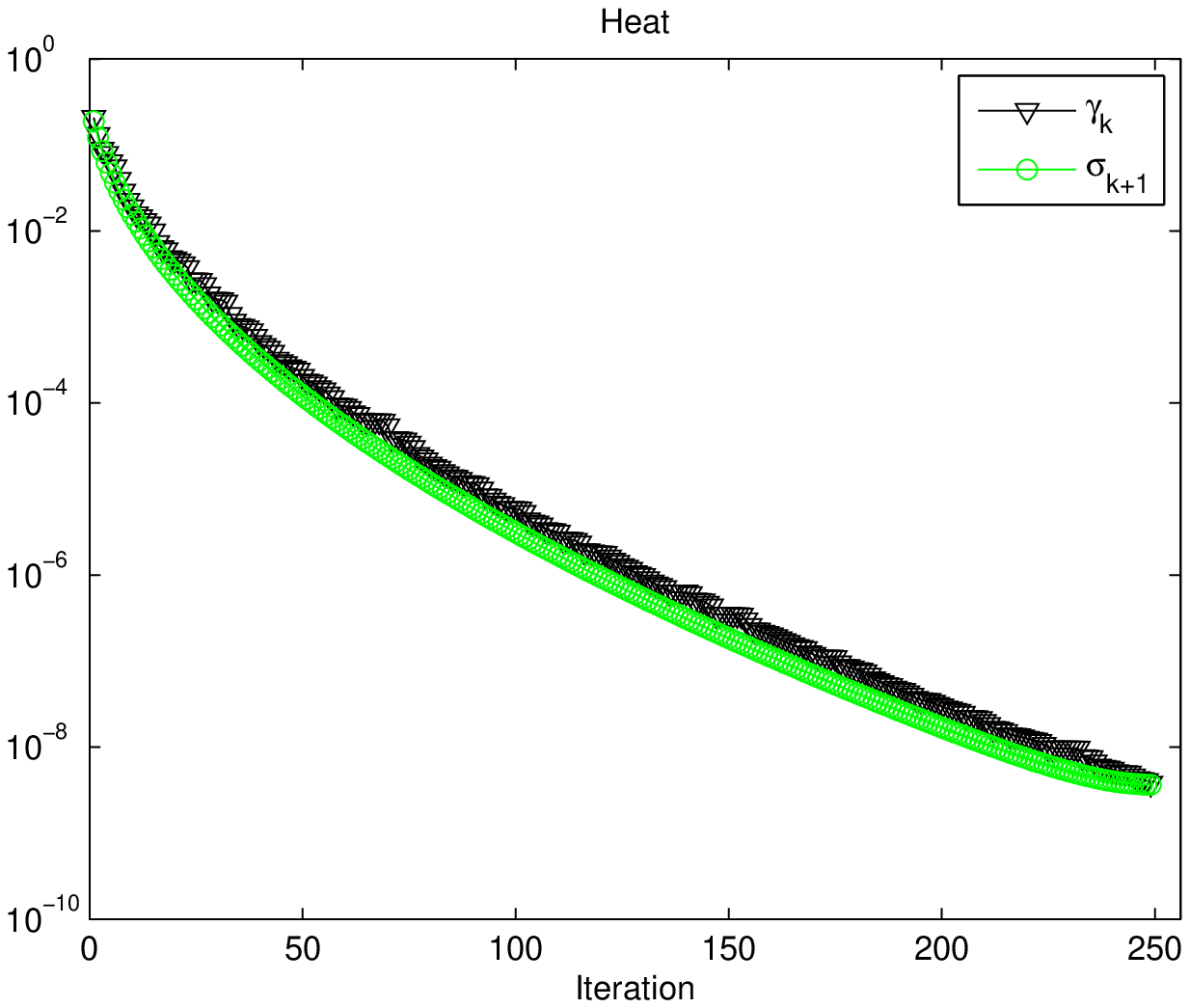}}
  \centerline{(a)}
\end{minipage}
\hfill
\begin{minipage}{0.48\linewidth}
  \centerline{\includegraphics[width=6.0cm,height=4.5cm]{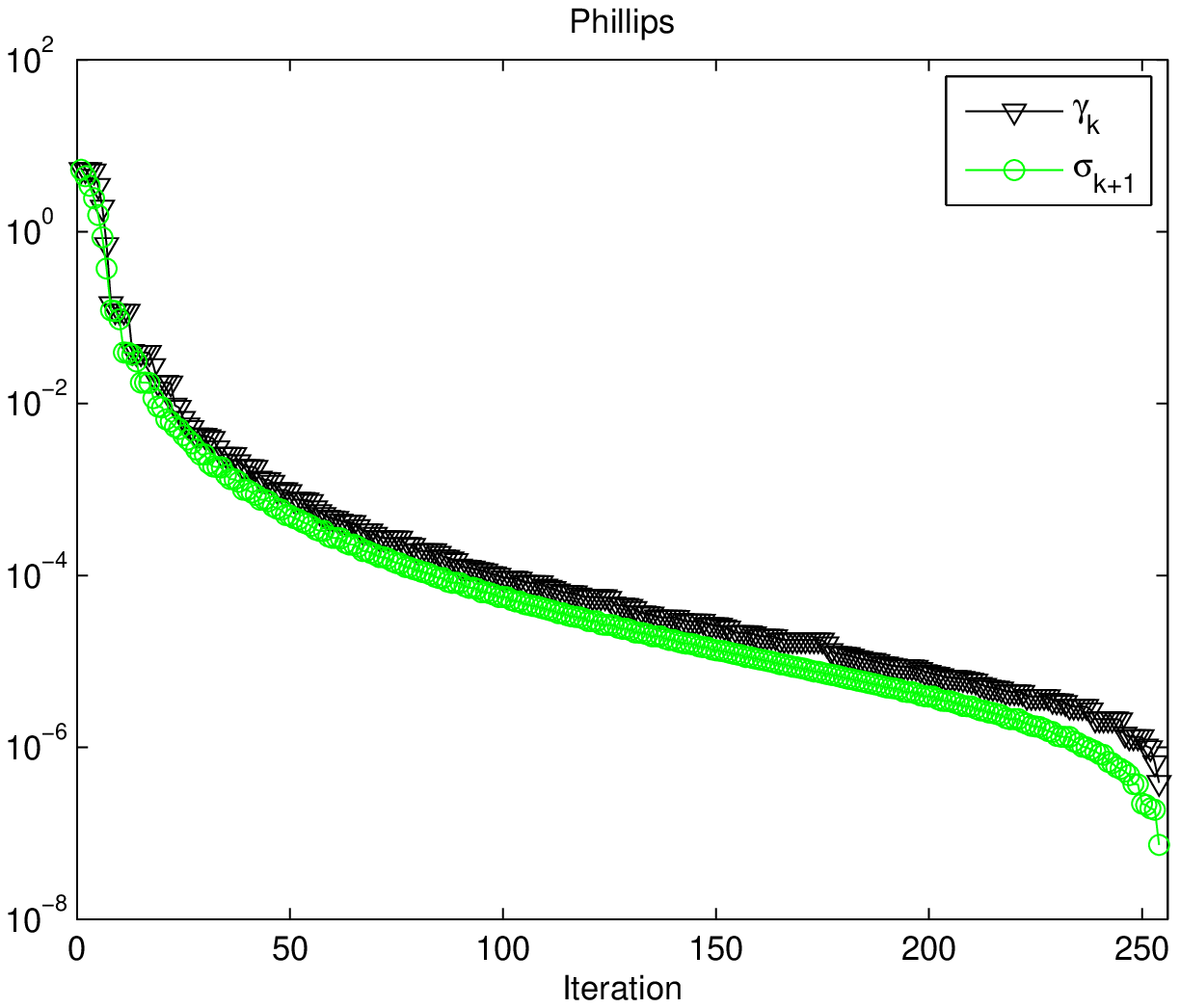}}
  \centerline{(b)}
\end{minipage}
\caption{(a): Decay curves of the sequences $\gamma_k$ and $\sigma_{k+1}$
for $\mathsf{heat}$ with (left) and (b): Decay curves of the
sequences $\gamma_k$ and $\sigma_{k+1}$ for
$\mathsf{phillips}$ with $\varepsilon=10^{-3}$ (right).}\label{fig3}
\end{figure}

From Figure~\ref{fig3}, we see that $\gamma_k$ decreases almost as fast as
$\sigma_{k+1}$ for the moderately ill-posed problems $\mathsf{heat}$ and
$\mathsf{phillips}$. However, slightly different
from severely ill-posed problems, $\gamma_k$, though
excellent approximations to $\sigma_{k+1}$, may not be so very accurate.
This is expected, as the constants $\eta_k$ in \eqref{const2} are generally
bigger than those in \eqref{const1} for severely ill-posed problems. Also,
different from Figure~\ref{fig1}, we observe from Figure~\ref{fig3} that
$\gamma_k$ deviates more from $\sigma_{k+1}$ with $k$ increasing, especially
for the problem $\mathsf{phillips}$. This confirms
Remarks~\ref{decayrate}--\ref{decayrate2} on moderately ill-posed problems.

\begin{figure}
\begin{minipage}{0.48\linewidth}
  \centerline{\includegraphics[width=6.0cm,height=4.5cm]{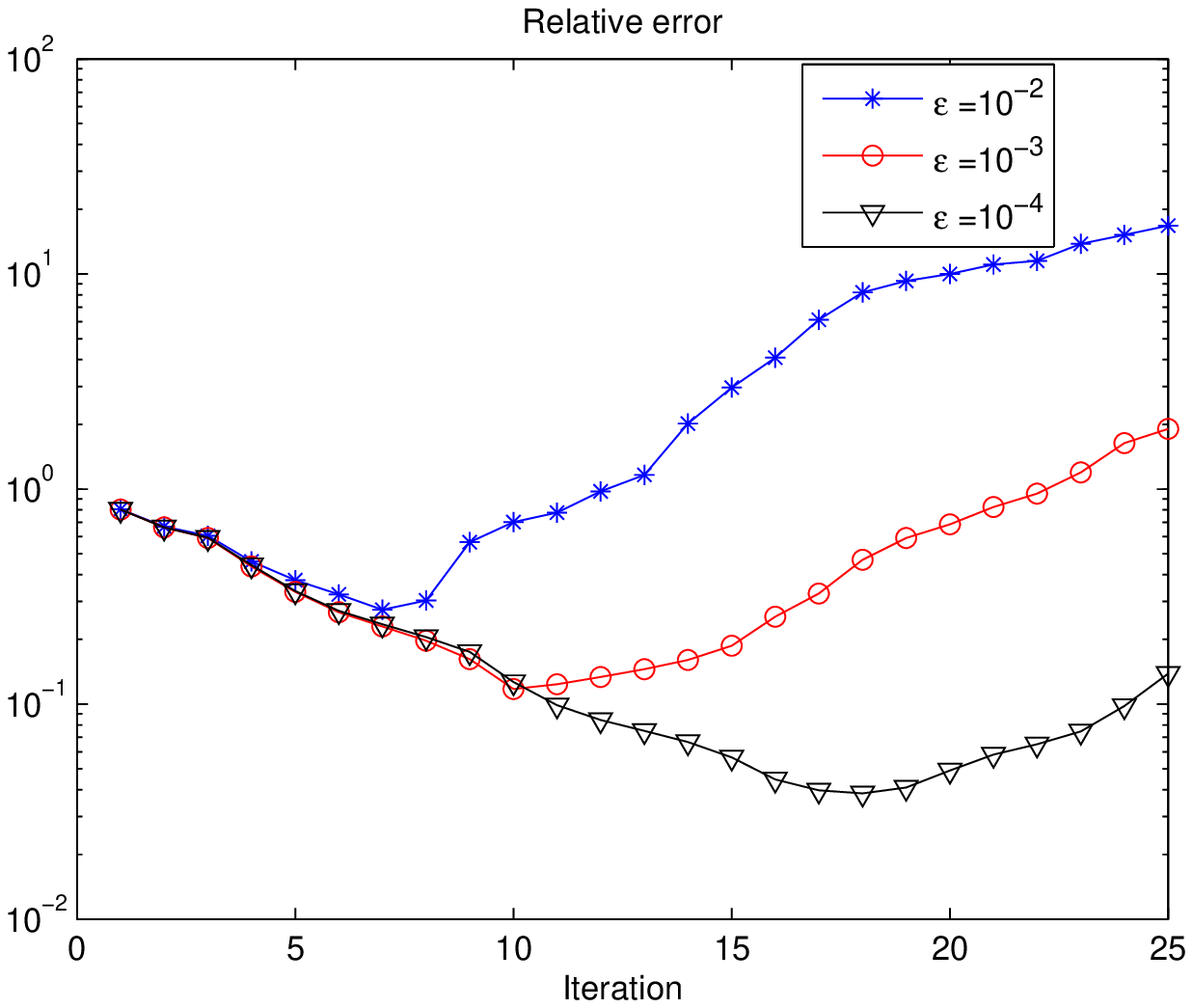}}
  \centerline{(a)}
\end{minipage}
\hfill
\begin{minipage}{0.48\linewidth}
  \centerline{\includegraphics[width=6.0cm,height=4.5cm]{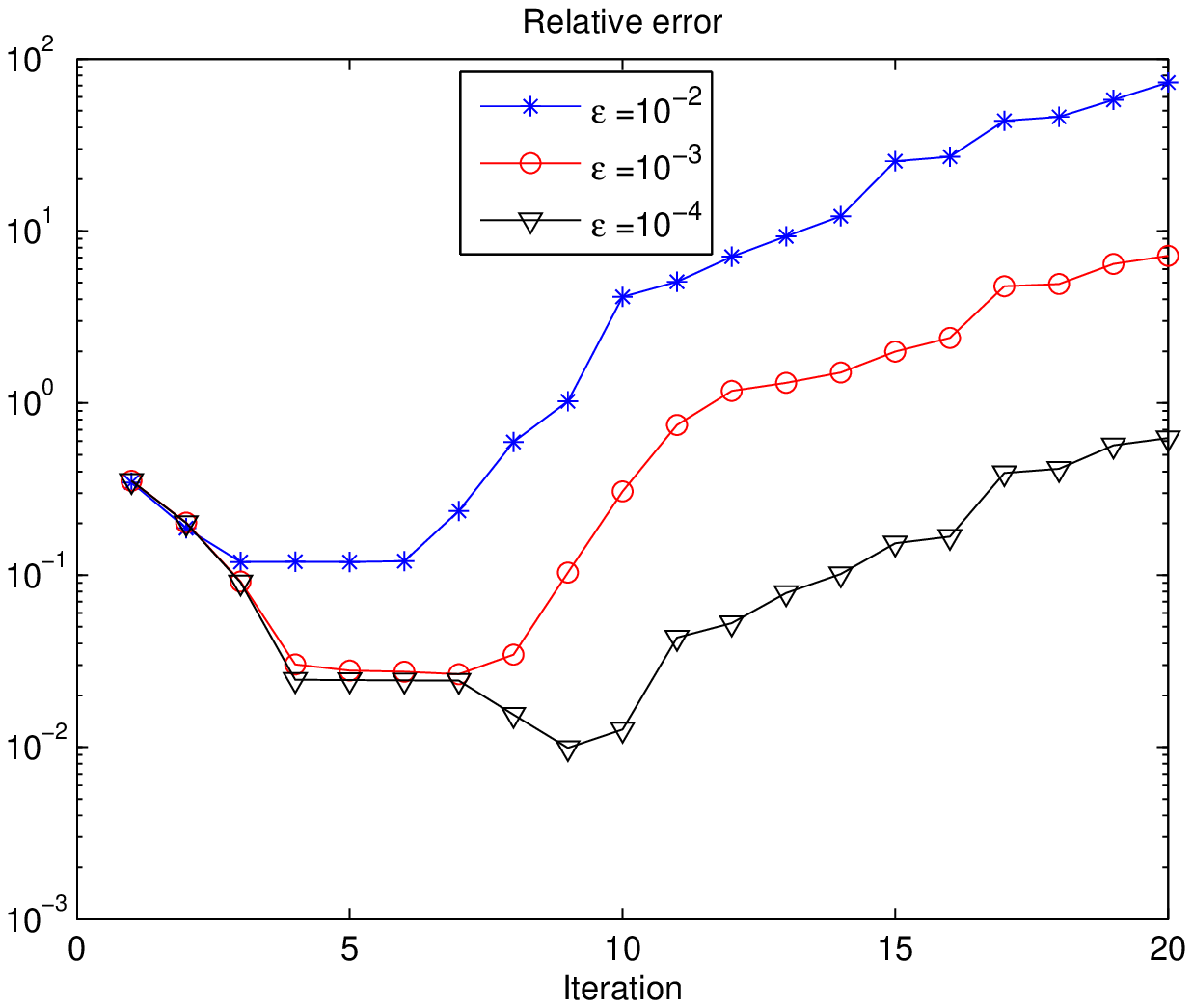}}
  \centerline{(b)}
\end{minipage}
\caption{ The relative errors $\|x^{(k)}-x_{true}\|/\|x_{true}\|$
with $\varepsilon=10^{-2}, 10^{-3}, 10^{-4}$ for
$\mathsf{heat}$ (left) and $\mathsf{phillips}$ (right).} \label{fig4}
\end{figure}

In Figure~\ref{fig4}, we depict the relative errors of $x^{(k)}$, and from them
we observe analogous phenomena to those for severely ill-posed problems.
The only distinction is that LSQR now needs more iterations, i.e.,
a bigger $k_0$ is needed for moderately ill-posed problems with the same
$\varepsilon$, as is seen from \eqref{picard} and \eqref{picard1}.

\begin{figure}
\begin{minipage}{0.48\linewidth}
  \centerline{\includegraphics[width=6.0cm,height=4.5cm]{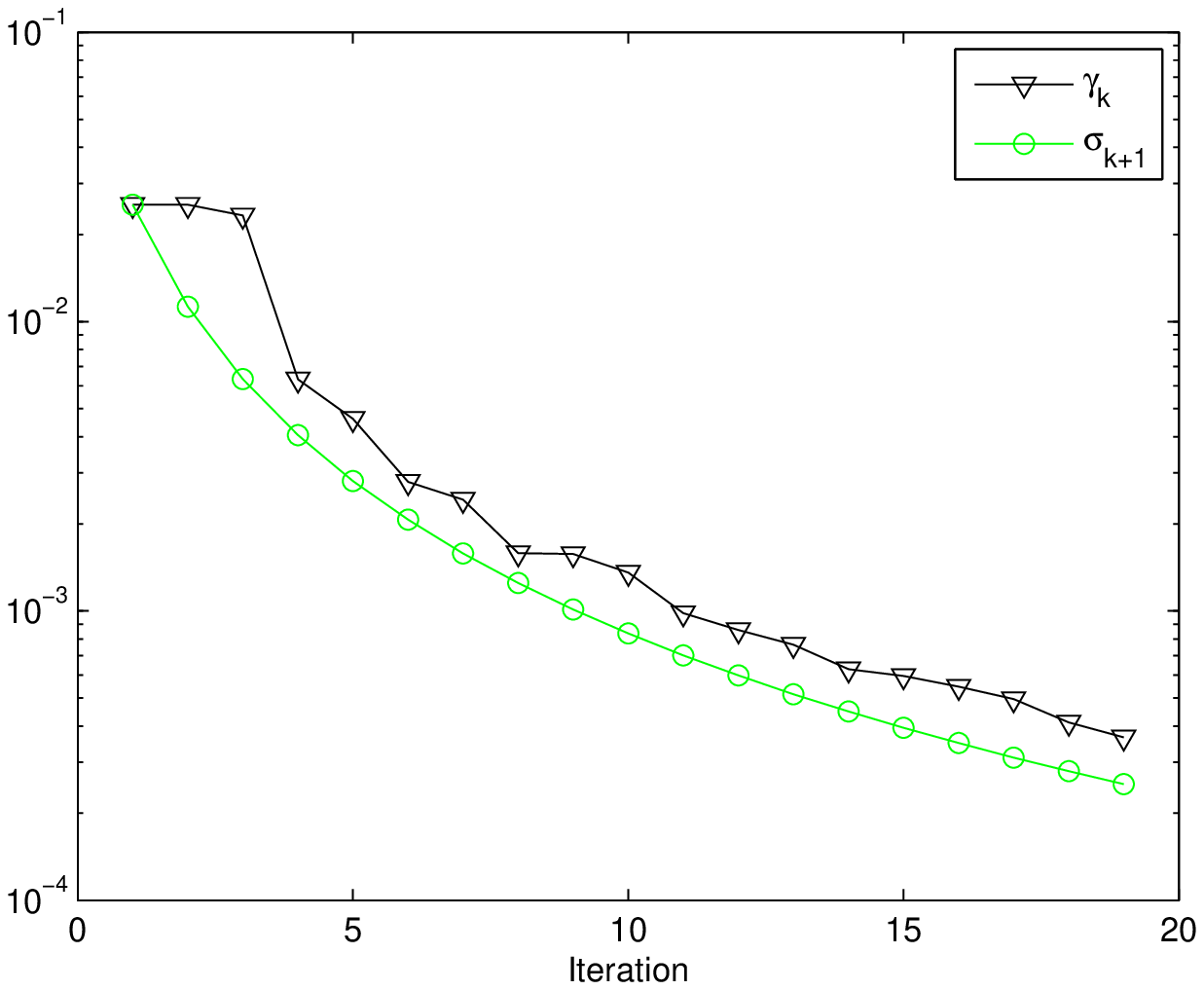}}
  \centerline{(a)}
\end{minipage}
\hfill
\begin{minipage}{0.48\linewidth}
  \centerline{\includegraphics[width=6.0cm,height=4.5cm]{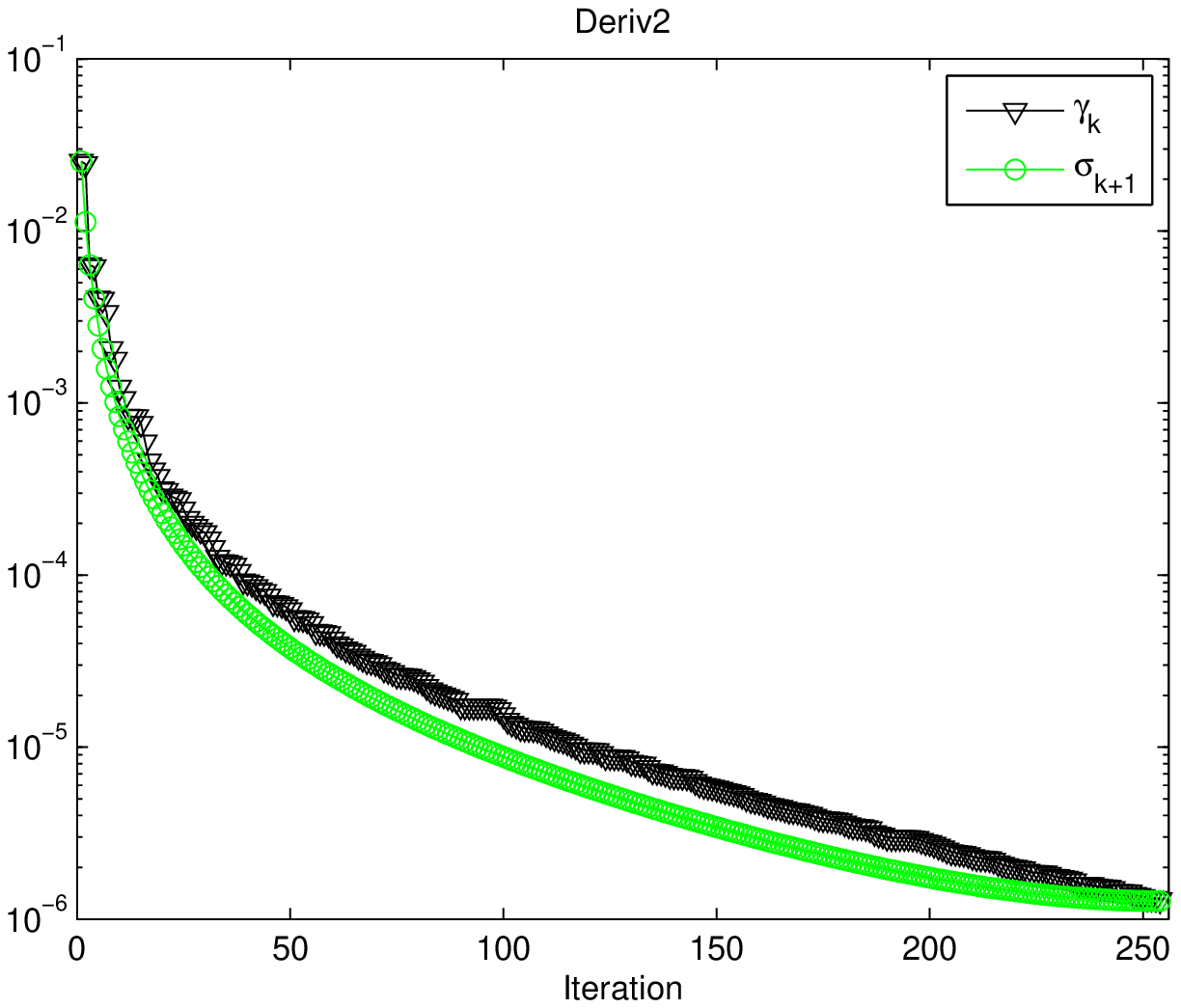}}
  \centerline{(b)}
\end{minipage}
\caption{(a)-(b): Decay curves of the partial and complete sequences
$\gamma_k$ and $\sigma_{k+1}$ for $\mathsf{deriv2}$ with
$\varepsilon=10^{-3}$}
\label{figmild}
\end{figure}

Figure~\ref{figmild} (a)-(b) display
the decay curves of the partial and complete sequences $\gamma_k$ and
$\sigma_{k+1}$ for the mildly ill-posed problem {\sf deriv2},
respectively. We see that, different from severely and
moderately ill-posed problems, $\gamma_k$ does not decay so fast as
$\sigma_{k+1}$ and deviates from $\sigma_{k+1}$ significantly.
These observations justify our theory and confirm
that the rank $k$ approximations to $A$ generated by Lanczos
bidiagonalization are not as accurate as those for severely and moderately
problems.

\subsection{Decay behavior of $\alpha_k$ and $\beta_{k+1}$}

For the severely ill-posed $\mathsf{shaw, wing}$ and the moderately
ill-posed $\mathsf{heat, phillips}$, we now
illustrate that $\alpha_k$ and $\beta_{k+1}$ decay as fast as the singular
values $\sigma_k$ of $A$. We take the noise level $\varepsilon=10^{-3}$. The results
are similar for $\varepsilon=10^{-2}$ and $10^{-4}$.

\begin{figure}
\begin{minipage}{0.48\linewidth}
  \centerline{\includegraphics[width=6.0cm,height=4.5cm]{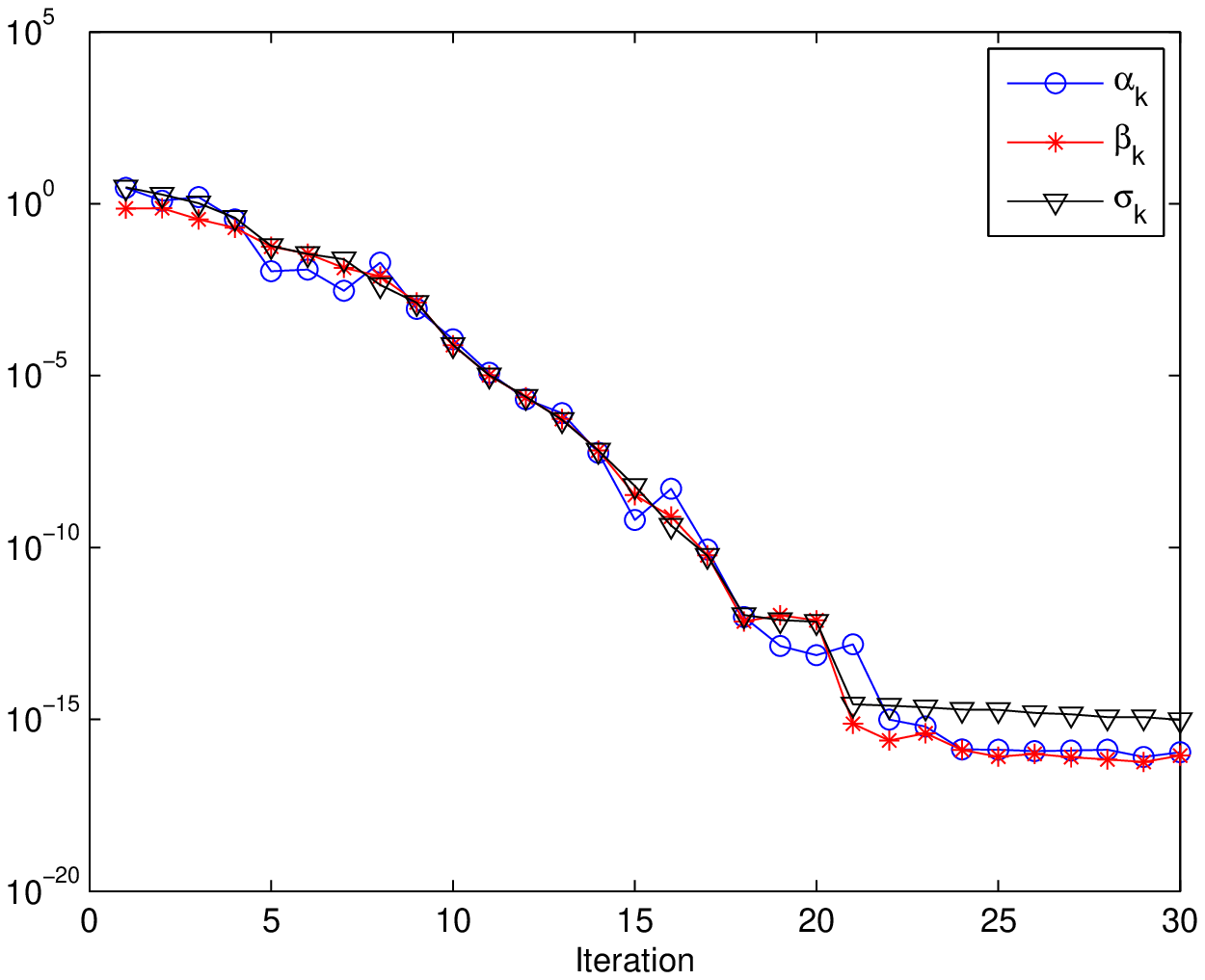}}
  \centerline{(a)}
\end{minipage}
\hfill
\begin{minipage}{0.48\linewidth}
  \centerline{\includegraphics[width=6.0cm,height=4.5cm]{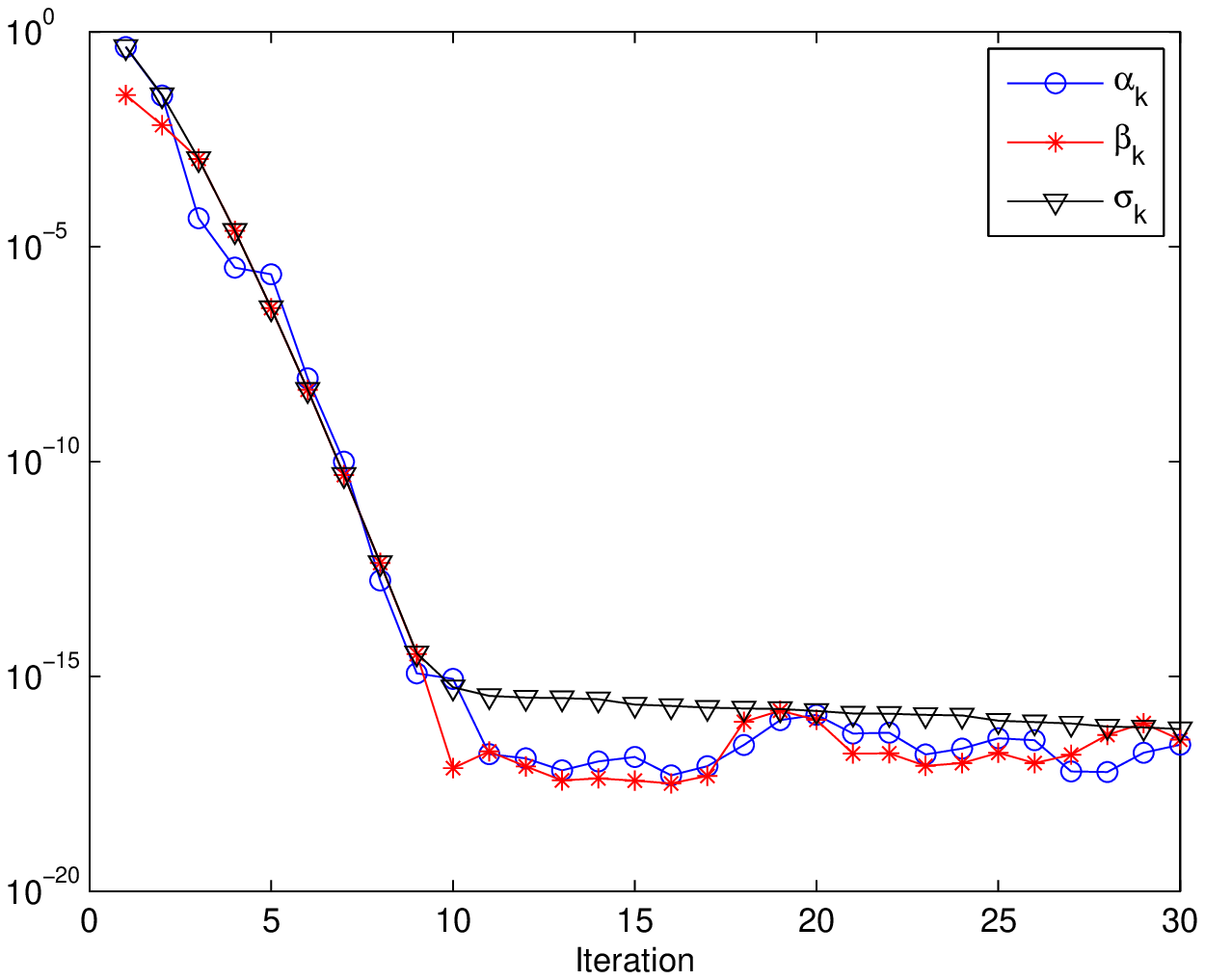}}
  \centerline{(b)}
\end{minipage}
\vfill
\begin{minipage}{0.48\linewidth}
  \centerline{\includegraphics[width=6.0cm,height=4.5cm]{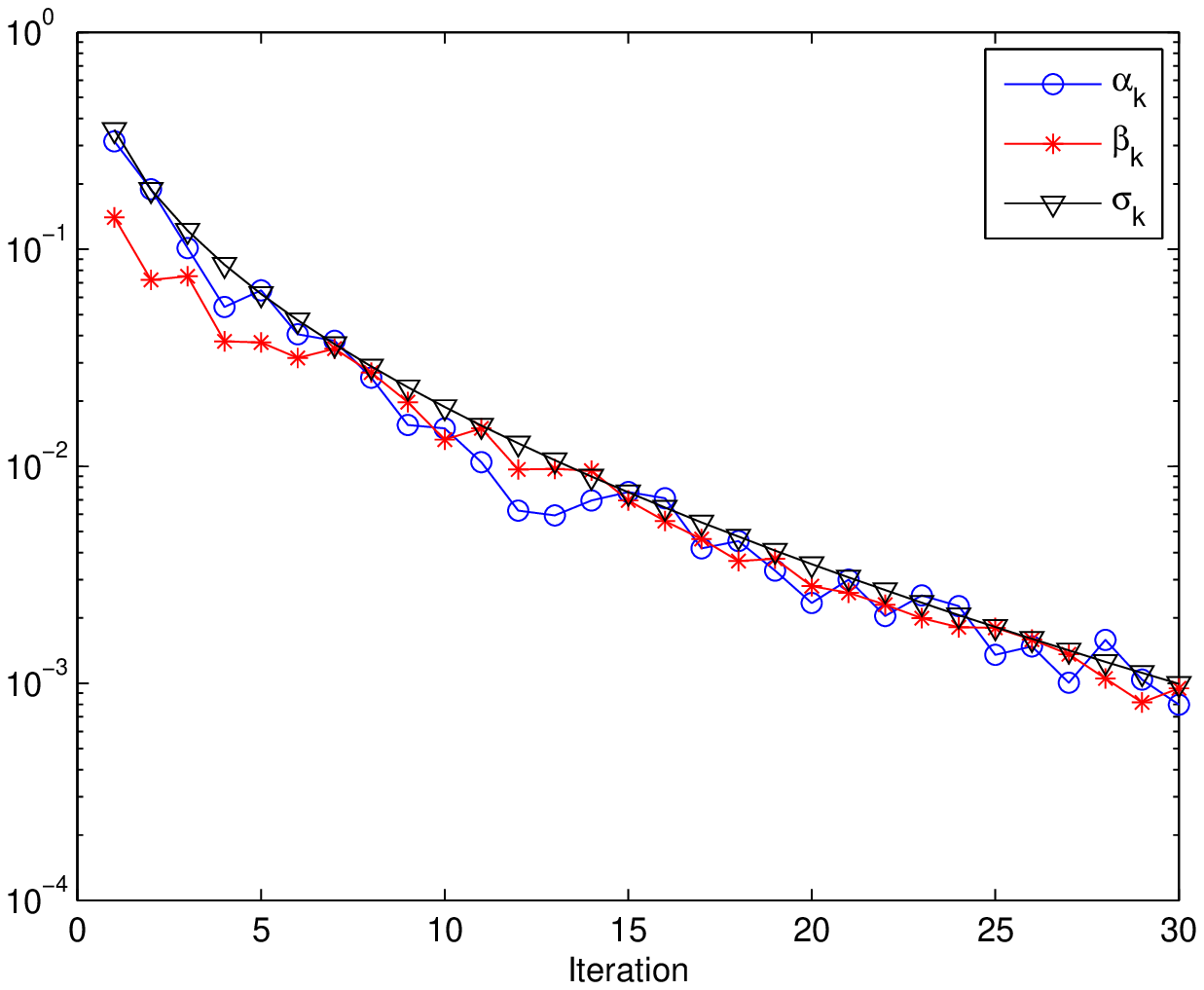}}
  \centerline{(c)}
\end{minipage}
\hfill
\begin{minipage}{0.48\linewidth}
  \centerline{\includegraphics[width=6.0cm,height=4.5cm]{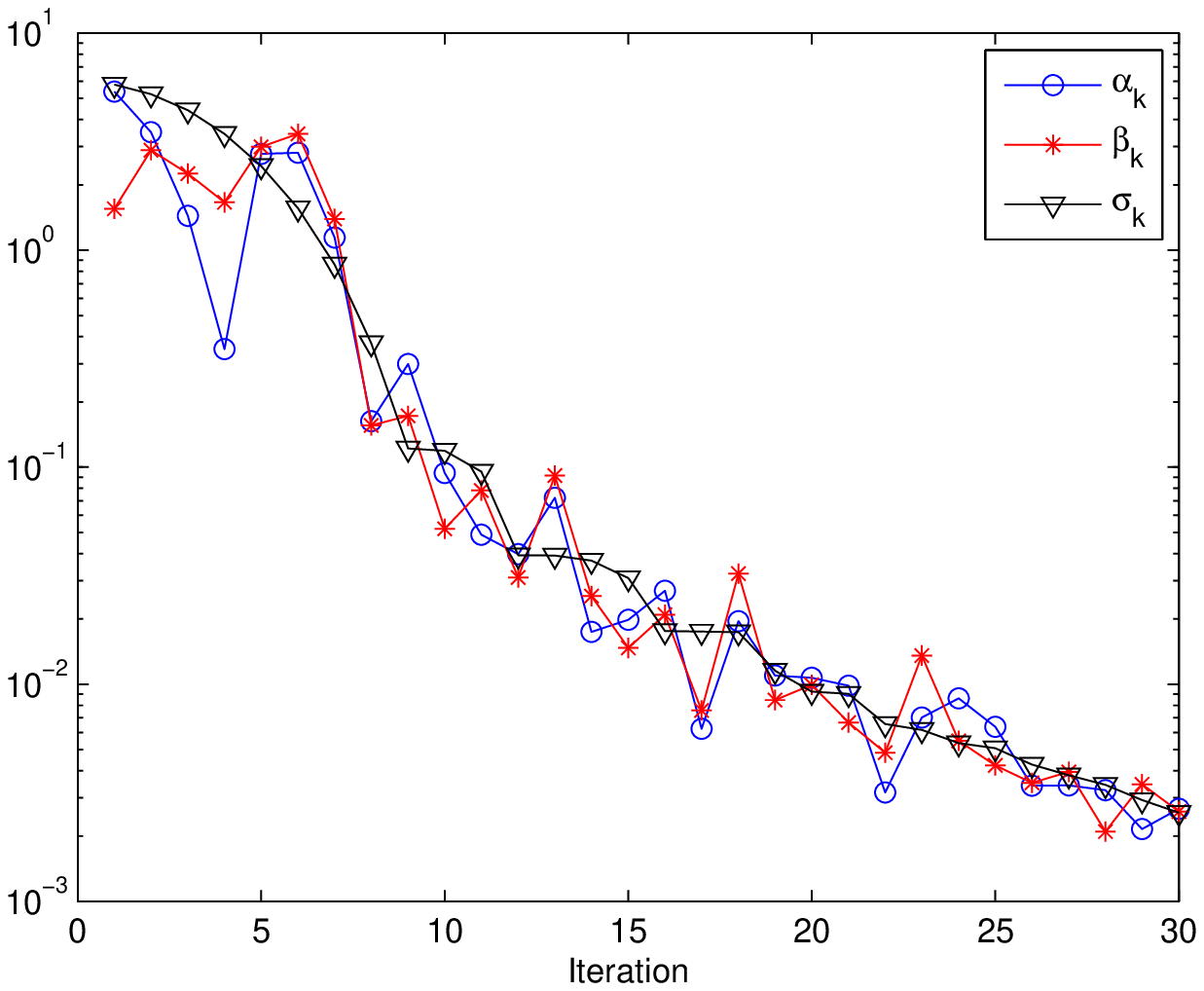}}
  \centerline{(d)}
\end{minipage}
\caption{(a)-(d): Decay curves of the sequences $\alpha_k$,
$\beta_{k+1}$ and $\sigma_k$ for $\mathsf{shaw, wing, i\_laplace}$ and $\mathsf{heat}$
(from top left to bottom right).} \label{fig7}
\end{figure}

Figure~\ref{fig7} illustrates that both $\alpha_k$ and $\beta_{k+1}$
decay as fast as $\sigma_k$, and for $\mathsf{shaw}$ and
$\mathsf{wing}$ all of them decay swiftly and level off at
$\epsilon_{\rm mach}$ due to round-off errors
in finite precision arithmetic. Precisely, they reach
the level of $\epsilon_{\rm mach}$ at $k=22$ and $k=8$ for $\mathsf{shaw}$ and
$\mathsf{wing}$, respectively. Such decay behavior has also
been observed in \cite{bazan14,gazzola14,gazzola-online}, but no theoretical
support was given. These experiments confirm Theorem~\ref{main1} and
Theorem~\ref{main2}, which have proved that $\gamma_k$ decreases as fast as
$\sigma_{k+1}$ and that $\alpha_k$, $\beta_{k+1}$ and $\alpha_k+\beta_{k+1}$
decay as fast as $\sigma_k$.

\subsection{A comparison of LSQR and the TSVD method}

We compare the performance of LSQR and the TSVD method
for the severely ill-posed $\mathsf{shaw,\ wing}$, the moderately ill-posed
$\mathsf{heat,\ phillips}$ and the mildly ill-posed problem
$\mathsf{deriv2}$ of $n=3,000$. We take $\varepsilon=10^{-3}$. For each problem,
we compute the relative errors of regularized solutions and
the residual norms obtained by the two methods. We will
demonstrate that LSQR has the full regularization for the
severely and moderately ill-posed problems, but it has only
the partial regularization for the mildly ill-posed problem. The
results on $\varepsilon=10^{-2},\ 10^{-4}$ are very similar, and we
thus omit them.

\begin{figure}
\begin{minipage}{0.48\linewidth}
  \centerline{\includegraphics[width=6.0cm,height=4.5cm]{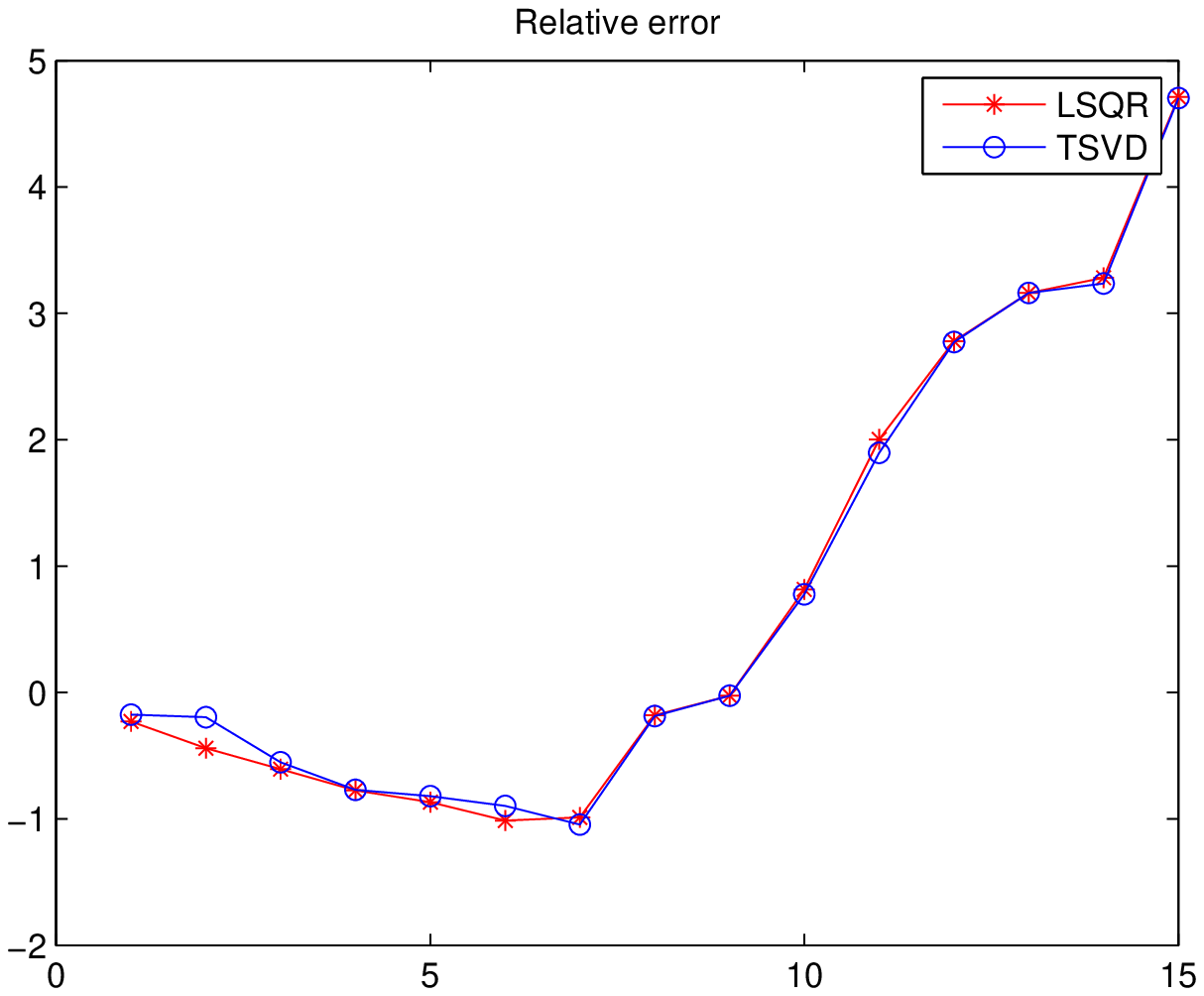}}
  \centerline{(a)}
\end{minipage}
\hfill
\begin{minipage}{0.48\linewidth}
  \centerline{\includegraphics[width=6.0cm,height=4.5cm]{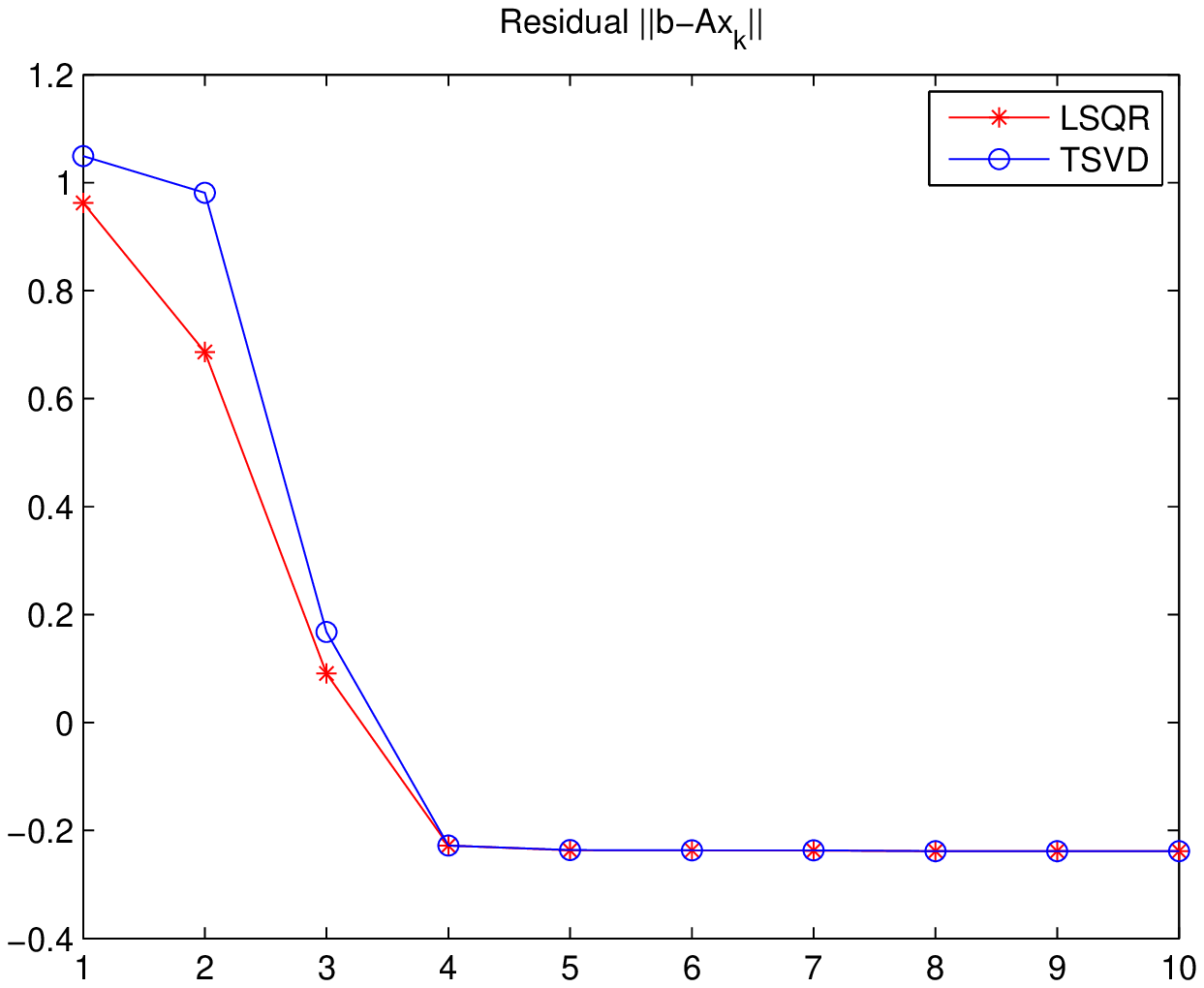}}
  \centerline{(b)}
\end{minipage}
\caption{Results for the severely ill-posed problem $\mathsf{shaw}$.}
\label{lsqrtsvd1}
\end{figure}

\begin{figure}
\begin{minipage}{0.48\linewidth}
  \centerline{\includegraphics[width=6.0cm,height=4.5cm]{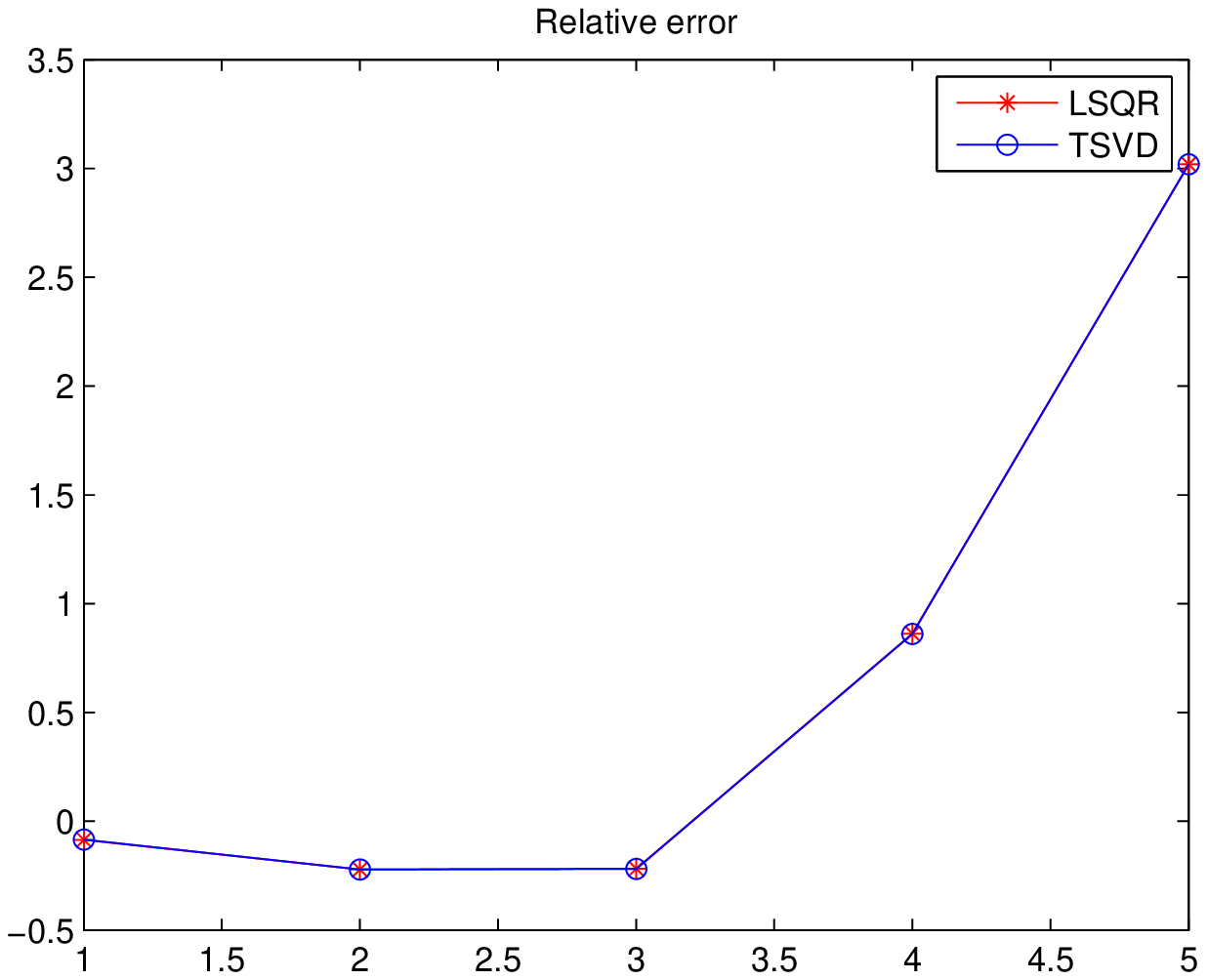}}
  \centerline{(a)}
\end{minipage}
\hfill
\begin{minipage}{0.48\linewidth}
  \centerline{\includegraphics[width=6.0cm,height=4.5cm]{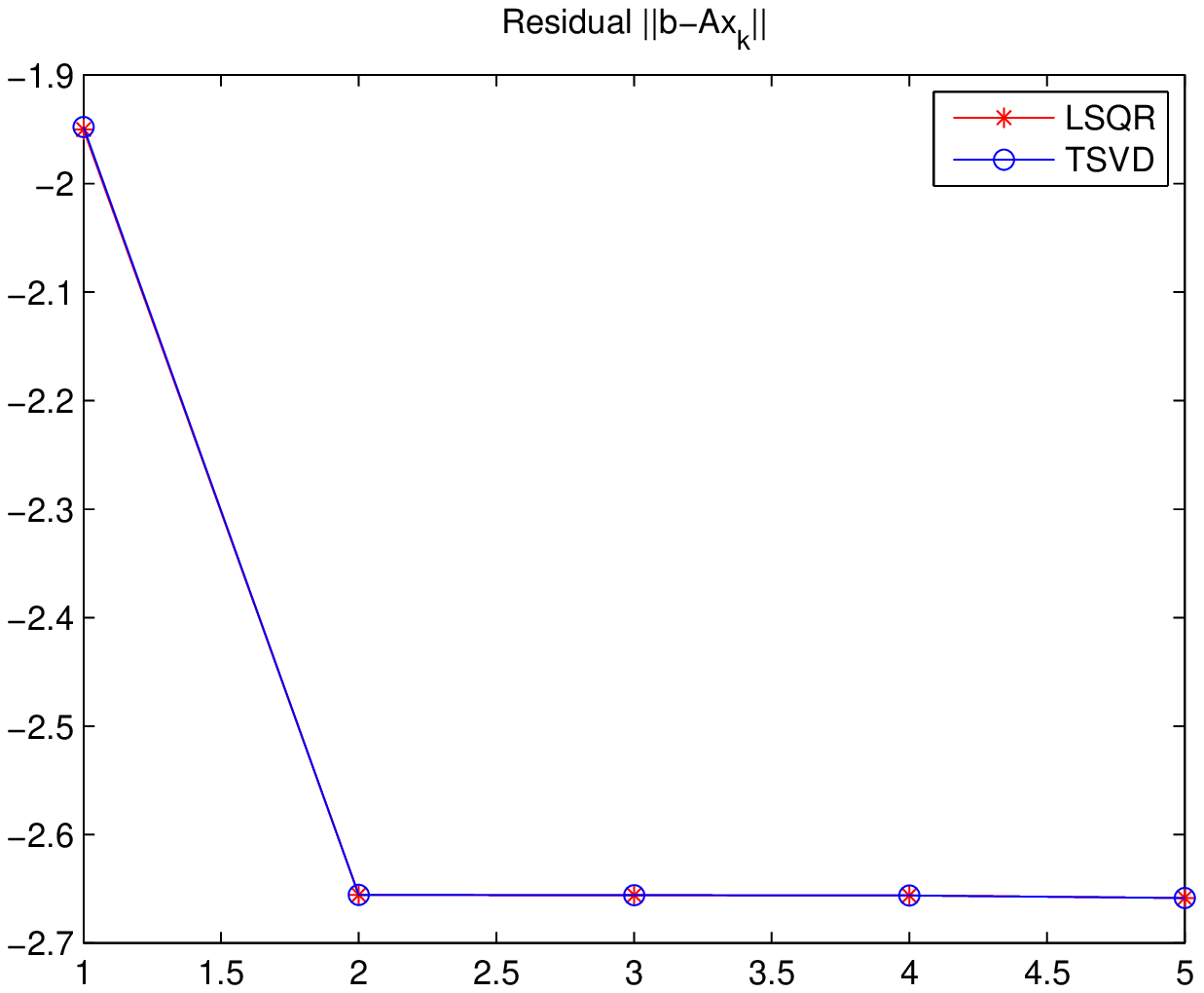}}
  \centerline{(b)}
\end{minipage}
\caption{Results for the severely ill-posed problem $\mathsf{wing}$.}
\label{lsqrtsvd2}
\end{figure}

Figures~\ref{lsqrtsvd1}--\ref{lsqrtsvd2} indicate LSQR and the TSVD method
behave very similarly for $\mathsf{shaw}$ and $\mathsf{wing}$. They
illustrate that, for $\mathsf{wing}$, the norms of approximate
solutions and the relative errors by the two methods
are almost indistinguishable for the same $k$, and, for $\mathsf{shaw}$, the
residual norms by LSQR decreases more quickly than the ones by the TSVD method
for $k=1,2,3$ and then they become almost identical starting from $k=4$.
These results demonstrate that LSQR has the full regularization.

\begin{figure}
\begin{minipage}{0.48\linewidth}
  \centerline{\includegraphics[width=6.0cm,height=4.5cm]{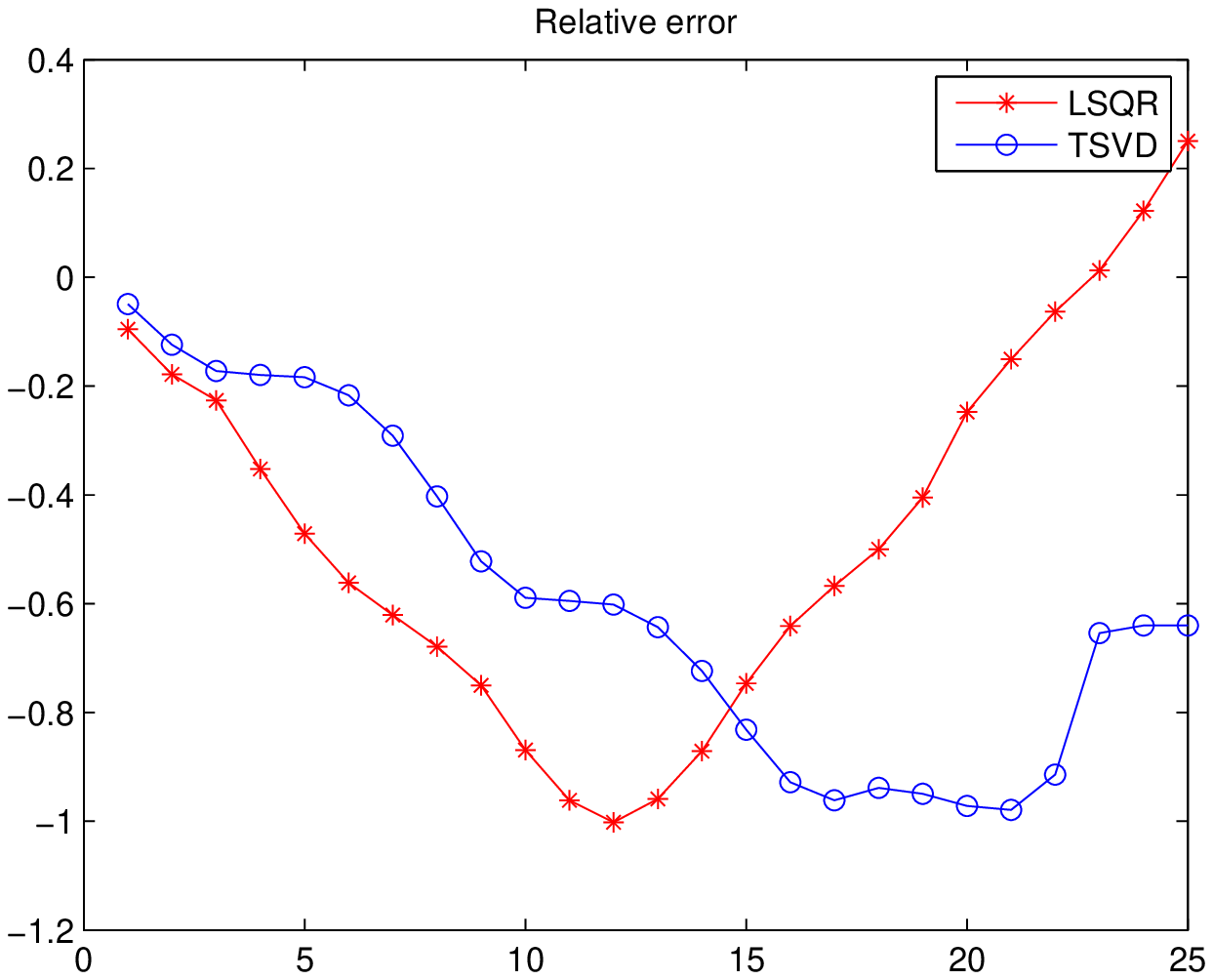}}
  \centerline{(a)}
\end{minipage}
\hfill
\begin{minipage}{0.48\linewidth}
  \centerline{\includegraphics[width=6.0cm,height=4.5cm]{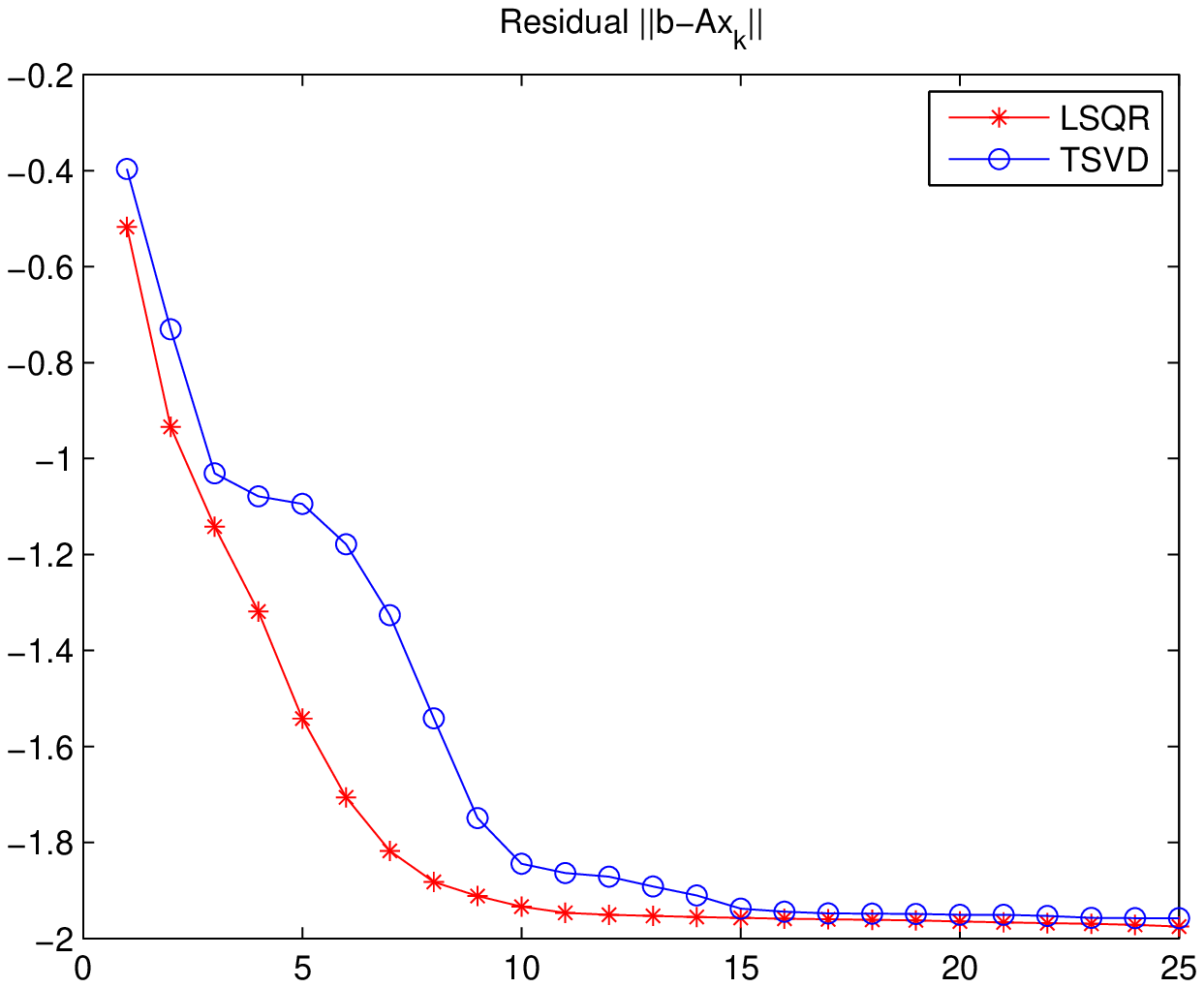}}
  \centerline{(b)}
\end{minipage}
\caption{Results for the moderately ill-posed problem $\mathsf{heat}$.}
\label{lsqrtsvd3}
\end{figure}

\begin{figure}
\begin{minipage}{0.48\linewidth}
  \centerline{\includegraphics[width=6.0cm,height=4.5cm]{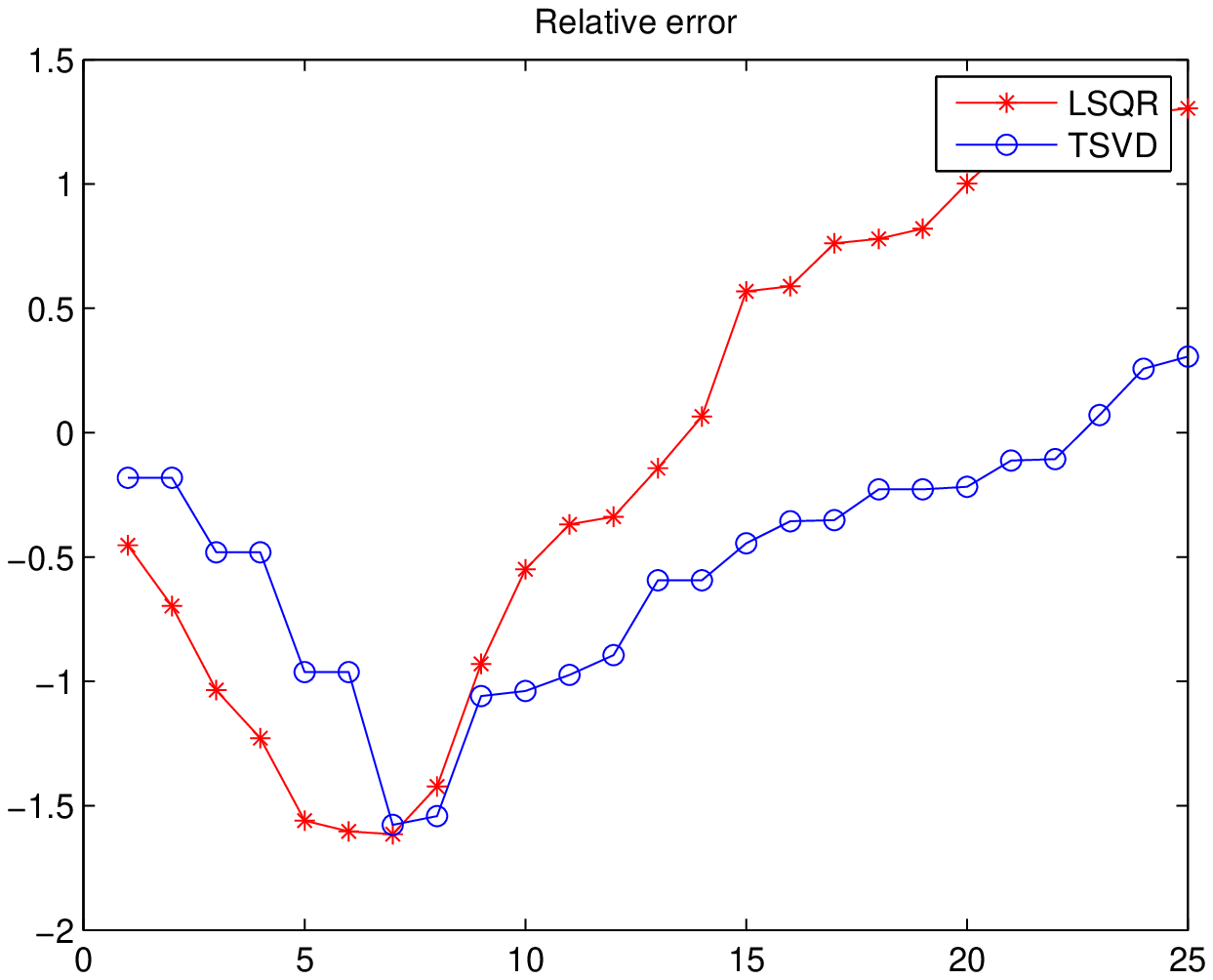}}
  \centerline{(a)}
\end{minipage}
\hfill
\begin{minipage}{0.48\linewidth}
  \centerline{\includegraphics[width=6.0cm,height=4.5cm]{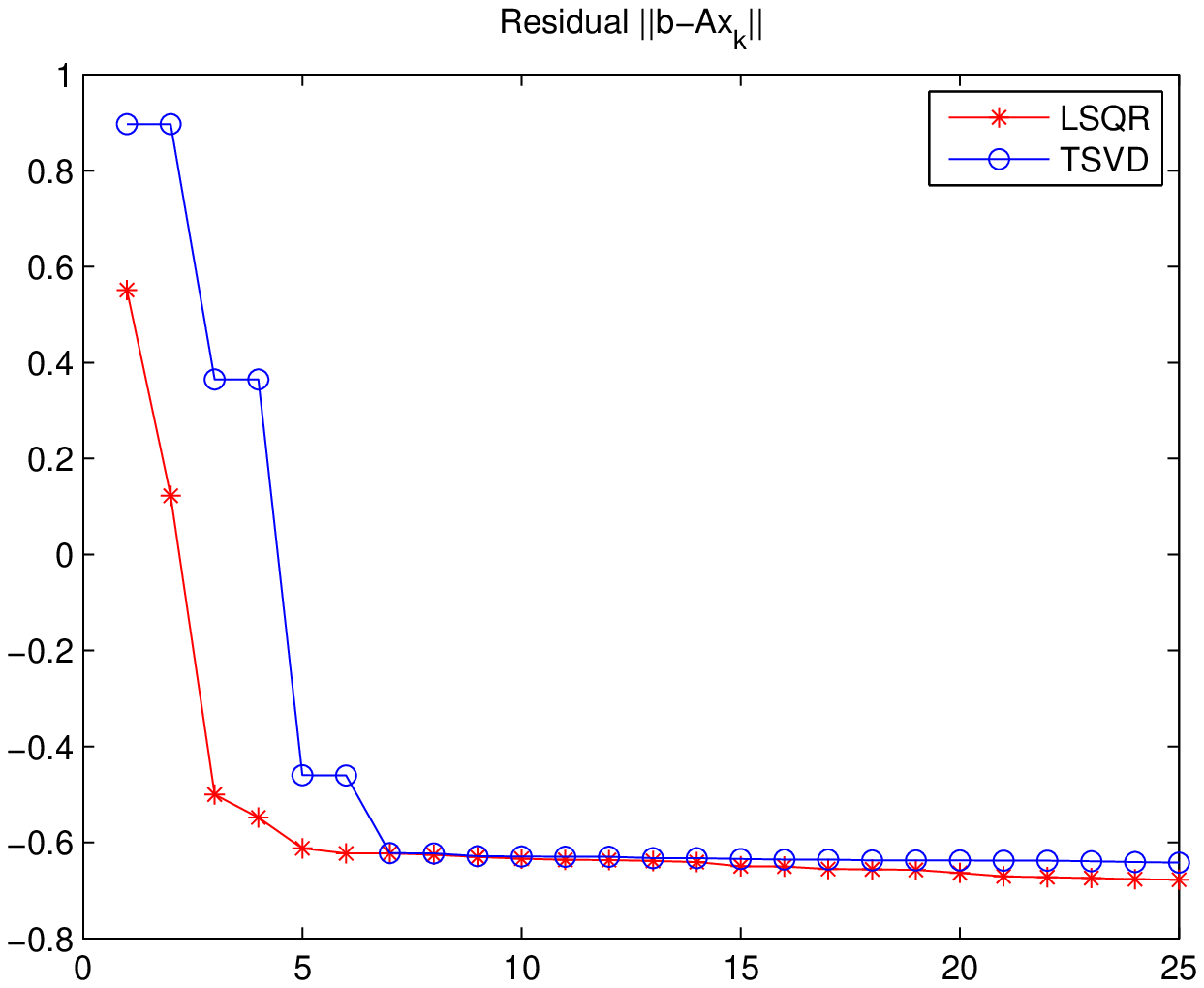}}
  \centerline{(b)}
\end{minipage}
\caption{Results for the moderately ill-posed problem $\mathsf{phillips}$.}
\label{lsqrtsvd4}
\end{figure}

For each of $\mathsf{heat}$ and $\mathsf{phillips}$,
Figures~\ref{lsqrtsvd3}--\ref{lsqrtsvd4} demonstrate that
the best regularized solution obtained by LSQR is at
least as accurate as, in fact, a little bit more accurate than
that by the TSVD method, and the corresponding
residual norms decreases and drop below at least the same level as those
by the TSVD method. The residual norms by the two methods then stagnate after
the best regularized solutions are found. All these confirm that
LSQR has the full regularization.

\begin{figure}
\begin{minipage}{0.48\linewidth}
  \centerline{\includegraphics[width=6.0cm,height=4.5cm]{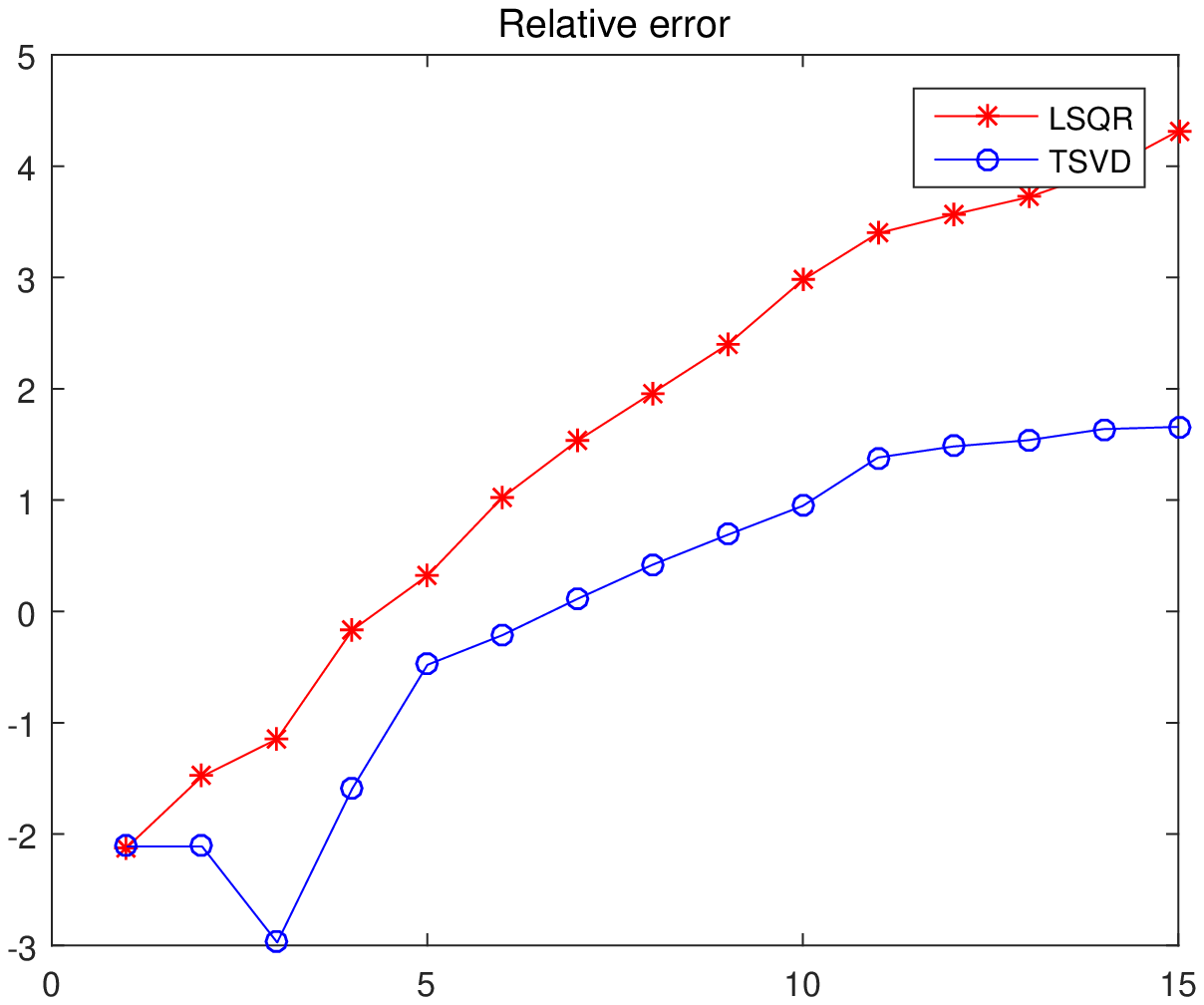}}
  \centerline{(a)}
\end{minipage}
\hfill
\begin{minipage}{0.48\linewidth}
  \centerline{\includegraphics[width=6.0cm,height=4.5cm]{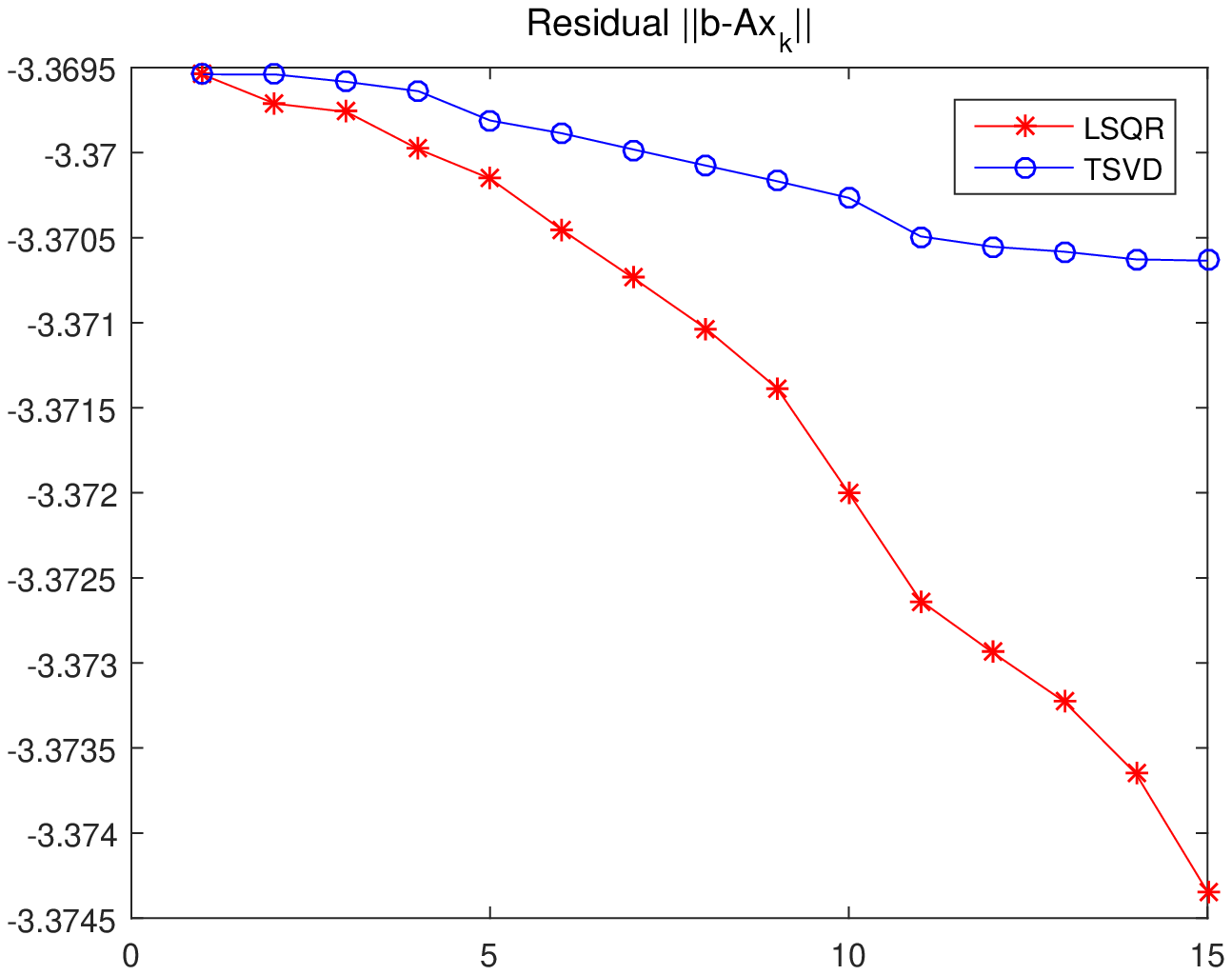}}
  \centerline{(b)}
\end{minipage}
\caption{Results for the moderately ill-posed problem $\mathsf{phillips}$.}
\label{lsqrtsvd5}
\end{figure}

To better illustrate the regularizing effects of LSQR, we test a
larger $\mathsf{deriv2}$ of $n=3000$ whose condition
number is $1.1\times 10^{7}$. Figure~\ref{lsqrtsvd5} demonstrates
that the best regularized solution by
LSQR at semi-convergence is considerably less accurate than $x_{k_0}^{tsvd}$.
Actually, the relative error of the former is $8.0\times 10^{-3}$, while that
of the latter is only $1.1\times 10^{-3}$, almost one order more accurate.
As we have observed, the semi-convergence of LSQR occurs at the very first
iteration, while the best regularized solution $x_{k_0}^{tsvd}$ consists
of three dominant SVD components of $A$.
The results clearly shows that LSQR has only the partial regularization for
mildly ill-posed problems.

From the figures we observe some obvious differences
between moderately and severely ill-posed problems. For $\mathsf{heat}$,
it is seen that the relative errors and
residual norms converge considerably more quickly for the LSQR solutions
than for the TSVD solutions. Figure~\ref{lsqrtsvd3} (a) tells us that
LSQR only uses 12 iterations to find the best regularized
solution, but the TSVD method finds the best regularized
solution for $k_0=21$.
Similar differences are observed for
$\mathsf{phillips}$, where Figure~\ref{lsqrtsvd4} (a) indicates
that both LSQR and the TSVD method find the
best regularized solutions at $k_0=7$.

We can observe more. Figure~\ref{lsqrtsvd3} shows that
the TSVD solutions improve little and their residual norms decrease
very slowly for the indices $i=4,5,11,12,18,19,20$.  This implies that
the $v_i$ corresponding to these indices $i$ make very little contribution
to the TSVD solutions. This is due to the fact that the Fourier coefficients
$|u_i^T\hat{b}|$ are very small relative to $\sigma_i$ for these indices $i$.
Note that $\mathcal{K}_k(A^TA,A^Tb)$ adapts itself in
an optimal way to the specific right-hand side $b$, while the TSVD method uses
all $v_1,v_2,\ldots,v_k$ to construct a regularized solution,
independent of $b$. Therefore, $\mathcal{K}_k(A^TA,A^Tb)$ picks up only those
SVD components making major contributions to $x_{true}$, such that LSQR uses
possibly fewer $k$ iterations than $k_0$ needed by the TSVD method to capture
those truly needed
dominant SVD components. The fact that LSQR (CGLS) includes fewer SVD components
than the TSVD solution with almost the same accuracy was first noticed
by Hanke \cite{hanke01}.
Generally, for severely and moderately ill-posed problems,
we may deduce that LSQR uses possibly fewer than $k_0$ iterations to compute a
best possible regularized solution if, in practice, some of
$|u_i^T b|$, $i=1,2,\ldots,k_0$ are considerably bigger than
the corresponding $\sigma_i$ and some of them are reverse.
For $\mathsf{phillips}$,
as noted by Hansen \cite[p.32, 123--125]{hansen10}, half of the SVD
components satisfy $u_i^T \hat{b}=v_i^Tx_{true}=0$ for $i$ even, only the
odd indexed $v_1,v_3,\ldots,$ make contributions to $x_{true}$. This is why
the relative errors and residual norms of TSVD solutions
do not decrease at even indices before $x_{k_0}^{tsvd}$ is found.

\section{Conclusions}\label{concl}

For the large-scale \eqref{eq1}, iterative solvers
are the only viable approaches. Of them, LSQR and CGLS are most popularly
used for general purposes, and CGME and LSMR are also choices.
They have general regularizing effects and exhibit semi-convergence.
However, if semi-convergence occurs before it
captures all the needed dominant SVD components, then best possible
regularized solutions are not yet found and the solvers have only the
partial regularization. In this case, their hybrid variants have often
been used to compute best possible regularized solutions.
If semi-convergence means that they have already found
best possible regularized solutions, they have the full regularization,
and we simply stop them after semi-convergence.

For the case that the singular values of $A$ are all simple,
we have considered the fundamental open question in depth: Do LSQR and CGLS
have the full or partial regularization for severely, moderately
and mildly ill-posed problems? We have first considered the case that
all the singular values of $A$ are simple. As a key and indispensable step, we
have established accurate bounds for the 2-norm
distances between the underlying $k$ dimensional Krylov subspace and the
$k$ dimensional dominant right singular subspace for the three kinds of
ill-posed problems under consideration. Then we have
provided other absolutely necessary background and ingredients. Based on them,
we have proved that, for severely or moderately ill-posed problems
with $\rho>1$ or $\alpha>1$ suitably, LSQR has the full regularization.
Precisely, for $k\leq k_0$
we have proved that a $k$-step Lanczos bidiagonalization produces a near best
rank $k$ approximation of $A$ and the $k$ Ritz values approximate
the first $k$ large singular values of $A$ in natural order, and
no small Ritz value smaller than $\sigma_{k_0+1}$ appears before
a best possible regularized solution has been found.
For mildly ill-posed problems, we have proved
that LSQR generally has only the partial regularization
since a small Ritz value generally appears before all
the needed dominant SVD components are captured. Since CGLS is mathematically
equivalent to LSQR, our assertions on the full or partial regularization
of LSQR apply to CGLS as well.

We have derived bounds for the diagonals and subdiagonals of bidiagonal
matrices generated by Lanczos bidiagonalization. Particularly,
we have proved that they decay as fast as the singular values of $A$
for severely ill-posed problems or moderately ill-posed problems
with $\rho>1$ or $\alpha>1$ suitably and decay more slowly
than the singular values of $A$ for mildly ill-posed problems.
These bounds are of theoretical and practical importance, and they
can be used to identify the degree of ill-posedness without extra cost
and decide the full or partial regularization of LSQR. We have made detailed
and illuminating numerical experiments, confirming our
theory.

Our analysis approach can be adapted to MR-II for symmetric ill-posed problems,
and certain definitive assertions are expected for
three kinds of symmetric ill-posed problems.
Our approach are applicable
to the preconditioned CGLS (PCGLS) and LSQR (PLSQR) \cite{hansen98,hansen10}
by exploiting the transformation technique originally proposed
in \cite{bjorck79} and advocated in \cite{hanke92,hanke93,hansen07}
or the preconditioned MR-II \cite{hansen10,hansen06},
all of which correspond to a general-form Tikhonov regularization involving the
matrix pair $\{A,L\}$, in which the regularization term $\|x\|^2$ is replaced by
$\|Lx\|^2$ with some $p\times n$ matrix $L\not=I$. It
should also be applicable to the mathematically equivalent LSQR
variant \cite{kilmer07} that is based on a joint bidiagonalization of
the matrix pair $\{A,L\}$ that corresponds to the above general-form Tikhonov
regularization. In this setting, the Generalized SVD (GSVD)
of $\{A,L\}$ or the mathematically equivalent SVD of $AL_A^{\dagger}$
will replace the SVD of $A$ to play a central role in analysis, where
$L_A^{\dagger}=\left(I-\left(A(I-L^{\dagger}L)^{\dagger}A\right)\right)^{\dagger}
L^{\dagger}$ is call the {\em $A$-weighted generalized inverse of $L$}
and $L_A^{\dagger}=L^{-1}$ if $L$ is square and invertible;
see \cite[p.38-40,137-38]{hansen98} and \cite[p.177-183]{hansen10}.

\begin{acknowledgements}
I thank Dr. Yi Huang and Mrs. Yanfei Yang for running the numerical
experiments. I am grateful to Professors \AA. Bj\"{o}rck,
P. C. Hansen, L. Reichel and
D.~P. O'Leary for their comments and suggestions that helped
improve the presentation of this paper.
\end{acknowledgements}

\end{document}